\documentclass[a4paper,reqno]{amsart}

\usepackage{amsmath}
\usepackage{paralist}
\usepackage{graphics}
\usepackage{epsfig}
\usepackage{amsfonts}
\usepackage{amssymb}
\usepackage{hyperref}
\usepackage{epstopdf}
\usepackage{color}

\newtheorem{theorem}{Theorem}[section]
\newtheorem{lemma}[theorem]{Lemma}

\def\ifl{\iffalse }

\def\bc{\begin{center}}       \def\ec{\end{center}}
\def\ba{\begin{array}}        \def\ea{\end{array}}
\def\be{\begin{equation}}     \def\ee{\end{equation}}
\def\bea{\begin{eqnarray}}    \def\eea{\end{eqnarray}}
\def\beaa{\begin{eqnarray*}}  \def\eeaa{\end{eqnarray*}}

\numberwithin{equation}{section}

\newtheorem{definition}[theorem]{Definition}

\newtheorem{remark}[theorem]{Remark}
\numberwithin{equation}{section}

\begin{document}
	
\title[Chemorepulsive-haptotaxis,  regularity, exponential convergence]
{Global solvability, regularity  and exponential convergence in a quasilinear chemotaxis-haptotaxis system with mostly repulsive signal production}


\author{Yuanlin Chen}
\address{School of Mathematics and Institute for Mathematical Sciences, Renmin University of China, Beijing, 100872, China}
\email{chenyuanlin@ruc.edu.cn}

\author{Tian Xiang$^*$}
\address{School of Mathematics and Institute for Mathematical Sciences, Renmin University of China, Beijing, 100872, China}
\email{txiang@ruc.edu.cn}
\thanks{$^*$ Corresponding author.}





\begin{abstract}
	In this work, motivated from biological and mathematical viewpoints, we systematically investigate the following quasilinear chemotaxis-haptotaxis system with mostly repulsive signal production:
	\be
	\begin{cases}
		u_t =  \Delta u+\chi \nabla\cdot(u\nabla v)-\xi\nabla\cdot(S(u)\nabla w)  &\text{in } \Omega\times(0,\infty), \\[0.2cm]
		\tau v_t = \Delta v-v+u  &\text{in } \Omega\times(0,\infty),\\[0.2cm]
		w_t=-vw 	 &\text{in } \Omega\times(0,\infty)
	\end{cases}
	\nonumber
	\ee
	in a bounded and smooth domain $\Omega\subset\mathbb{R}^n(n\geq2)$ with nonnegative initial data and no-flux  boundary conditions, where $\tau\in\{0,1\}$, the given sensitivity parameters $\chi,  \xi$ are nonnegative, and the $C^2$-smooth nonnegative hapto-sensitivity function $S(u)$ generalizes the prototype  $S(u)=u(u+1)^{\theta-1}$ with $\theta\in\mathbb{R}$.
	
	First, for $\theta<\frac{3}{2}$,  we establish  uniform-in-time boundedness and exponential convergence of classical solutions to its constant equilibrium  for either  $n=2$ or $\tau=0$. Then for  $\theta\leq1$,  we establish   uniform boundedness and exponential convergence of classical solutions under suitable smallness conditions on initial data in optimal Lebesgue spaces irrespective of the sign of $\chi$. Finally, for the physical dimension  $n=3$, we first show  global existence of weak solutions for $\theta<\frac{7}{6}$,  and then, we  show  intermediate and eventual regularity and exponential convergence of the derived global weak solutions for $\theta\leq\frac{1}{2}$. These findings greatly extend the well-known  knowledge on the corresponding   Keller-Segel repulsive chemotaxis-only model upon setting $w=0$ to the repulsive chemotaxis-haptotaxis model, showing positive effect of chemo-repulsion on global dynamics in such  chemotaxis-haptotaxis settings.
\end{abstract}

\maketitle

\section{Introduction and statement of  main results}

	Chemotaxis, the directed movement of cells in response to a concentration gradient of  chemical stimulus in their environment, is an integral part of many biological processes - for instance, bacteria monition, cellular organism development, immune response, etc. Thus, mathematical modeling of chemotaxis has received great interests  in  many fields. The following well-known Neumann-initial boundary value problem (IBVP) is the most famous Keller-Segel minimal chemotaxis model \cite{Keller-Segel70}:
	\be\label{KS-model}
	\begin{cases}
		u_t =  \Delta u-\chi\nabla\cdot\left(u\nabla v\right)  &\text{in } \Omega\times(0,\infty), \\[0.2cm]
		\tau v_t = \Delta v-v+u  &\text{in } \Omega\times(0,\infty),\\[0.2cm]
		\frac{\partial u}{\partial\nu}=\frac{\partial v}{\partial\nu}=0  &\text{on } \partial\Omega\times(0,\infty), \\[0.2cm]
		(u,\tau v)=(u_0,\tau v_0)  &\text{on } \Omega\times\{0\},
	\end{cases}
	\ee
	where $\tau\in\{0,1\}$,  $\chi>0$, and the unknown functions $u$ and $v$ denote respectively the cell population density and chemical concentration, and $\frac{\partial}{\partial\nu}$ represents the outward normal derivative with respect to  $\partial\Omega$. A striking feature of this classical  chemotaxis model \eqref{KS-model} is the cell aggregation  phenomenon, which is exhibited by a critical mass of finite/infinite time blow-up of solutions in 2D (cf. \cite{HW01, Nagai95, Nagai-Senba-Yoshida97}) and a generic mass of finite time blow-up of solutions in 3D or higher  (cf. \cite{Bellomo-Bellouquid-Tao-Winkler15, Winkler10JDE, Winkler13}).
	
	The chemotaxis system \eqref{KS-model} indicates that cells move toward  higher concentration of chemical stimuli, which is called chemo-attraction (or positive chemotaxis). Nevertheless, if  chemical stimuli are harmful/dangerous, cells typically escape from the unfavorable environment and move in the opposite direction.  For example, E. coli cells swim away from potentially noxious chemicals, such as alcohols and fatty acids; preys escapes from predators, etc.   This process is called chemo-repulsion (or negative chemotaxis). Simply switching  $-\chi$ into $+\chi$ in \eqref{KS-model}   turns into the following chemo-repulsive Keller-Segel model:
	\be\label{KS-model-repulsion}
	\begin{cases}
		u_t =  \Delta u+\chi\nabla\cdot\left(u\nabla v\right)  &\text{in } \Omega\times(0,\infty), \\[0.2cm]
		\tau v_t = \Delta v-v+u  &\text{in } \Omega\times(0,\infty).
	\end{cases}
	\ee
	This change of sign  in \eqref{KS-model} drastically  alters the dynamics of  \eqref{KS-model}. More specifically, in the parabolic-elliptic case, i.e., $\tau=0$, the chemo-repulsion system \eqref{KS-model-repulsion} possesses globally bounded solutions, which decays exponentially to  constant steady states, revealing that negative chemotaxis effect enhances global boundedness and thus prevent  blow-up \cite{Mock74, Mock75}. In fully parabolic case, i.e. $\tau=1$, the solutions exist globally in time and converge exponentially to the steady states in 2D, and global weak solutions and their convergence in a certain weak sense are obtained when $n=3,4$ \cite{Cieslak-Laurencot-MM08}. In 2019, Jiang \cite{Jiang19} derives  eventual smoothness and exponential stabilization of such global weak solutions in 3D.   When $\Delta u$ is replaced by  $\epsilon\Delta u$ with $\epsilon>0$ being sufficiently small \cite{CF19} or when $\Omega\subset\mathbb{R}^3$ is convex and assuming $\nabla \sqrt{u}\in L^2_{\text{loc}}(0,\infty; L^2(\Omega))$ \cite{CFHS24},  the global existence (and boundedness in the former case) of classical solutions is known.  However,  for the chemo-repulsion model \eqref{KS-model-repulsion} with $\tau=1$ and $n\geq 3$,  without any further conditions, the global existence (not to  mention boundedness) of classical solutions remains widely open, though it is expected. Further illustrations of boundedness/regularity-enhancement effect of chemo-repulsion in more complex chemo-repulsion models with possibly nonlinear diffusion, nonlinear (singular) sensitivity or nonlinear signal production can be found in  \cite{Chen-Kong-Wang21, Freit-18, Jiang22, Lin-Xiang21, Tao13}.

	On the other hand, let's recall another important type of directed movement, called  haptotaxis,  a type of cell motility in which cells move in a directed manner toward the gradient of non-diffusible extracellular adhesion sites or substrate-bound chemoattractants or repellents. Combining chemotaxis- and haptotaxis-movement, Chaplain and Lalas \cite{Chaplain-Lolas05,Chaplain-Lolas06} extend the KS chemotaxis  model to  describe the cancer cell invasion and metastatic spread.   In that procedure, alongside random motion, tumor cells (with density by $u$) are presumed to align their motion through two factors: firstly, a specific diffusible chemical substance (with density by $v$) known as matrix degrading enzymes (MDEs) secreted by tumor cells (i.e., chemotaxis); and secondly, the non-diffusive extracellular matrix (ECM) (with density by $w$) decomposed by the MDEs (i.e., haptotaxis). In such a way,  a classical Chaplain-Lolas model fulfills the following chemotaxis-haptotaxis model:
	\be\label{Chaplain-Lolas model}
	\begin{cases}
		u_t =  \Delta u-\chi\nabla\cdot\left(u\nabla v\right)-\xi\nabla\cdot\left(u\nabla w\right)+\mu u(1-u-w)  &\text{in } \Omega\times(0,\infty), \\[0.2cm]
		\tau v_t = \Delta v-v+u  &\text{in } \Omega\times(0,\infty), \\[0.2cm]
		w_t = -vw &\text{in } \Omega\times(0,\infty),
	\end{cases}
	\ee
where the new parameters $\mu,\xi\geq0$. This model contains  an interplay of   diffusion, chemotactic cross-diffusion, haptotactic cross-diffusion, and logistic damping. A lot of studies are concerned about the global existence, boundedness and large time behavior of solutions to the corresponding  IBVP  \eqref{Chaplain-Lolas model} and its various variants. Firstly, without haptotaxis, i.e.,  $w\equiv0$, the model \eqref{Chaplain-Lolas model} reduces to a chemotaxis-only system with logistic source, in which blow-up is ruled out  for $n\leq2$ \cite{Osaki-Tsujikawa-Yagi-Mimura02, Xiang18JMP}. Suitably strong logistic damping ensures  boundedness and convergence for  $n\geq3$  \cite{Winkler10CPDE, Xiang18JMAA}. While, weak logistic damping may induce blow-up  for $n\geq 5$ \cite{Fuest21}. Secondly, with the presence of haptotaxis, existing literature still shows that the logistic damping role is  significant: roughly speaking,  any presence of (even sub-) logistic damping is sufficient to prevent blow-up in $\leq$2D \cite{Tao14JO, Tao-Wang08, Tao-Wang09, Xiang-Zheng19Non} and properly strong logistic damping prevents blow-up in $\geq$ 3D \cite{Cao16, Li-Mu-Yi19, Tao-Wang09} and further strong logistic damping ensures convergence to constant equilibrium \cite{Hillen-Painter-Winkler13, Tao-Winkler14PRSE, Tao-Winkler14Non, Tao-Winkler15,Zheng-Ke19}.  Finally, we observe there exists  vast literature investigating  dynamical properties of solutions to a general framework of \eqref{Chaplain-Lolas model} with more complex mechanisms like nonlinear diffusion, nonlinear chemo-sensitivity,
remodeling effects and generalized logistic source, etc, cf.  \cite{Dai-Liu20, Jin-ch23, LL16-non, YFW20} and the references therein for more details.

In the presence of logistic damping, upon the observation that haptotaxis effect seems to be overbalanced by  chemotaxis effect,  several recent studies \cite{Jin-Xiang21, Liu-Li23, Xiang22} have been focusing on the competitions between chemotactic cross-diffusion and  haptotactic  cross-diffusion in the  model  \eqref{Chaplain-Lolas model} without logistic source:
	\be\label{Chaplain-Lolas model+}
	\begin{cases}
		u_t =  \Delta u-\chi\nabla\cdot\left(u\nabla v\right)-\xi\nabla\cdot\left(u\nabla w\right)  &\text{in } \Omega\times(0,\infty), \\[0.2cm]
		\tau v_t = \Delta v-v+u  &\text{in } \Omega\times(0,\infty), \\[0.2cm]
		w_t = -vw  &\text{in } \Omega\times(0,\infty).
	\end{cases}
	\ee
  Indeed, even with re-modeling of ECM, i.e., the third equation in \eqref{Chaplain-Lolas model+} is replaced by $w_t=-vw+\eta w(1-w)$, the negligibility of haptotaxis vs. chemotaxis is comprehensively examined in 2D  \cite{Jin-Xiang21}: in particular, when  $\chi\|u_0\|_{L^1(\Omega)}<4\pi$,  it is shown that the corresponding IBVP \eqref{Chaplain-Lolas model+} admits only global-in-time solutions for all $\eta\geq0$ and remains uniformly bounded in time for small $\eta\geq0$ and that solution converges exponentially to that of chemotaxis-only system if $\chi>0$ is small; on the other hand, when $\chi\|u_0\|_{L^1(\Omega)}>4\pi$, the emanating solution may blow-up.  When $\Omega=\mathbb{R}^2$ and $\tau=1$, under $\chi\|u_0\|_{L^1(\Omega)}<8\pi$, a unique global strong solution to the Cauchy problem to \eqref{Chaplain-Lolas model+} is derived in \cite{Liu-Li23}. Recently, it is shown that pure haptotactic cross-diffusion in \eqref{Chaplain-Lolas model+}  cannot induce blow-up and pattern for either  $n\leq 3$ or $\tau=0$ and that generic blow-up occurs in radial setting when $\chi\|u_0\|_{L^1(B(R))}>A(n,R) \xi\|w_0\|_{L^\infty(B(R))}$ for $n\geq 3$ \cite{Xiang22}, where $B(R)\subset \mathbb{R}^n$ is the ball centered at the origin with radius $R>0$ and $A(n, R)>0$ is a constant depending only on $n$ and $R$. In summary, we find  that solutions may form singularity for large  $\chi\|u_0\|_{L^1(\Omega)}>0$  in 2D and for even small $\chi\|u_0\|_{L^1(\Omega)}>0$ in $\geq 3$D, which show  that the solution behaviors of   \eqref{Chaplain-Lolas model+} are more or less the same as that of  \eqref{KS-model}, especially in terms of global boundedness. Finally, we observe that  a quasilinear version of \eqref{Chaplain-Lolas model+} is very recently considered in \cite{RT24JDE}, where the first equation is replaced by
  $$
  u_t=\nabla\cdot(D(u)\nabla u)-\chi\nabla\cdot\left(S(u)\nabla v\right)-\xi\nabla\cdot\left(u\nabla w\right);
  $$
  boundedness of classical solutions is derived under enhanced diffusion and subcritical chemosesitivity, and small-data boundedness for $D(u)\equiv 1$ and $S(u)\equiv u$ is derived  under rough smallness on initial data compared to Theorem \ref{theorem3} below. When reducing to the case of $w\equiv0$, the derived results are worse than that of the corresponding chemotaxis-only model. In the supercritical quasilinear chemotaxis-haptotaxis model where the chemo-sensitivity and hapto-sensitivity are assumed to be the same as $S(u)$, the claimed blow-up phenomenon  \cite[Theorem 1.6]{RT24JDE} should be at most formulated as "almost" blow-up in the spirit of \cite[Theorem 1.1 (B3)]{Jin-Xiang21}. More recently, using essentially the same trick as  in \cite{RT24JDE} and adapting roughly from the blow-up
techniques for chemotaxis-only system in \cite{Winkler13} and chemotaxis-haptotaxis system in \cite{Jin-Xiang21, Xiang22}, Rani
and Tyagi \cite{RT25JDE} tried to extend our blow-up results for parabolic-elliptic model  in \cite{Xiang22} to fully parabolic case, while, their blow-up arguments have serious gaps, especially, the positivity of $\mathcal{D}$ in \cite[Lemma 3.6]{RT25JDE}.

 In 1974, Tso and Adler \cite{Tso-Adler74} investigated  chemo-repulsion in \textit{E. coli} cells  and proposed several approaches to measure and study this important negative chemotaxis mechanism, and they performed  successive  \textit{in vivo} experiments to suggest that no aggregation appears in chemo-repulsion systems. The above surveyed results  about \eqref{KS-model-repulsion} confirm  their experiment observation in 2D in the classical sense and higher dimensions in the weak solution sense.  On the other hand,  we have also reported  the negligibility of haptotaxis in the chemotaxis-haptotaxis model \eqref{Chaplain-Lolas model+} compared to the chemotaxis-only model \eqref{KS-model}. Accordingly, we are wondering: \textbf{how and to what extent does chemo-repulsion change the dynamics of the model \eqref{Chaplain-Lolas model+} and whether   haptotaxis is also overbalanced by repulsive chemotaxis}? To explore these ambitions, we consider the following no-flux initial-boundary value problem for a chemotaxis-haptotaxis model with mostly repulsive signal production and density-dependent haptotactic sensitivity:
	\be\label{main equation}
	\begin{cases}
		u_t =  \Delta u+\chi \nabla\cdot(u\nabla v)-\xi\nabla\cdot(S(u)\nabla w)  &\text{in } \Omega\times(0,\infty), \\[0.2cm]
		\tau v_t = \Delta v-v+u  &\text{in } \Omega\times(0,\infty),\\[0.2cm]
		w_t=-vw 	 &\text{in } \Omega\times(0,\infty),\\[0.2cm]
		\frac{\partial u}{\partial\nu} +\chi u\frac{\partial v}{\partial\nu}-\xi S(u) \frac{\partial w}{\partial\nu}=\frac{\partial v}{\partial\nu}=0 &\text{on } \partial\Omega\times(0,\infty), \\[0.2cm]
		(u,\tau v,w)=(u_0,\tau v_0,w_0)  &\text{on } \Omega\times\{0\},
	\end{cases}
	\ee
	where the newly introduced  $C^2$-smooth hapto-sensitivity function $S: [0,\infty)\mapsto [0,\infty)$ satisfies, for some $\theta\in\mathbb{R}$, that
 \be\label{S-con}
0=S(0)\leq S(u)\leq u(u+1)^{\theta-1}, \ \ \forall u\geq 0.
 \ee
 We allow nonlinear choice of $S$   to extend the prototypical choice of $S(u)=u$. We mainly consider the case of  $\chi>0$, which may be  biologically interpreted that MDEs become toxic/harmful so that tumor cells stay far away from the toxic MDEs, cf. a similar phenomenon in \cite{DLT23}. Or, in food-chain environment, denoting by $(u,v,w)$  the respective density of (prey, predator, nutrient), or  (insect, rat, cereal) for specific  instance, we may interpret  the model \eqref{main equation} as that prey  diffuses randomly and it stays away from higher concentration of predators and stays toward higher concentration of nutrients, predators diffuse  randomly and decay linearly and the presence of prey is beneficial to the growth of predators, and nutrients  don't diffuse and are  consumed by predators.

Throughout this paper, we  assume that the initial data satisfy
	\be\label{initial-data}
	\begin{cases}
		(u_0,\tau v_0,w_0)\in C(\bar{\Omega})\times W^{1,\infty}(\Omega)\times C^2(\bar{\Omega}),	\\[0.2cm]
		u_0\geq,\not\equiv 0,\ \ \tau v_0\geq 0,\ \ w_0\geq 0 \text{ on } \Omega,\ \ \frac{\partial w_0}{\partial\nu}=0 \text{ on } \partial\Omega.
	\end{cases}
	\ee
Solving the $w$-equation and then taking gradient, one can readily see that
$$
		\nabla w(t)=\nabla w_0e^{-\int_{0}^{t}v(s)ds}-w_0e^{-\int_{0}^{t}v(s)ds}\int_{0}^{t}\nabla v(s)ds,
	$$
and thus, the fact $\nabla w_0\cdot\nu=0$ on $\partial \Omega$ in \eqref{initial-data} quickly implies that the no-flux boundary conditions  for $u$ and $v$   in \eqref{main equation} are equivalent to   the homogeneous Neumann boundary conditions for $u,v$ and $w$:
\be\label{nbc}
 \frac{\partial u}{\partial\nu}=\frac{\partial v}{\partial\nu}=\frac{\partial w}{\partial\nu}=0 \text{ on } \partial \Omega\times (0, \infty).
\ee
In the rest, we shall tacitly use \eqref{nbc} for the boundary conditions of \eqref{main equation}.

	Now, we are in a position to state our main findings on the global dynamical properties of the IBVP  \eqref{main equation}. Our first set of results show the 2D unconditional global existence and boundedness, and then,  illustrate that the haptotaxis-only mechanism cannot  induce any blow-up, and finally, demonstrate the exponential convergence of any supposedly globally bounded solutions to equilibria.

\begin{theorem}\label{theorem1}
	Let $\chi,\xi\geq 0$, $S$ satisfy \eqref{S-con}, $\Omega\subset\mathbb{R}^n(n\geq2)$ be a bounded and smooth domain and let the initial datum $(u_0,\tau v_0,w_0)$ satisfy \eqref{initial-data}.
	\begin{itemize}
		\item[(B1)](\textbf{Global existence and boundedness in 2D})
			Let $n=2$ and $\theta<\frac{3}{2}$ or $\theta<2$ with $\nabla\sqrt{w_0}\in L^\infty(\Omega)$. Then the IBVP \eqref{main equation} possesses a unique global and positive classical solution which is uniformly bounded in time in the sense there exists $C_1>0$ such that
			\be\label{global-boundedness}
				\|u(t)\|_{L^\infty(\Omega)}+\|v(t)\|_{W^{1,\infty}(\Omega)}+\|w(t)\|_{W^{1,\infty}(\Omega)}\leq C_1,\quad\forall t\geq0.
			\ee	
		\item[(B2)](\textbf{Haptotaxis-only can not induce  blow-up})
		For $n\geq 3, \tau=0$, assume $\theta<1+\frac{1}{n}$ or $\theta<1+\frac{2}{n}$ with
	$\nabla\sqrt{w_0}\in L^\infty (\Omega)$. For $n\geq 3, \tau=1$,  assume
\begin{equation*}
\theta< \begin{cases}
 \frac{4}{n}, &\text{ if } 3\leq n\leq 4, \\[0.2cm]
\frac{6n-n^2-4}{2(n-2)}, &\text{ if } n\geq 5.
\end{cases}
 \end{equation*}
 Then the unique classical positive solution $(u,v,w)$ of the IBVP \eqref{main equation} with $\chi=0$ exists globally in time and remains bounded in the sense of \eqref{global-boundedness}.
	\item[(B3)](\textbf{Exponential convergence to constant steady states})
			Any given globally  bounded classical solution component $w$ of   \eqref{main equation} decays exponentially in the sense, for any $\gamma\in(0,v_*)$ with $v_*>0$ defined by \eqref{v-lower-bound-estimate}, there exists $C_2>0$ such that
			\be\label{w-exponential-convergence}
				\|w(t)\|_{W^{1,\infty}(\Omega)}\leq C_2e^{-\gamma t},\quad\forall t>0.
			\nonumber
			\ee
			And  the globally bounded  classical solution component $(u,v)$ of  \eqref{main equation} converges exponentially to its steady state $(\bar{u}_0,\bar{u}_0)$ in the sense, for $p>n$, there exist $C_3>0$ and $\delta>0$ determined by \eqref{delta-def}, such that
			\be\label{2D-convergence2}
				\|u(t)-\bar{u}_0\|_{L^\infty(\Omega)}+\|v(t)-\bar{u}_0\|_{W^{1,\infty}(\Omega)}\leq C_3e^{-\min\{\frac{\delta}{2p},1\} t},\quad\forall t>0.
			\nonumber
			\ee
	\end{itemize}
\end{theorem}
	Motivating by the well-known fact that the chemo-repulsion system \eqref{KS-model-repulsion} with $\tau=0$ admits only global-in-time bounded solutions in any dimensions (cf. \cite{Mock74, Mock75}), we are wondering how such  boundedness can continue to hold in the presence of  haptotactic  interaction.  Our next  result gives a satisfactory answer and shows  that all classical solutions to the IBVP \eqref{main equation} with $\tau=0$ exist globally in time and remain uniformly bounded for a range of  $\theta$ containing $\theta=1$.

\begin{theorem}[\textbf{Blow-up prevention by chemo-repulsion in nD with $\tau=0$}]\label{theorem2}
	Let  $\chi>0,\xi\geq 0$, $S$ satisfy \eqref{S-con}, $\Omega\subset\mathbb{R}^n(n\geq2)$ be a bounded and smooth domain and let the initial datum $(u_0,\tau v_0,w_0)$ satisfy \eqref{initial-data} with $\theta<\frac{3}{2}$ or $\theta<2$ with $\nabla\sqrt{w_0}\in L^\infty(\Omega)$. Then the IBVP \eqref{main equation} with $\tau=0$ admits a unique   positive and classical solution, which is globally bounded in the sense of \eqref{global-boundedness}.
\end{theorem}	
	
	The results in Theorems \ref{theorem1} and \ref{theorem2} reveal, even in the presence of haptotaxis,   chemo-repulsion prevents any formation of   non-trivial patterns in any dimensions and prevents  blow-up of classical solutions in 2D, showing the negligibility of  haptotaxis effect versus chemo-repulsion effect in these respects.
	
	Recall also that the 3D  chemo-repulsion model \eqref{KS-model-repulsion} has  global   weak solutions and they will become  smooth eventually, as seen in Jiang's works \cite{Jiang19, Jiang22}. Here, inspired from the ideas therein, we extend to derive global existence of  weak solutions,  and then their eventual regularity and global stability for the chemo-repulsion-haptotaxis model \eqref{main equation}. To achieve our goals and as a by-product, we first show the exponential stability of constant equilibrium  under some smallness condition on initial data in optimal Lebesgue spaces, which also greatly relax the conditions set in \cite[Theorem 1.7]{RT24JDE}, and extend the corresponding result for chemotaxis-only model \cite[Theorem 1]{Cao15} to chemotaxis-haptotaxis one.

\begin{theorem}[\textbf{Small-data boundedness and exponential convergence}] \label{theorem3}
	Let  $\chi,\xi \in\mathbb{R}$, $S$ satisfy \eqref{S-con} with $\theta\leq 1$, $\Omega\subset\mathbb{R}^n(n\geq2)$ be a bounded and smooth domain, and let the initial datum $(u_0,  \tau v_0,w_0)$ satisfy \eqref{initial-data}. Then there exists $\eta_0>0$ such that whenever
	\be\label{small-initial-data}
		\|u_0\|_{L^\frac{n}{2}(\Omega)}<\eta_0,\  \tau \|\nabla v_0\|_{L^n(\Omega)}<\eta_0\text{ and }\| \nabla w_0\|_{L^\infty(\Omega)}<\eta_0,
	\ee
	the unique classical solution $(u,v,w)$ of the IBVP \eqref{main equation} is globally bounded and decays exponentially to its steady state $(\bar{u}_0,\bar{u}_0,0)$ in the sense, for $t>0$,
	\be\label{3D-convergence1}	\|u(t)-\bar{u}_0\|_{L^\infty(\Omega)}+\|v(t)-\bar{u}_0\|_{W^{1,\infty}(\Omega)}+\|w(t)\|_{W^{1,\infty}(\Omega)}\leq C_4e^{-\min\{\mu,1\}t},
	\ee
where $C_4>0$ is a constant  and  $0<\mu<\min\{\lambda_1,\bar{u}_0\}$ with $\lambda_1$ being the first positive eigenvalue to $-\Delta$ under homogeneous Neumann boundary condition.
\end{theorem}
\begin{remark}Under the following  smallness condition: for some $q>\frac{n}{2}, p>n$,
 \be\label{small-initial-data-rt}
		\|u_0\|_{L^q(\Omega)}<\eta_0,\    \|\nabla v_0\|_{L^p(\Omega)}<\eta_0\text{ and }\|  w_0\|_{W^{1,\infty(\Omega)}}<\eta_0,
	\ee
global boundedness of solutions to \eqref{main equation} with $\tau=1$ and $\chi<0$ is derived in \cite[Theorem 1.7]{RT24JDE}.  By comparison, we not only relax  the smallness condition \eqref{small-initial-data-rt} to   \eqref{small-initial-data} but also offer exponential convergence for such bounded solutions. We especially mention that no smallness on $w_0$ is required in \eqref{small-initial-data}. Theorem \ref{theorem3} also extends \cite{Winkler10JDE, Cao15}'s work on chemotaxis-only model to chemotaxis-haptotaxis model.
\end{remark}
	Next,  we  establish  the global existence of 3D weak solutions to the IBVP \eqref{main equation} with $\tau=1$.  The crucial observation here is that  we detect that the  $L^2$-norm of $\nabla w$ decays exponentially (c.f. Lemma \ref{gradw-l2-decay} below), which in turn allows us to  employ a conditional Lyapunov functional to derive the claimed eventual regularity  and exponential stabilization of the global weak solutions.

\begin{theorem}[\textbf{Global existence, intermediate and eventual regularity,  and exponential convergence of 3D weak solutions}] \label{theorem4}
	Let $\tau=1$,  $\chi,\xi\geq 0$, $S$ satisfy \eqref{S-con} with $\theta<\frac{7}{6}$, $\Omega\subset\mathbb{R}^3$ be a bounded and smooth domain and let the initial datum $(u_0,v_0,w_0)$ satisfy \eqref{initial-data}. Then the IBVP \eqref{main equation} admits a global weak solution $(u,v,w)$ in the sense of Definition \ref{weak-definition} below, which fulfills, for any $T>0$,
\be\label{reg-weak1}
\begin{cases}
	u\in L^\frac{5}{3}(\Omega\times (0, T))\cap L^\frac{5}{4}(0,T;W^{1,\frac{5}{4}}(\Omega)), \\[0.2cm]
  v\in L^\infty(0, T; W^{1,2}(\Omega))\cap  L^2(0,T; W^{1,6}(\Omega))\cap L^2(0,T; W^{2,2}(\Omega)), \\[0.2cm]
  w\in  L^\infty(\Omega\times (0, T))\cap L^\infty(0, T; W^{1,2}(\Omega))\cap L^2(0,T; W^{1,6}(\Omega)).
  \end{cases}
	\ee	
	Moreover, when $\theta\leq\frac{1}{2}$, $S^\prime$ is bounded on $\mathbb{R}^+$  and  $\inf_\Omega v_0>0$,  besides \eqref{reg-weak1}, there exist $t_0>0$  and $\delta\in(0, 1)$ such that the derived global weak solution $(u,v,w)$ satisfies the following refined eventual regularities:
\be\label{reg-weak2}
\begin{cases}
	&u\in  L^\infty(\Omega\times(t_0,\infty))\cap C_{\text{loc}}^{\delta, \frac{\delta}{2}}(\bar\Omega\times [t_0,\infty)),\\[0.2cm]
& v\in  L^\infty(t_0,\infty;W^{1,\infty}(\Omega))\cap  C_{\text{loc}}^{2+\delta, 1+\frac{\delta}{2}}(\bar\Omega\times [t_0,\infty)), \\[0.2cm]
& w \in L^\infty(\Omega\times(0,\infty))\cap L^\infty(t_0,\infty;W^{1,6}(\Omega));
\end{cases}
	\ee	
 and the global weak solution $(u,v,w)$ converges exponentially to the constant equilibrium $(\bar{u}_0,\bar{u}_0,0)$: there exist $\mu>0$ and $C_5>0$ such that, for all $t\geq t_0$,
	\be\label{3D-convergence2}
\begin{split}
&\|u(t)-\bar{u}_0\|_{L^\infty(\Omega)}+\|v(t)-\bar{u}_0\|_{W^{1,\infty}(\Omega)} \\[0.2cm]
&+\|w(t)\|_{L^\infty(\Omega)}+\|\nabla w(t)\|_{L^6(\Omega)}\leq C_5e^{-\mu(t-t_0)}.
\end{split}
	\ee
\end{theorem}
\begin{remark} \label{reg-rek}
 In the absence of haptotaxis, the eventual regularity of $u$ in \eqref{reg-weak2} can be further  improved to be $C_{\text{loc}}^{2+\delta, 1+\frac{\delta}{2}}(\bar\Omega\times [t_0,\infty))$ by the standard   parabolic Schauder theory, and so we have eventual smoothness of weak solutions. In the presence of the ODE interaction of $w$ (which has no spatial regularizing effect), using the best global version of the integrability regularity $v\in L^2(0,\infty; W^{1,6}(\Omega))$ in \eqref{reg-weak1}, we  can only derive the eventual regularity $w\in L^\infty(t_0,\infty; W^{1,6}(\Omega))$ in \eqref{reg-weak2}.
\end{remark}

 By comparison, it can be easily seen that our findings in   Theorems \ref{theorem1},  \ref{theorem2} and \ref{theorem4} (c.f. Remark \ref{reg-rek} for eventual smoothness)  greatly extend the  knowledge about  the corresponding repulsive chemotaxis-only model \eqref{KS-model-repulsion} to the repulsive chemotaxis-haptotaxis model \eqref{main equation}.
	
 In the remaining of this section, we outline the layout as well as main points of this article,
which consists of six main sections. In this section, stemming mainly from  the  regularity-enhancement of chemo-repulsion in  \eqref{KS-model-repulsion}  as well as the almost negligibility of haptotaxis in chemotaxis-haptotaxis models like \eqref{Chaplain-Lolas model} and \eqref{Chaplain-Lolas model+},  we extend  to explore the chemo-repulsion effect in our selected chemotaxis-haptotaxis model  \eqref{main equation}, which also has biological bases. Finally, we
formulate our main findings on the strong regularity-enhancement effect of chemo-repulsion on global dynamics in \eqref{main equation}, in particular, the  blow-up prevention of classical solutions in 2D, the blow-up prevention of weak solutions in 3D and the prevention of pattern formation, in Theorems \ref{theorem1},  \ref{theorem2}, \ref{theorem3} and \ref{theorem4}.

 In section 2, we state the local-in-time existence of classical solutions and present some useful lemmas needed thereafter, some of which (e.g. Lemma \ref{ext-ck}) is of interest for its own sake. Section 3 mainly focuses on the 2D existence and boundedness of classical solutions to \eqref{main equation} by constructing a crucial  energy  identity (cf. Lemma  \ref{lyapunov-lemma}), and we also consider the elliptic case $\tau=0 $ and haptotaxis-only case $\chi=0$ by  subtle  energy methods. In Section 4,  thanks to the positive lower bound of $v$, we   first establish the   $W^{1,\infty}$-exponential decay of $w$ (cf. Lemma \ref{w-convergence-lemma}), and then, we construct a  Lyapunov functional (cf. Lemma \ref{g-functional-lemma}) and use interpolation and semigroup arguments to derive  the global stability of globally bounded solutions. Section 5 mainly uses loop arguments (opennes and closedness argument) to investigate the global  existence, boundedness and exponential decay of solutions to \eqref{main equation} for $\theta\leq1$ under the smallness condition \eqref{small-initial-data}. In Section 6, we first introduce an  approximating system \eqref{regularized-equation}, and then, we show enough uniform a-priori estimates so as to employ  compactness arguments to establish the global existence of weak solutions  to   \eqref{main equation} in 3D. Thereafter, by detecting the key observation  that the  $L^2$-norm of $\nabla w$ decays exponentially (c.f. Lemma \ref{gradw-l2-decay} below), we obtain the stability properties for the linearized system of the approximating system  \eqref{regularized-equation} around   $(\bar{u}_0,\bar{u}_0,0)$ and the stability of   $(\bar{u}_0,\bar{u}_0,0)$ for the approximating system  \eqref{regularized-equation}   in scaling-invariant spaces, with such adequate information, we finally achieve  the eventual regularity and exponential convergence of weak solutions to \eqref{main equation} in 3D.
	
	In the above and below texts, for the sake of convenience and without ambiguity, we use $C$ and $C_i (i=1,2,...)$ to denote positive generic constants depending at most on $u_0, \tau v_0, w_0, \chi, \xi, \theta, \Omega$ and we shall omit integration variables like $dx$, $dt$ etc., identify $u(t) = u(\cdot, t)$ and use frequently abbreviations like, for a function $f$,
	$$
		\|f(t)\|_{L^p}=\|f(\cdot, t)\|_{L^p(\Omega)}=\left(\int_\Omega|f(x,t)|^pdx\right)^\frac{1}{p}, \ \  \ \  \bar f=\frac{1}{|\Omega|}\int_\Omega f.
	$$

\section{Local well-posedness and preliminaries}

	We first state the local well-posedness, extensibility criterion and some basic estimates of local classical solutions to the IBVP \eqref{main equation}.

\begin{lemma}\label{local-existence-lemma}
	 Let $\chi,\xi \in\mathbb{R}$, $S$ satisfy \eqref{S-con},  $\Omega\subset \mathbb{R}^n(n\geq 1)$ be a bounded and smooth domain, and let the initial datum $(u_0, \tau v_0, w_0)$ satisfy \eqref{initial-data}. Then there exists $T_m\in (0,\infty]$ and a uniquely determined pair of nonnegative functions $
		(u,v,w)\in \left[C^{0}(\bar{\Omega}\times [0, T_m))\cap C^{2,1}(\bar{\Omega}\times (0, T_m))\right]^3$ solving  the IBVP \eqref{main equation} classically on $\Omega\times (0, T_m)$ and such that, for $(x,t)\in \bar\Omega\times (0, T_m)$,
	\be\label{u-l1-estimate}
		0<u(x,t),\ \ 0<v(x,t), \ \   w(x,t)=w_0(x)e^{-\int_{0}^{t}v(x,s)ds},
	\ee
	and
	\be\label{v-l1-estimate}
		\|u(t)\|_{L^1}=\|u_0\|_{L^1},  \   \|v(t)\|_{L^1}\leq \max\left\{\tau\|v_0\|_{L^1},    \|u_0\|_{L^1}\right\}, \   0\leq w\leq  \|w_0\|_{L^\infty}.
	\ee
	Moreover, if the maximal existence time $T_m<\infty$, then
	\be\label{extensibility-criterion}
		\limsup_{t\nearrow T_m}\left(\|(u(t)\|_{L^{\infty}}+\|v(t)\|_{W^{1,\infty}}+\|w(t)\|_{W^{1,\infty}}\right)=\infty.
	\ee
\end{lemma}
\begin{proof}
	The local well-posedness and extensibility criterion of classical solutions to the IBVP \eqref{main equation} are  well-established   via a proper contraction mapping argument and the standard parabolic regularity theory, cf.   \cite{Bellomo-Bellouquid-Tao-Winkler15, Jin-Xiang21, MR08, TW11SIAM, Tao-Winkler14JDE, Tao-Winkler14PRSE}.  The positivity of $u$ and $v$ follows from the strong maximum principle and the Hopf boundary point lemma.   Solving the ODE $w_t=-vw$ fulfilled by $w$ in \eqref{main equation}
gives the expression for $w$ in \eqref{u-l1-estimate}, and then  the bounds of $w$ in \eqref{v-l1-estimate}. An integration by parts over $\Omega$ from the first and second equation in \eqref{main equation} and a use of the no-flux boundary conditions directly entail the mass identity of $u$ and the $L^1$-estimate of $v$ in \eqref{v-l1-estimate}.
\end{proof}
	
	The next lemma displays some useful $W^{1,q}$-estimates of $v$  in terms of $L^p$-estimate of $u$, which are derived by  the elliptic estimates if $\tau=0$ or   the standard Neumann heat semigroup estimates for $(e^{t\Delta})_{t\geq0}$ in $\Omega$ if $\tau=1$ \cite{Winkler10JDE, Cao15}, one can refer to \cite[Lemma 2.3]{Jin-Xiang21} for details when $n=2$.

\begin{lemma}\label{v-lq-norm-lemma}
	Let $p\geq1$ and $(u,v,w)$ be the local-in-time classical solution of  \eqref{main equation}. Then for
	\be\label{q-index}
	\begin{cases}
		q_1\in[1,\frac{np}{n-2p}), \quad&\text{if } n\geq 2p, \\[0.2cm]
		q_1\in[1,\infty], \quad&\text{if } n<2p,
	\end{cases}
	\nonumber
	\ee
	there exists $C_1>0$ such that
	\be\label{v-lq-norm-estimate1}
		\|v(t)\|_{L^{q_1}}\leq C_1\left(1+\sup_{s\in(0,t)}\|u(s)\|_{L^p}\right),\quad\forall t\in(0,T_m).
	\ee
	And, for
	\be\label{r-exp-m}
	\begin{cases}
		q_2\in [1, \frac{np}{n-p}), \quad & \text{if } n\geq p,\\[0.2cm]
		q_2\in [1, \infty], \quad & \text{if }  n<p,
	\end{cases}
	\nonumber
	\ee
	there exists $C_2>0$ such that
	\be\label{gradv-lrm}
		\|\nabla v(t)\|_{L^{q_2}}\leq C_2\left(1+\sup_{s\in(0,t)}\|u(s)\|_{L^p}\right),\quad\forall t\in(0, T_m).
	\nonumber
	\ee
	Especially, when $n=2$, there exists $C_3>0$, such that, for $1\leq q<\infty$,
	\be\label{v-lq-bdd2d}
		\|v(t)\|_{L^q}\leq C_3,\quad\forall t\in(0,T_m).
	\ee
\end{lemma}
The following widely known Gagliardo-Nirenberg interpolation inequality \cite{Friedman69, JX18-DCDSB, LL16-non, Xiang-Zheng19Non} as well as its logarithmic  version \cite{Tao-Winkler14JDE, Xiang-Zheng19Non} will be  used in several places in  our forthcoming analysis.

\begin{lemma}[Gagliardo-Nirenberg  interpolation inequality ]\label{GN-lemma}
	Let $\Omega\subset \mathbb{R}^n(n\geq 1)$ be a bounded and smooth domain, let $i$ and $j$ be two integers satisfying $0\leq i<j$, and let  $0<p, q,r\leq \infty$ and $\frac{i}{j}\leq \delta\leq 1$ fulfill
$$
		\frac{1}{p}-\frac{i}{n}=\delta(\frac{1}{q}-\frac{j}{n})+(1-\delta)\frac{1}{r}.
$$
Then there exists a positive constant $C_{GN}=C(i,j, p,q, r,n,\Omega)$ such that
	\be\label{GN-intro}
		\|D^i\varphi\|_{L^p}\leq C_{GN}\left(\|D^j\varphi\|_{L^q}^\delta\|\varphi\|_{L^r}^{1-\delta}+\|\varphi\|_{L^r}\right), \ \ \forall \varphi\in W^{j,q}(\Omega)\cap L^r(\Omega),
	\ee
	with the following exception: if $1 <q<\infty$ and $j-i-\frac{n}{q}$ is a nonnegative integer,
then \eqref{GN-intro} holds only for $\delta$ satisfying $\frac{i}{j}\leq \delta< 1$.
	
In particular, for $n=2$, there exists $C>0$ such that for each $\epsilon>0$ one can pick $C_\epsilon=C(\epsilon,\Omega)>0$ with the property  that
	\be\label{GN-log}
		\|\varphi\|_{L^p}^p\leq \epsilon\|\nabla\varphi\|_{L^2}^{p-q}\|\varphi\ln |\varphi|\|_{L^q}^q+C\|\varphi\|_{L^r}^p+C_\epsilon, \ \ \forall \varphi\in W^{1,2}(\Omega).
	\ee
\end{lemma}

	When  $\tau=0$, then the second equation in \eqref{main equation}  becomes  $-\Delta v+v=u$. Then by $W^{2,p}$-elliptic  estimates and the Gagliardo-Nirenberg interpolation inequality, one gets the following $L^p$-estimate for  $v$, cf. \cite[Lemmas 2.2]{Winkler14} or \cite[Lemma 2.7]{JX18-DCDSB}.

\begin{lemma}
	Let $p>0$, then the local-in-time solution $(u,v,w)$ of \eqref{main equation} with $\tau=0$ satisfies, for any $\eta>0$ and some constant $C(\eta,p)>0$, that
	\be\label{v-lp+1}
		\int_\Omega v^{p+1}\leq\eta\int_\Omega u^{p+1}+C(\eta,p)\left(\int_\Omega u\right)^{p+1},\quad\forall t\in(0,T_m).
	\ee
\end{lemma}

Next, we have  the following  point-wise lower bound for $v$, which follows from the fact that $\|u(t)\|_{L^1}=\|u_0\|_{L^1}$  and the positivity of the associated Green function (case of $\tau=0$) or the heat kernel or lower bound of the Neumann heat semigroup in $\Omega$ (case of $\tau=1$), cf.  \cite[Lemma 2.3]{Xiang22} and  \cite{COY06, Fujie-Senba18, FWY15,  Jin-Xiang21}.

\begin{lemma}\label{v-lower-bound-lemma}
	For $\tau\in\{0,1\}$,   the local-in-time classical solution component $v$ of \eqref{main equation} has a point-wise lower bound in the sense that
	\be\label{v-lbs-gen}
		v(x,t) \geq \tau \|u_0\|_{L^1}\zeta(t)+(1-\tau)\|u_0\|_{L^1}\zeta_0,\quad\forall (x,t)\in \Omega\times(0,T_m),
	\ee
	where  the nonnegative bounded  function  $\zeta: [0,\infty)\mapsto [0,\infty)$ is strictly increasing with $\zeta(0)=0$ and $\zeta_0=\zeta(\infty)$. In particular,  if $\Omega$ is convex-plus (containing convex) \cite[Definition 2.1]{Jin-Xiang21}, then $\zeta$ is given by
\begin{equation}\label{zeta-def}
\zeta(t)=\int_0^t \frac{1}{(4\pi s)^{\frac{n}{2}}}e^{-\left(s+\frac{d^2}{4s}\right)}ds, \ \ \ \ d:=\text{diam}(\Omega)=\sup_{x,y\in\Omega}|x-y|, \ \ \zeta_0=\zeta(\infty).
\end{equation}
It then follows from \eqref{v-lbs-gen}, for any $\sigma\in(0, T_m)$,  that $v$ has a positive lower bound over $\Omega\times (\sigma, T_m)$ as follows:
\begin{equation}\label{v-lower-bound-estimate0}
 v(x,t)\geq \tau \|u_0\|_{L^1}\zeta(\sigma)+(1-\tau)\|u_0\|_{L^1}\zeta_0=:v_*(\sigma)>0.
\end{equation}
\end{lemma}
\begin{remark}By the  extensibility criterion in \eqref{extensibility-criterion} of Lemma \ref{local-existence-lemma}, upon shitting the initial time $t=0$ to $t=\sigma$, we shall henceforth assume for convenience that \eqref{v-lower-bound-estimate0} holds on $\Omega\times(0,T_m)$; namely, there exists $v_*>0$ such that
 \be\label{v-lower-bound-estimate}
 v(x,t)\geq v_*>0,\quad\forall (x,t)\in\Omega\times(0,T_m).
\ee
\end{remark}
 By the expression of $w$ in \eqref{u-l1-estimate} and the positive lower bound of $v$, we derive several more properties of $w$, cf. \cite{Cao16, Jin-Xiang21, Tao14JO, Xiang22}.

\begin{lemma}\label{w-lemma}
	For  $t\in(0,T_m)$, the classical solution component $w$ satisfies
	\be\label{w-expression}
		0\leq w(x,t)=w_0(x)e^{-\int_{0}^{t}v(x,s)ds}\leq \|w_0\|_{L^\infty}e^{-v_*t},
	\ee
	and for given $p\in(1,\infty]$,
	\be\label{w-lp}
		\|\nabla w(t)\|_{L^p}\leq\|\nabla w_0\|_{L^p}e^{-v_*t}+\|w_0\|_{L^\infty} e^{-v_*t}\int_0^t\|\nabla v(s)\|_{L^p}ds.
	\ee
In particular, for any $0<\mu<v_*$, there holds
	\be\label{grad w-lp0}
		\int_\Omega|\nabla w(t)|^p\leq e^{-p\mu t}\int_\Omega|\nabla w_0|^p+\frac{\|w_0\|_{L^\infty}^p}{(v_*-\mu)^{p-1}}\int_{0}^{t}e^{-p\mu(t-s)}\int_\Omega|\nabla v(s)|^pds.
	\ee
	Moreover, $-\Delta w$ has a general one-sided point-wise upper bound as follows:
	\be\label{delta-w-estimate0}
		-\Delta w\leq\tau\|w_0\|_{L^\infty}v(t)+\frac{\|w_0\|_{L^\infty}}{e}+2e^{-\int_{0}^{t}v(s)ds}\nabla w_0\int_{0}^{t}\nabla v(s)ds+\|\Delta w_0\|_{L^\infty}.
	\ee	
	Furthermore, if $\nabla \sqrt{w_0}\in L^\infty(\Omega)$, then the estimate \eqref{delta-w-estimate0} sharpens as follows:
	\be\label{delta-w-estimate}
	-\Delta w(t)\leq\tau\|w_0\|_{L^\infty}v(t)+\frac{\|w_0\|_{L^\infty}}{e}+4\|\nabla\sqrt{w_0}\|_{L^\infty}^2+\|\Delta w_0\|_{L^\infty}.
	\ee
\end{lemma}
\begin{proof}
	Solving the $w$-equation in \eqref{main equation} along with $v\geq v_*$ in \eqref{v-lower-bound-estimate} on $\Omega\times(0,T_m)$  yields \eqref{w-expression}. Then we compute from the expression for $w$ in \eqref{w-expression} that
	$$
		\nabla w(t)=\nabla w_0e^{-\int_{0}^{t}v(s)ds}-w_0e^{-\int_{0}^{t}v(s)ds}\int_{0}^{t}\nabla v(s)ds,
	$$
	which along with $0\leq w\leq \|w_0\|_{L^\infty}$ in \eqref{v-l1-estimate} and $v\geq v_*$ in \eqref{v-lower-bound-estimate} on $\Omega\times(0,T_m)$ readily implies the $L^p$-estimate \eqref{w-lp}.

	To illustrate \eqref{grad w-lp0}, for  $\epsilon>0$, since $w_t=-vw$, we calculate by the lower bound of $v$ in \eqref{v-lower-bound-estimate} and Young's inequality with epsilon that
	\be
	\begin{split}
		\frac{1}{p}\frac{d}{dt}\int_\Omega|\nabla w|^p&=-\int_\Omega w|\nabla w|^{p-2}\nabla v\cdot\nabla w-\int_\Omega v|\nabla w|^p\\
		&\leq-(v_*-\epsilon)\int_\Omega|\nabla w|^p+\frac{\|w_0\|_{L^\infty}^p}{p\epsilon^{p-1}}\int_\Omega|\nabla v|^p, 	
	\end{split}
	\nonumber
	\ee
which upon setting $\epsilon=v_*-\mu>0$ implies
$$
\frac{d}{dt}\left(e^{p\mu t}\int_\Omega|\nabla w|^p\right)\leq \frac{\|w_0\|_{L^\infty}^p}{(v_*-\mu)^{p-1}}e^{p\mu t}\int_\Omega|\nabla v|^p.
$$
 Integrating  this ODI (ordinary differential inequality) from $0$ to $t$ entails  \eqref{grad w-lp0}.

 Finally, the one-sided point-wise upper bounds of $-\Delta w$ in \eqref{delta-w-estimate0} and \eqref{delta-w-estimate} can be obtained by simple calculations from \eqref{w-expression}, see \cite[Lemma 2.2]{Xiang22}.
\end{proof}

	In several places, we also need to use the following  estimate for integrals \cite{Cao-Lankeit16, Winkler10JDE}, which arises in calculations involving semigroup representations of solutions.
	
\begin{lemma}\label{integral-lemma}
	For all $\eta>0$, there exists $C=C(\eta)>0$ such that for all $\alpha\in[0,1-\eta]$, $\beta\in[\eta,1-\eta]$, $\gamma, \delta\in\mathbb{R}$ fulfilling $\eta\leq\gamma-\delta\leq\frac{1}{\eta}$ and for all $t>0$, we have
	\be\label{algebraic-estimate}
	\int_0^t(1+s^{-\alpha})(1+(t-s)^{-\beta})e^{-\gamma s}e^{-\delta(t-s)}ds\leq C(\eta)\left(1+t^{\min\{0,1-\alpha-\beta\}}\right)e^{-\min\{\gamma,\delta\}t}.
	\nonumber
	\ee
\end{lemma}

To establish the eventual regularity of weak solutions in 3D, we need a crucial logarithmic-type Sobolev inequality, which is of interest for its own sake.
\begin{lemma}\label{ext-ck}
Let $\Omega\subset \mathbb{R}^n$ be a bounded and smooth domain with $1\leq n\leq 6$. Then
\be\label{sharp-ck}
\frac{1}{k}:=\inf\left\{\frac{\int_\Omega |\nabla \phi^\frac{1}{2}|^2}{\int_\Omega \phi\ln \phi}:  \ \ \phi> 0, \int_\Omega \phi=|\Omega|\right\}>0,
\ee
so that
\be\label{sharp-ck+}
 \int_\Omega \psi\ln \frac{\psi}{\bar\psi}\leq k\int_\Omega |\nabla \psi^\frac{1}{2}|^2 , \ \ \forall \psi>0, \ \ \psi\in L^1(\Omega).
\ee
\end{lemma}
\begin{proof}For $n\leq2$, the proof can be shown via the  continuous Sobolev embedding $W^{1,1}(\Omega)\hookrightarrow L^2(\Omega)$, Poincar\'{e}'s  inequality as well as the algebraic inequality:
		$$
			s\ln s\leq s-1+(s-1)^2,\quad\forall s>0.
		$$
In particular, for $\phi>0$ with $\int_\Omega \phi=|\Omega|$, one has, for some $k>0$, that
		$$
			\int_\Omega \phi \ln \phi\leq\int_\Omega(\phi-1)+\int_\Omega(\phi-1)^2\leq k \int_\Omega|\nabla \phi^\frac{1}{2}|^2.
		$$
In general, for $1\leq n\leq 6$ and  $z\geq 0$ with $\int_\Omega z^2=|\Omega|$, since the function $s\in (0, \infty)\mapsto \ln s$ is concave,  it follows from Jensen's inequality that
\be\label{z-jensen}
\int_\Omega z^2\ln z^2=2|\Omega|\int_\Omega \frac{z^2}{|\Omega|}\ln z\leq 2|\Omega|\ln \left(\int_\Omega \frac{z^3}{|\Omega|}\right).
\ee
To estimate the integral on the right-hand side of \eqref{z-jensen}, we compute that
\be\label{z3-int}
\begin{split}
\int_\Omega \frac{z^3}{|\Omega|}&=\frac{1}{|\Omega|}\int_\Omega\left[(z-\bar z)^3+3\bar z( z-\bar z)^2+3\bar z^2( z-\bar z)+\bar z^3\right]\\
&=\frac{1}{|\Omega|}\int_\Omega(z-\bar z)^3+\frac{3\bar z}{|\Omega|}\int_\Omega ( z-\bar z)^2+ \bar z^3.
\end{split}
\ee
Recalling $\int_\Omega z^2=|\Omega|$, then H\"{o}lder's  inequality shows that
\be\label{hl-inequality}
\bar z=\frac{1}{|\Omega|}\int_\Omega z\leq \frac{1}{|\Omega|^\frac{1}{2}}\left(\int_\Omega z^2\right)^\frac{1}{2}=1.
\ee
Based on the continuous Sobolev embedding $W^{1,2}(\Omega)\hookrightarrow L^q(\Omega)$ for bounded smooth domain $\Omega\subset \mathbb{R}^n(n\leq6)$ and $q\in \{2,3\}$, one can argue by contradiction there exist  $\mu_1, \mu_2>0$ such that the following Poincar\'{e} type inequalities hold:
\be\label{poincare-ineq23}
\begin{cases}
	\|z-\bar z\|_{L^2}^2\leq \mu_1 \|\nabla z\|_{L^2}^2,\\[0.25cm]
\|z-\bar z\|_{L^3}^3\leq \mu_2 \|\nabla z\|_{L^2}^3.
\end{cases}
\ee
Plugging \eqref{z3-int},  \eqref{hl-inequality} and  \eqref{poincare-ineq23} into \eqref{z-jensen}, we infer that
\be\label{z3-int+}
\begin{split}
\int_\Omega z^2\ln z^2&\leq 2|\Omega|\ln \left(\int_\Omega \frac{z^3}{|\Omega|}\right)\\
& \leq  2|\Omega|\ln\left(1+\frac{3\mu_1}{|\Omega|}\|\nabla z\|_{L^2}^2+\frac{\mu_2}{|\Omega|}\|\nabla z\|_{L^2}^3\right)\\
&\leq 2|\Omega|\mu \|\nabla z\|_{L^2}^2=2|\Omega|\mu\int_\Omega |\nabla z|^2,
\end{split}
\ee
where $\mu$ is a positive and finite number  and is given by
$$
\mu:=\sup_{s>0}\frac{\ln\left(1+\frac{3\mu_1}{|\Omega|}s+\frac{\mu_2}{|\Omega|}s^\frac{3}{2}\right)}{s}.
$$
Finally, given  $\phi>0$ with $\int_\Omega \phi=|\Omega|$, a simple substitution of $z=\phi^\frac{1}{2}$ in \eqref{z3-int+} yields \eqref{sharp-ck},  and then, the estimate \eqref{sharp-ck} follows upon setting $\phi=\frac{\psi}{\bar\psi}$ in \eqref{z3-int+}.
\end{proof}

When deriving the eventual regularity of weak solutions in 3D, the following  elementary  result of analysis obtained via contradictory argument  will also be  used \cite[Lemma 2.8]{Jiang19}.

\begin{lemma}\label{limit-lemma}
	Assume that $a\geq 0$ and $0\leq f_k(t)\in L^1(a,\infty)$ with $k\in\mathbb{N}$ is a sequence of functions such that
	$$
		\sup_{k\in\mathbb{N}}\int_a^{\infty}f_k(s)ds<\infty.
	$$
 Then, for arbitrary $\sigma>0$, there is $t_\sigma>a$ such that for any $k$, there exists a time $t_\sigma^k\in(a,t_\sigma)$ such that $f_k(t_\sigma^k)\leq\sigma$.
\end{lemma}

\section{Global existence and boundedness for $n\geq2$}
In this section, we aim to show global existence and boundedness of classical solutions to the IBVP \eqref{main equation} for $n\geq2$. Our analysis relies on the following key evolution identity associated with \eqref{main equation}.
\begin{lemma}\label{lyapunov-lemma}
	The local-in-time classical solution $(u,v,w)$ of the  IBVP  \eqref{main equation} fulfills, for $t\in(0,T_m)$,
	\be\label{functional1}
		\mathcal{F}^\prime(t)+4\int_\Omega |\nabla u^\frac{1}{2}|^2+\chi\int_\Omega\left(|\Delta v|^2+|\nabla v|^2\right)
		=\xi\int_\Omega \frac{S(u)}{u}\nabla u\cdot\nabla w,
	\ee
 	where $\mathcal{F}(t)$ is defined by
	\be\label{F-definition}
	\mathcal{F}(t)=\int_\Omega u\ln u+\frac{\tau\chi}{2}\int_\Omega|\nabla v|^2.
	\ee	
\end{lemma}
\begin{proof}
	We first test the $u$-equation in \eqref{main equation} by $\ln u$, and then integrate by parts over $\Omega$ to see  that
	\be\label{ulnu}
		\frac{d}{dt}\int_\Omega u\ln u=-4\int_\Omega |\nabla u^\frac{1}{2}|^2-\chi\int_\Omega\nabla u\cdot\nabla v+\xi\int_\Omega\frac{S(u)}{u}\nabla u\cdot\nabla w.
	\ee
	Next, taking $\nabla$ on both sides of the $v$-equation in \eqref{main equation}, we find  $\tau\nabla v_t=\nabla\Delta v-\nabla v+\nabla u$; then multiplying by $\nabla v$, we infer from integration by parts that
	\be\label{grad v-l2}
		\frac{\tau}{2}\frac{d}{dt}\int_\Omega |\nabla v|^2=-\int_\Omega|\Delta v|^2-\int_\Omega|\nabla v|^2+\int_\Omega\nabla u\cdot\nabla v.
	\ee
	Combining \eqref{ulnu} and \eqref{grad v-l2}, we readily derive \eqref{functional1} with $\mathcal{F}(t)$ given by \eqref{F-definition}.
\end{proof}

\subsection{Global existence and boundedness in 2D with $\tau=1$}

	In this subsection, we explore the 2D global boundedness of classical solutions to the IBVP \eqref{main equation} with $\tau=1$, thus proving Theorem \ref{theorem1} (B1). Thanks to \eqref{functional1}, we use the one-sided point-wise upper bound of $-\Delta w$ to derive a Gronwall inequality, which enables us to improve the $L^1$-boundedness of $u$ to the $L^1$-boundedness of $u\ln u$. Thus invoking the Moser iteration or Neumann heat semigroup, we finally conclude the global boundedness of solutions in 2D.  We proceed to estimate the term on the right-hand side of \eqref{functional1} by means of  the one-sided point-wise upper bounds of $-\Delta w$ shown in \eqref{delta-w-estimate0} and \eqref{delta-w-estimate}.

\begin{lemma}\label{grad u-grad w-lemma}
	Let $n=2$ and $\tau=1$. Assume $\theta<2$  with  $\nabla \sqrt{w_0}\in L^\infty(\Omega)$ or $\theta<\frac{3}{2}$ without any further condition on the initial datum $w_0$. Then for any $\eta>0$, there exists  $C_\eta>0$ such that the solution component $(u,w)$ satisfies
	\be\label{grad u-grad w-estimate}
		\xi\int_\Omega \frac{S(u)}{u}\nabla u\cdot\nabla w\leq  \eta\int_\Omega|\nabla u^\frac{1}{2}|^2+C_\eta.
	\ee	
\end{lemma}
\begin{proof}
By the  growth condition  $0\leq S(u)\leq u(u+1)^{\theta-1}$, it follows that
	\begin{equation*}				
0\leq\Phi_1(u):=\int_{0}^{u}\frac{S(\sigma)}{\sigma} d\sigma\leq\begin{cases}
\frac{(u+1)^\theta-1}{\theta}, &\text{ if } \theta>0, \\[0.2cm]
\ln (1+u), &\text{ if } \theta=0, \\[0.2cm]
\frac{1-(u+1)^\theta}{-\theta}, &\text{ if } \theta<0.
\end{cases}
	\end{equation*}
Notice that $\theta \ln (1+u)\leq (1+u)^\theta-1$ for all $\theta>0, u\geq0$; we therefore shall  proceed with the assumption $\theta>0$ and then we have
\be\label{phi1}						
0\leq\Phi_1(u)=\int_{0}^{u}\frac{S(\sigma)}{\sigma} d\sigma\leq
\frac{(u+1)^\theta}{\theta}.
\ee
Then an integration by parts gives
	\be\label{hap-est2d}
		\xi\int_\Omega \frac{S(u)}{u}\nabla u\cdot\nabla w=-\xi\int_\Omega\Phi_1(u)\Delta w.
	\ee
	
Next, we proceed to  estimate the right-hand side with/without  condition on $w_0$.
	
	\textbf{Case 1:} $\theta<2$ and $\nabla\sqrt{w_0}\in L^\infty(\Omega)$. In this case, by the estimate \eqref{phi1}	and the pointwise upper bound of $-\Delta w$ in \eqref{delta-w-estimate} with $\tau=1$, we get that
	\be\label{3-estimate1}
		-\xi\int_\Omega\Phi_1(u)\Delta w\leq \xi\|w_0\|_{L^\infty}\theta^{-1}\int_\Omega(u+1)^\theta v+C_1\int_\Omega(u+1)^\theta.
	\ee
	Then applying the 2D Gagliardo-Nirenberg inequality in Lemma \ref{GN-lemma} and the mass conservation $\|u\|_{L^1}=\|u_0\|_{L^1}$, we infer  that
	\be\label{u-l2-estimate}
	\begin{split}
		\|u+1\|_{L^2}^2&=\|\sqrt{u}\|_{L^4}^4+2\|u_0\|_{L^1}+|\Omega|\\
		&\leq C_2\left(\|\nabla\sqrt{u}\|_{L^2}^2\|\sqrt{u}\|_{L^2}^2+ \|u\|_{L^1}^4\right)+2\|u_0\|_{L^1}+|\Omega|\\
		&\leq C_3\|\nabla\sqrt{u}\|_{L^2}^2+C_3,
	\end{split}
	\ee
	which along with the Young's inequality and $\theta<2$ implies, for any $\eta>0$, that
	\be\label{3-estimate2}
	\begin{split}
		\xi\|w_0\|_{L^\infty}\theta^{-1}\int_\Omega (u+1)^\theta v&\leq  \frac{\eta}{2C_3}\int_\Omega(u+1)^2+C_{\eta}\int_\Omega v^{\frac{2}{2-\theta}}\\
		&\leq\frac{\eta}{2}\int_\Omega|\nabla u^{\frac{1}{2}}|^2+\frac{\eta}{2}+ C_{\eta}\int_\Omega v^{\frac{2}{2-\theta}}
	\end{split}
	\ee
	and
	\be\label{3-estimate3}
	\begin{split}
		C_1\int_\Omega(u+1)^\theta&\leq\frac{\eta}{2C_3}\int_\Omega(u+1)^2+C_{\eta}\\
		&\leq\frac{\eta}{2}\int_\Omega|\nabla u^{\frac{1}{2}}|^2+\frac{\eta}{2}+C_\eta.
	\end{split}
	\ee
	Substituting \eqref{3-estimate1}, \eqref{3-estimate2} and  \eqref{3-estimate3} into  \eqref{hap-est2d}, and using the $L^q$-bound of $v$ with $q=\frac{2}{2-\theta}\in[1,\infty)$ in \eqref{v-lq-bdd2d} of  Lemma \ref{v-lq-norm-lemma} with $n=2$, we conclude the estimate \eqref{grad u-grad w-estimate}.
	
	\textbf{Case 2:} $\theta<\frac{3}{2}$ and without any further condition on $w_0$. In light of the general upper bound for $-\Delta w$ in \eqref{delta-w-estimate0}, we derive from \eqref{phi1} that
	\be\label{est1}
	\begin{split}
		-\xi\int_\Omega\Phi_1(u)\Delta w\leq&\xi\theta^{-1}\|w_0\|_{L^\infty}\int_\Omega (u+1)^\theta v\\
& \  +\xi\theta^{-1}\left(e^{-1}\|w_0\|_{L^\infty} +\|\Delta w_0\|_{L^\infty}\right)\int_\Omega (u+1)^\theta \\
		&-2\xi\theta^{-1}\int_\Omega (u+1)^\theta\nabla e^{-\int_{0}^{t} v(s)ds}\cdot\nabla w_0.
	\end{split}
	\ee
To estimate the last term, we employ  integration by parts due to  $\frac{\partial w_0}{\partial\nu}=0$ on $\partial\Omega$, the Young's inequality and  the nonnegativity of $v$ to infer,  for any $\eta>0$, that
	\be\label{est2}
	\begin{split}
		&-2\xi\theta^{-1}\int_\Omega (u+1)^\theta\nabla e^{-\int_{0}^{t} v(s)ds}\cdot\nabla w_0	\\
		&=2\xi\int_\Omega (u+1)^{\theta-1} e^{-\int_{0}^{t} v(s)ds}\nabla u\cdot\nabla w_0\\
&\ +2\xi\theta^{-1}\int_\Omega(u+1)^\theta e^{-\int_{0}^{t} v(s)ds}\Delta w_0 \\
		&\leq2\xi\|\nabla w_0\|_{L^\infty}\int_\Omega (u+1)^{\theta-1}|\nabla u|+2\xi\theta^{-1}\|\Delta w_0\|_{L^\infty}\int_\Omega (u+1)^\theta	\\
		&\leq \frac{\eta}{2}\int_\Omega\frac{|\nabla u|^2}{u}+\frac{2}{\eta}\xi^2\|\nabla w_0\|_{L^\infty}^2\int_\Omega(u+1)^{2\theta-1}\\
&+2\xi\theta^{-1}\|\Delta w_0\|_{L^\infty}\int_\Omega (u+1)^\theta.
	\end{split}
	\ee
Since $\max\{\theta, \ 2\theta-1\}<2$, similar to  the estimates \eqref{3-estimate2} and \eqref{3-estimate3}, we use  Young's inequality and \eqref{u-l2-estimate} to deduce, for  any $\eta>0$,  that
	\be\label{est3}
	\begin{split}
	&\xi\theta^{-1}\|w_0\|_{L^\infty}\int_\Omega (u+1)^\theta v+\xi\theta^{-1}\left(e^{-1}\|w_0\|_{L^\infty} +\|\Delta w_0\|_{L^\infty}\right)\int_\Omega (u+1)^\theta \\
&+\frac{2}{\eta}\xi^2\|\nabla w_0\|_{L^\infty}^2\int_\Omega(u+1)^{2\theta-1}+2\xi\theta^{-1}\|\Delta w_0\|_{L^\infty}\int_\Omega (u+1)^\theta\\
&\leq \frac{\eta}{2}\int_\Omega\frac{|\nabla u|^2}{u}+C_\eta.
	\end{split}		
	\ee
	Substituting \eqref{est1},  \eqref{est2} and \eqref{est3} into \eqref{hap-est2d}, we end up with  \eqref{grad u-grad w-estimate}.
\end{proof}

	With the key estimate provided by Lemma \ref{grad u-grad w-lemma}, we are  now in a position to achieve the $L^1$-boundedness of $u\ln u$.
	
\begin{lemma}\label{ulnu-lemma}
	Let $n=2$ and $\tau=1$. Assume $\theta<2$  with  $\nabla \sqrt{w_0}\in L^\infty(\Omega)$ or $\theta<\frac{3}{2}$. Then there exists $C>0$ such that
	\be\label{ulnu-grad v-delta v-estimate}
		\|(u\ln u)(t)\|_{L^1}+\|\nabla v(t)\|_{L^2}^2\leq C,\quad\forall t\in(0,T_m).
	\ee
\end{lemma}
\begin{proof}
Submitting the estimate \eqref{grad u-grad w-estimate} with $\eta=1$ into the evolution identity \eqref{functional1}, we find, for all $t\in (0, T_m)$, that
	\be\label{ulnu1}
		\mathcal{F}^\prime(t)+3\int_\Omega|\nabla u^\frac{1}{2}|^2+\chi\int_\Omega\left(|\Delta v|^2+|\nabla v|^2\right)
		\leq C_1.
	\ee
	For any $\epsilon>0$, notice the elementary inequality that
	\be\label{slns}
		s\ln s\leq \epsilon (s+1)^2+C_\epsilon, \ \ \  C_\epsilon=\sup_{s>0}\{s\ln s-\epsilon (s+1)^2\}\in \mathbb{R}.
	\ee
Choosing sufficiently small  $\epsilon>0$ in \eqref{slns}, we can easily infer from  \eqref{u-l2-estimate} that,
	\be\label{ulnu2}
		\int_\Omega u\ln u\leq \int_\Omega|\nabla u^{\frac{1}{2}}|^2+C_2.
	\ee
	Adding \eqref{ulnu1} and \eqref{ulnu2} and recalling the definition of $\mathcal{F}$ in \eqref{F-definition}, we derive that
$$
		\mathcal{F}^\prime(t)+\mathcal{F}(t)+2\int_\Omega|\nabla u^{\frac{1}{2}}|^2+\chi\int_\Omega|\Delta v|^2+\frac{\chi}{2}\int_\Omega|\nabla v|^2\leq C_3,
$$
	which upon being multiplied by $e^t$ and then being integrated, yields immediately
	\be\label{ulnu4}
		\int_\Omega u\ln u+\frac{\chi}{2}\int_\Omega|\nabla v|^2=\mathcal{F}(t)\leq C_4,\quad\forall t\in(0,T_m).
	\ee
	In light of  \eqref{ulnu4} and the fact that $-s\ln s\leq e^{-1}$ for $s>0$, we deduce that
	\be\label{ulnu5}
		\int_\Omega|u\ln u|=\int_\Omega u\ln u-2\int_{\{u\leq1\}}u\ln u\leq\mathcal{F}(t)+2e^{-1}|\Omega|\leq C_5,
	\nonumber
	\ee
	which along with \eqref{ulnu4} yields our desired result \eqref{ulnu-grad v-delta v-estimate}.
\end{proof}	
	
	In the 2D setting, based on the compound estimate obtained in \eqref{ulnu-grad v-delta v-estimate}, it is quite standard to derive the uniform $L^p$-boundedness of $u$ through energy estimates and then $L^\infty$-norm via Moser iteration or Neumann heat semigroup estimate.
	
\begin{lemma}\label{global-boundedness-lemma}
	Let $n=2$, $\tau=1$ and $\theta<2$  with  $\nabla \sqrt{w_0}\in L^\infty(\Omega)$ or $\theta<\frac{3}{2}$. Then the solution  of \eqref{main equation}  exists globally and is uniformly bounded according to \eqref{global-boundedness}.
\end{lemma}
\begin{proof}
	\textbf{Step 1:} $L^2$-boundedness of $u$.  Testing the $u$-equation in \eqref{main equation} by $u$ and then integrating by parts over $\Omega$, we derive that
	\be\label{u-l2-1}
	\begin{split}
		\frac{d}{dt}\int_\Omega u^2+2\int_\Omega|\nabla u|^2&=-2\chi\int_\Omega u\nabla u\cdot\nabla v+2\xi\int_\Omega S(u)\nabla u\cdot\nabla w	\\
		&=-2\chi\int_\Omega u\nabla u\cdot\nabla v-2\xi\int_\Omega\Phi_2(u)\Delta w\\
&\leq \frac{1}{2}\int_\Omega|\nabla u|^2+2\chi^2\int_\Omega u^2|\nabla v|^2-2\xi\int_\Omega\Phi_2(u)\Delta w,
	\end{split}	
	\ee
	where, as \eqref{phi1}, the growth condition $0\leq S(u)\leq u(u+1)^{\theta-1}$ and the  elementary inequality $(u+1)^{\theta+1}\leq 2^\theta(u^{\theta+1}+1)$ shows that
\be\label{phi2}
0\leq \Phi_2(u):=\int_{0}^{u}S(\sigma)d\sigma\leq\int_{0}^{u}(\sigma+1)^\theta d\sigma\leq\frac{2^\theta}{\theta+1}(u^{\theta+1}+1).
\ee

	\textbf{Case 1:} $\theta<2$ and $\nabla\sqrt{w_0}\in L^\infty(\Omega)$. In this case, by  the pointwise upper bound of $-\Delta w$ in \eqref{delta-w-estimate}, we conclude from \eqref{phi2} and the $L^1$-boundedness of $v$  that
	\be\label{u-l2-30}
	\begin{split}
		&-2\xi\int_\Omega\Phi_2(u)\Delta w\\
&\leq\frac{2^\theta\xi}{\theta+1}\int_\Omega\left(u^{\theta+1}+1\right) \|w_0\|_{L^\infty}\cdot v\\
&+\frac{2^\theta\xi}{\theta+1}\int_\Omega\left(u^{\theta+1}+1\right)\left(\frac{\|w_0\|_{L^\infty}}{e}+4\|\nabla\sqrt{w_0}\|_{L^\infty}^2+\|\Delta w_0\|_{L^\infty}\right)	\\
		&\leq C_1\int_\Omega u^{\theta+1}v+C_2\int_\Omega u^{\theta+1}+C_3.
	\end{split}
	\ee
	Invoking the Young's inequality, we deduce, for $\theta<2$, that
	\be\label{u-l2-31}
		C_1\int_\Omega u^{\theta+1}v\leq\frac{1}{4}\int_\Omega u^3+C_4\int_\Omega v^\frac{3}{2-\theta}
	\ee
	and
	\be\label{u-l2-32}
		C_2\int_\Omega u^{\theta+1}\leq \frac{1}{4}\int_\Omega u^3+C_5.
	\ee
	Noticing   the $L^q$-bound of $v$ with $q=\frac{3}{2-\theta}\in[1,\infty)$ in \eqref{v-lq-norm-estimate1},  and substituting \eqref{u-l2-31} and \eqref{u-l2-32} into \eqref{u-l2-30}, we then deduce that
	\be\label{u-l2-3}
		-2\xi\int_\Omega\Phi_2(u)\Delta w\leq\frac{1}{2}\int_\Omega u^3+C_6.
	\ee	
	
	\textbf{Case 2:} $\theta<\frac{3}{2}$ and without further condition on $w_0$. We use the bound of $\Phi_2$ in \eqref{phi2} and the general point-wise upper bound of $-\Delta w$ in \eqref{delta-w-estimate0} to estimate that
	\be\label{u-l2-5}
	\begin{split}
		-2\xi\int_\Omega\Phi_2(u)\Delta w&\leq\frac{2^\theta\xi\|w_0\|_{L^\infty}}{\theta+1}\int_\Omega(u^{\theta+1}+1)v	\\
		&\ \ \ +\frac{2^\theta\xi}{\theta+1}\left(\frac{\|w_0\|_{L^\infty}}{e}+\|\Delta w_0\|_{L^\infty}\right)\int_\Omega(u^{\theta+1}+1)	\\
		&\ \ \ -\frac{2^{\theta+1}\xi}{\theta+1}\int_\Omega(u^{\theta+1}+1)\nabla w_0\cdot\nabla e^{-\int_{0}^{t}v(s)ds}.
	\end{split}
	\ee
	Again, using  integration by parts, Young's inequality and the nonnegativity of $v$, we derive that
	\be\label{u-l2-6}
	\begin{split}
		&-\frac{2^{\theta+1}\xi}{\theta+1}\int_\Omega(u^{\theta+1}+1)\nabla w_0\cdot\nabla e^{-\int_{0}^{t}v(s)ds}	\\
		&=2^{\theta+1}\xi\int_\Omega u^\theta e^{-\int_{0}^{t}v(s)ds}\nabla u\cdot\nabla w_0\\
&\ \ +\frac{2^{\theta+1}\xi}{\theta+1}\int_\Omega(u^{\theta+1}+1)e^{-\int_{0}^{t}v(s)ds}\Delta w_0	\\
		&\leq 2^{\theta+1}\xi\|\nabla w_0\|_{L^\infty}\int_\Omega u^\theta|\nabla u|+\frac{2^{\theta+1}\xi}{\theta+1}\|\Delta w_0\|_{L^\infty}\int_\Omega(u^{\theta+1}+1)	\\
		&\leq\frac{1}{2}\int_\Omega |\nabla u|^2+C_7\int_\Omega u^{2\theta}+C_8\int_\Omega u^{\theta+1}+C_9\\
		&\leq\frac{1}{2}\int_\Omega |\nabla u|^2+\frac{1}{2}\int_\Omega u^3+C_{10},
	\end{split}
	\ee
where the fact that  $\max\{2\theta, \theta+1\}<3$ is used.

	Substituting \eqref{u-l2-6} into \eqref{u-l2-5} and applying similar estimates to \eqref{u-l2-31} and \eqref{u-l2-32}, we quickly infer that
	\be\label{u-l2-33}
		-2\xi\int_\Omega\Phi_2(u)\Delta w\leq\frac{1}{2}\int_\Omega |\nabla u|^2+\frac{1}{2}\int_\Omega u^3+C_{11}.
	\ee
	Then plugging  \eqref{u-l2-3} in \textbf{Case 1} or \eqref{u-l2-33} in \textbf{Case 2} into \eqref{u-l2-1} and using  the elementary fact that $s^2\leq \frac{1}{2}s^3+1$ for $s\geq0$, we arrive at
\be\label{u-l2-4}
	\frac{d}{dt}\int_\Omega u^2+\int_\Omega u^2+\int_\Omega|\nabla u|^2\leq2\chi^2\int_\Omega u^2|\nabla v|^2+\int_\Omega u^3+C_{12}.
	\ee
	To digest the first term on the right, we first take gradient of the  $v$-equation in \eqref{main equation} and multiply it by $\nabla v|\nabla v|^2$, and then integrate by parts over $\Omega$, further note the identity $2\nabla v\cdot\nabla\Delta v=\Delta|\nabla v|^2-2|D^2v|^2$ to find that
	\be\label{grad v-l4-1}
	\begin{split}
		&\frac{d}{dt}\int_\Omega|\nabla v|^4+4\int_\Omega|\nabla v|^2|D^2v|^2+4\int_\Omega|\nabla v|^4+2\int_\Omega|\nabla|\nabla v|^2|^2\\
		&=-4\int_\Omega u\Delta v|\nabla v|^2-4\int_\Omega u\nabla v\cdot\nabla|\nabla v|^2+2\int_{\partial\Omega}|\nabla v|^2\frac{\partial}{\partial\nu}|\nabla v|^2.
	\end{split}
	\ee
	Again using  Young's inequality and the fact that $|\Delta v|^2\leq2|D^2v|^2$, we deduce
	\be\label{grad v-l4-2}
	\begin{split}
		&-4\int_\Omega u\Delta v|\nabla v|^2-4\int_\Omega u\nabla v\cdot\nabla|\nabla v|^2	\\
		&\leq2\int_\Omega|\nabla v|^2|\Delta v|^2+2\int_\Omega u^2|\nabla v|^2+\int_\Omega|\nabla|\nabla v|^2|^2+4\int_\Omega u^2|\nabla v|^2\\
		&\leq4\int_\Omega|\nabla v|^2|D^2v|^2+\int_\Omega|\nabla|\nabla v|^2|^2+6\int_\Omega u^2|\nabla v|^2.
	\end{split}
	\ee
	Given the boundedness of $\|\nabla v\|_{L^2}$  in  \eqref{ulnu-grad v-delta v-estimate},   the boundary integral in \eqref{grad v-l4-1} is well-handed  by  Sobolev boundary trace embedding as follows (\cite{Jin-Xiang21, Xiang18JMP, Xiang18JMAA, Xiang-Zheng19Non}):
	\be\label{grad v-l4-3}
	\begin{split}
		2\int_{\partial\Omega}|\nabla v|^2\frac{\partial|\nabla v|^2}{\partial\nu}&\leq \frac{1}{2}\int_\Omega|\nabla|\nabla v|^2|^2+C_{13}\left(\int_\Omega|\nabla v|^2\right)^2	\\
		&\leq \frac{1}{2}\int_\Omega|\nabla|\nabla v|^2|^2+C_{14}.
	\end{split}			
	\ee
	Plugging \eqref{grad v-l4-2} and \eqref{grad v-l4-3} into \eqref{grad v-l4-1}, we obtain that
$$
		 \frac{d}{dt}\int_\Omega|\nabla v|^4+4\int_\Omega|\nabla v|^4+\frac{1}{2}\int_\Omega|\nabla|\nabla v|^2|^2\leq6\int_\Omega u^2|\nabla v|^2+C_{14},
$$
which along with \eqref{u-l2-4} show directly that
\be\label{u-l2+gradv-l4}
\begin{split}
	&\frac{d}{dt}\int_\Omega \left(u^2+|\nabla v|^4\right)+\int_\Omega \left(u^2+|\nabla v|^4\right)+\int_\Omega|\nabla u|^2+\frac{1}{2}\int_\Omega|\nabla|\nabla v|^2|^2\\
&\leq2(3+\chi^2)\int_\Omega u^2|\nabla v|^2+\int_\Omega u^3+C_{12}+C_{14}.
\end{split}
\ee
Now, to estimate the term $\int_\Omega u^2|\nabla v|^2$, we apply again  Young's inequality to obtain, for any $\eta_1>0$, that
	\be\label{u2-grad v2-1}
		\int_\Omega u^2|\nabla v|^2\leq\eta_1\int_\Omega|\nabla v|^6+\frac{2}{3(3\eta_1)^\frac{1}{2}}\int_\Omega u^3.
	\ee
	Then the 2D  Gagliardo-Nirenberg inequality and boundedness of $\|\nabla v\|_{L^2}^2$  show that
	\be\label{grad v-l6}
	\begin{split}
		\int_\Omega|\nabla v|^6=\||\nabla v|^2\|_{L^3}^3&\leq C_{15}\|\nabla|\nabla v|^2\|_{L^2}^2\||\nabla v|^2\|_{L^1}+C_{15}\||\nabla v|^2\|_{L^1}^3	\\
		&\leq C_{16}\int_\Omega|\nabla|\nabla v|^2|^2+C_{16}.
	\end{split}
	\ee
	As for $\int_\Omega u^3$, based on the boundedness  $\|u\|_{L^1}+\|u\ln u\|_{L^1}$ guaranteed  in
 \eqref{v-l1-estimate} and \eqref{ulnu-grad v-delta v-estimate},  for any $\eta_2>0$, we infer from the 2D Gagliardo-Nirenberg inequality involving logarithmic function in \eqref{GN-log} \cite{Tao-Winkler14JDE, Xiang-Zheng19Non}   that
	\be\label{u-l3}
	\begin{split}
		\int_\Omega u^3=\|u\|_{L^3}^3&\leq\eta_2\|\nabla u\|_{L^2}^2\|u\ln u\|_{L^1}+C_{17}\|u\|_{L^1}^3+C_{\eta_2}	\\
		&\leq\eta_2C_{18}\int_\Omega |\nabla u|^2+C_{\eta_2}.
	\end{split}
	\ee
	Then we conclude from \eqref{u2-grad v2-1}, \eqref{grad v-l6} and \eqref{u-l3} that
	\be\label{u2-grad v2-2}
		\int_\Omega u^2|\nabla v|^2\leq\frac{2\eta_2C_{18}}{3(3\eta_1)^\frac{1}{2}}\int_\Omega |\nabla u|^2+\eta_1C_{16}\int_\Omega|\nabla|\nabla v|^2|^2+C_{\eta_1\eta_2}.
	\ee
	Now, fixing $\eta_6,\eta_7>0$ in \eqref{u2-grad v2-2} sufficiently small such that
	$$
		2\eta_1C_{16}(3+\chi^2)=\frac{1}{2},\quad  \quad  \eta_2C_{18}\left(1+\frac{4(3+\chi^2)}{3(3\eta_1)^\frac{1}{2}}\right)=1,
	$$
then we substitute   \eqref{u2-grad v2-2} and \eqref{u-l3} into \eqref{u-l2+gradv-l4} to  derive that
	\be\label{u-l2-grad v-l4-grad w-l6-0}
		\frac{d}{dt}\int_\Omega\left(u^2+|\nabla v|^4\right)+\int_\Omega\left(u^2+|\nabla v|^4\right)\leq C_{19}.
	\nonumber
	\ee
	This upon being multiplied by $e^t$ and then being integrated, yields quickly
	\be\label{u-l2-grad v-l4-grad w-l6-1}
		\int_\Omega\left(u^2+ |\nabla v|^4\right)(t)\\
		\leq \int_\Omega\left(u_0^2+ |\nabla v_0|^4\right)+C_{19}, \quad\forall t\in (0, T_m).
	\nonumber
	\ee
	Applying  Lemma \ref{v-lq-norm-lemma} with $p=2$ and $n=2$ and noticing $ v_*>0, 0\leq w\leq \|w_0\|_{L^\infty}$ and  the $L^p$-bound of $\nabla w$ by $\nabla v$ in \eqref{w-lp}, we further have, for any $q\in[1,\infty)$, that
	\be\label{v-l1q-w-l1q-norm}
		\|v(t)\|_{W^{1,q}}+\|w(t)\|_{W^{1,q}}\leq C_{20}(q),\quad\forall t\in (0, T_m).
	\ee
	
	\textbf{Step 2:} $L^p$-boundedness of $u$. For any $p>\max\{1,1-\theta\} $, we multiply both sides of the $u$-equation in \eqref{main equation} by $u^{p-1}$ and then integrate by parts over $\Omega$ to derive that
	\be\label{u-lp}
	\begin{split}
		&\frac{1}{p(p-1)}\frac{d}{dt}\int_\Omega u^p+\frac{4}{p^2}\int_\Omega |\nabla u^{\frac{p}{2}}|^2\\
&=\chi\int_\Omega u^{p-1}\nabla u\cdot\nabla v+\xi\int_\Omega S(u)u^{p-2}\nabla u\cdot\nabla w	\\
		&=-\chi\int_\Omega u^{p-1}\nabla u\cdot\nabla v-\xi\int_\Omega\Phi_3(u)\Delta w,
	\end{split}			
	\ee
	where, as before in \eqref{phi2}, the growth condition $0\leq S(u)\leq u(u+1)^{\theta-1}$ shows
	\be\label{phi3} 		
		0\leq \Phi_3(u):=\int_{0}^{u}\sigma^{p-2}S(\sigma)d\sigma\leq\frac{2^{p+\theta-2}}{p+\theta-1}(u^{p+\theta-1}+1).
	\ee
	
	For the first term on the right hand side of \eqref{u-lp}, by Young's inequality, H\"{o}lder's  inequality and Gagliardo-Nirenberg inequality in Lemma \ref{GN-lemma} and note the mass identity $\|u\|_{L^1}=\|u_0\|_{L^1}$ and  $W^{1,4}$-boundedness of $v$ in \eqref{v-l1q-w-l1q-norm}, we have that
	\be\label{2term}
	\begin{split}
		&-\chi\int_\Omega u^{p-1}\nabla u\cdot\nabla v\\
&\leq\frac{1}{4}\int_\Omega u^{p-2}|\nabla u|^2+\chi^2\int_\Omega u^p|\nabla v|^2	\\
		&\ \leq\frac{1}{p^2}\int_\Omega |\nabla u^{\frac{p}{2}}|^2+\chi^2\|u^p\|_{L^2}\|\nabla v\|_{L^4}^2	\\
		&\ \leq\frac{1}{p^2}\int_\Omega |\nabla u^{\frac{p}{2}}|^2+C_{21}\|u^{\frac{p}{2}}\|_{L^4}^2	\\
		&\ \leq\frac{1}{p^2}\int_\Omega |\nabla u^\frac{p}{2}|^2+C_{22}\left(\|\nabla u^\frac{p}{2}\|_{L^2}^{\frac{2p-1}{p}}\|u^\frac{p}{2}\|_{L^\frac{2}{p}}^\frac{1}{p}+ \|u^\frac{p}{2}\|_{L^\frac{2}{p}}^2\right)	\\
		&\ \leq\frac{2}{p^2}\int_\Omega |\nabla u^\frac{p}{2}|^2+C_{23}.
	\end{split}
	\ee
	Proceeding in a similar way leading to \eqref{u-l2-3},  for $\theta<2$ with $\nabla\sqrt{w_0}\in L^\infty(\Omega)$, we utilize the estimate \eqref{phi3},   the upper bound of $-\Delta w$ in \eqref{delta-w-estimate} with $\tau=1$ and the $(L^2, L^q)$-boundedness of $(u,v)$ for any finite $q$ to estimate the second term on the right hand side of \eqref{u-lp} as follows:
	\be\label{3term}
	\begin{split}
		-\xi\int_\Omega\Phi_3(u)\Delta w&\leq\frac{2^{p+\theta-2}\xi}{p+\theta-1}\int_\Omega (u^{p+\theta-1}+1)(\|w_0\|_{L^\infty}\cdot v+C_{24})	\\
		&\leq C_{25}\int_\Omega u^{p+\theta-1}v+C_{25}\int_\Omega u^{p+\theta-1}+C_{25}\\
		&\leq\int_\Omega u^{p+1}+C_{26}\\
        &\leq C_{27}\left(\|\nabla u^\frac{p}{2}\|_{L^2}^\frac{2(p-1)}{p}\|u^\frac{p}{2}\|_{L^\frac{4}{p}}^\frac{4}{p}
        +\|u^\frac{p}{2}\|_{L^\frac{4}{p}}^\frac{2(p+1)}{p}\right)\\
        &\ \leq\frac{1}{p^2}\int_\Omega |\nabla u^\frac{p}{2}|^2+C_{28}.
	\end{split}
	\ee
	In the case that $\theta<\frac{3}{2}$ while without  further assumption on $w_0$,   we apply the bound of $\Phi_3$ in \eqref{phi3} and the general point-wise upper bound of $-\Delta w$ in \eqref{delta-w-estimate0} to derive the same type  estimate \eqref{3term} by adapting the arguments in \eqref{u-l2-5}-\eqref{u-l2-6}.
	
	A use of the GN inequality and Young's inequality simply show that
	\be\label{4term}
		\frac{1}{p(p-1)}\int_\Omega u^p\leq\frac{1}{p^2}\int_\Omega|\nabla u^\frac{p}{2}|^2+C_{29}.
	\ee
	Substituting \eqref{2term}, \eqref{3term} and \eqref{4term} into \eqref{u-lp}, we obtain the following key ODI:
	$$
		\frac{d}{dt}\int_\Omega u^p+\int_\Omega u^p\leq C_{30},\quad\forall t\in (0, T_m),
	$$
	which quickly shows the $L^p$-boundedness of $u$.

	\textbf{Step 3:} With this arbitrary global $L^p$-boundedness of $u$ at hand, we conclude from Lemma \ref{v-lq-norm-lemma} that $v\in L^\infty(0, T_m; W^{1,\infty}(\Omega))$, which along with the facts that $ v_*>0, 0\leq w\leq \|w_0\|_{L^\infty}$ and  the $L^\infty$-bound of $w$ by $v$ in \eqref{w-lp} shows that $w\in L^\infty(0, T_m; W^{1,\infty}(\Omega))$.  Then, it is quite standard to derive the global existence and boundedness via either Moser iteration or Neumann heat semigroup estimates, cf. \cite{Jin-Xiang21,Tao-Winkler14JDE, Xiang15, Xiang-Zheng19Non, Xiang22,Zheng-Ke19}. Here, due to the homogeneous boundary condition for $u$ in \eqref{nbc}, we use the later to give a short proof for completeness. Indeed, by the variation-of-constants formula for $u$, we have
\be\label{u-exp}
u(t)=e^{t\Delta}u_0+\chi\int_0^te^{(t-s)\Delta}\nabla\cdot(u\nabla v)(s)ds-\xi\int_0^te^{(t-s)\Delta}\nabla\cdot(S(u)\nabla w)(s)ds.
\ee
Then applying the well-known $L^p$-$L^q$ smoothing estimates of the Neumann heat semigroup $\{e^{t\Delta }\}_{t\geq0}$ (\cite{Winkler10JDE, Cao15}), we use the derived  $(L^p, L^\infty, L^\infty)$-boundedness of $(u,\nabla v, \nabla w)$ and the growth condition $0\leq S(u)\leq u(u+1)^{\theta-1}$ to infer that
\begin{equation*}
\begin{split}
\|u(t)\|_{L^\infty}&\leq \|e^{t\Delta }u_0\|_{L^\infty}+ \chi\int_0^t\left\|e^{(t-s)\Delta }\nabla \cdot((u\nabla v)(s))\right\|_{L^\infty}ds\\
&\ \ \ +\xi\int_0^t\left\|e^{(t-s)\Delta }\nabla \cdot((S(u)\nabla w)(s))\right\|_{L^\infty}ds\\
  &\leq  \|u_0\|_{L^\infty}+C_{31}\chi\int_0^t\left(1+(t-s)^{-\frac{1}{2}-\frac{1}{4}}\right)e^{-\lambda_1(t-s) }\left\|(u\nabla v)(s)\right\|_{L^{2n}}ds\\
  &\ \ \ +C_{31}\xi\int_0^t\left(1+(t-s)^{-\frac{1}{2}-\frac{1}{4}}\right)e^{-\lambda_1(t-s) }\left\|(S(u)\nabla w)(s)\right\|_{L^{2n}}ds\\
 &\leq  \|u_0\|_{L^\infty}+C_{32}\chi\int_0^t\left(1+(t-s)^{-\frac{3}{4}}\right)e^{-\lambda_1(t-s) }\left\|u(s)\right\|_{L^{2n}}\left\|\nabla v(s)\right\|_{L^\infty}ds\\
 &\ \ \ +C_{32}\xi\int_0^t\left(1+(t-s)^{-\frac{3}{4}}\right)e^{-\lambda_1(t-s) }\left\|1+u(s)\right\|_{L^{2n\max\{1,\theta\}}}^{\max\{1,\theta\}}\left\|\nabla w(s)\right\|_{L^\infty}ds\\
  &\leq \|u_0\|_{L^\infty}+C_{33}(\chi+\xi)\left[\int_0^1\left(1+z^{-\frac{3}{4}}\right)dz
  +2\int_1^\infty
  e^{-\lambda_1z }dz\right]\\
   &\leq \|u_0\|_{L^\infty}+C_{34}(\chi+\xi)\left(5+\frac{2}{\lambda_1}\right).
\end{split}
\end{equation*}
 Therefore, to sum up, we have shown that
 $$
 \sup_{t\in(0, T_m)}\left(\|u(t)\|_{L^\infty}+\|v(t)\|_{W^{1,\infty}}+\|w(t)\|_{W^{1,\infty}}\right)<\infty,
 $$
 and then, by the extensibility criterion \eqref{extensibility-criterion}, we see  that $T_m=\infty$, that is, the classical solution $(u,v,w)$ of \eqref{main equation} exists globally in time and is uniformly bounded in the sense of \eqref{global-boundedness} .	
\end{proof}
	
\subsection{Global existence and boundedness in nD $(n\geq2)$ for $\tau=0$}
	
	In this subsection, we study the global existence and boundedness of solutions to  \eqref{main equation} with $\tau=0$. We shall see that chemo-repulsion still prevents blow-up of classical solutions to the combined negative chemotaxis-haptotaxis model \eqref{main equation} with $\tau=0$ and suitably strong hapto-sensitivity $S$, thus proving Theorem \ref{theorem2}.
	
\begin{lemma}\label{u-lp2-lemma}
	Let $n\geq 2$, $\tau=0$ and $\chi>0$. Assume that $\theta<2$ with $\nabla\sqrt{w_0}\in L^\infty(\Omega)$ or $\theta<\frac{3}{2}$ without further condition on $w_0$. Then the solution of  \eqref{main equation}  exists globally is uniformly bounded in the sense of \eqref{global-boundedness}.
\end{lemma}
\begin{proof}
	Noticing $\Delta v=v-u$ due to $\tau=0$ in \eqref{main equation},  we  compute from \eqref{u-lp}   that
	\be\label{u-lp2}
	\begin{split}
		&\frac{1}{p(p-1)}\frac{d}{dt}\int_\Omega u^p+\frac{4}{p^2}\int_\Omega |\nabla u^{\frac{p}{2}}|^2\\
		&=-\frac{\chi}{p}\int_\Omega \nabla u^p\cdot\nabla v-\xi\int_\Omega\Phi_3(u)\Delta w	\\
		&=\frac{\chi}{p}\int_\Omega u^p\Delta v-\xi\int_\Omega\Phi_3(u)\Delta w	\\
		&=\frac{\chi}{p}\int_\Omega u^pv-\frac{\chi}{p}\int_\Omega u^{p+1}-\xi\int_\Omega\Phi_3(u)\Delta w,
	\end{split}			
	\ee
	where $\Phi_3$ is prescribed as in \eqref{phi3}.
	
	For the first two terms, we  use  Young's inequality and the elliptic estimate \eqref{v-lp+1} along with the mass conservation $\|u\|_{L^1}=\|u_0\|_{L^1}$ to bound that
	\be\label{u-lp-v}
	\begin{split}
		\frac{\chi}{p}\int_\Omega u^pv-\frac{\chi}{p}\int_\Omega u^{p+1}&\leq-\frac{3\chi}{4p}\int_\Omega u^{p+1}+C_1\int_\Omega v^{p+1}\\
		&\leq-\frac{\chi}{2p}\int_\Omega u^{p+1}+C_2.
	\end{split}
	\ee

\textbf{Case 1:} $\theta<2$ and $\nabla\sqrt{w_0}\in L^\infty(\Omega)$. In this case, by the estimate of $\Phi_3$ in  \eqref{phi3}	and the pointwise upper bound of $-\Delta w$ in \eqref{delta-w-estimate} with $\tau=0$, we obtain that
\be\label{hap-taxis1}
		-\xi\int_\Omega\Phi_3(u)\Delta w\leq C_3\int_\Omega (u^{p+\theta-1}+1)\leq\frac{\chi}{2p}\int_\Omega u^{p+1}+C_{4}.
	\ee

\textbf{Case 2:} $\theta<\frac{3}{2}$ and without further condition on $w_0$. We then  employ the bound of $\Phi_3$ in \eqref{phi3} and the general point-wise upper bound of $-\Delta w$ in \eqref{delta-w-estimate0} with $\tau=0$ to derive  that
	\be\label{u-lp-1}
	\begin{split}
		-\xi\int_\Omega\Phi_3(u)\Delta w&\leq\frac{2^{p+\theta-2}\xi}{p+\theta-1}\left(\frac{\|w_0\|_{L^\infty}}{e}+\|\Delta w_0\|_{L^\infty}\right)\int_\Omega(u^{p+\theta-1}+1)	\\
		&\ \ \ -\frac{2^{p+\theta-1}\xi}{p+\theta-1}\int_\Omega(u^{p+\theta-1}+1)\nabla w_0\cdot\nabla e^{-\int_{0}^{t}v(s)ds}.
	\end{split}
	\ee
	Applying Young's inequality infers, for $\theta<2$ and any $\eta_1>0$, that
	\be\label{u-lp-2}
		\int_\Omega u^{p+\theta-1}\leq\eta_1\int_\Omega u^{p+1}+C_{\eta_1}.
	\ee	
	In light of the Young's inequality and the nonnegativity of $v$, we deduce for $\theta<\frac{3}{2}$ and $\eta_2,\eta_3>0$, that
	\be\label{u-lp-3}
	\begin{split}
		&-\int_\Omega(u^{p+\theta-1}+1)\nabla w_0\cdot\nabla e^{-\int_{0}^{t}v(s)ds}\\
		&=(p+\theta-1)\int_\Omega u^{p+\theta-2} e^{-\int_{0}^{t}v(s)ds}\nabla u\cdot\nabla w_0	\\
		&\ \ \ +\int_\Omega(u^{p+\theta-1}+1)e^{-\int_{0}^{t}v(s)ds}\Delta w_0	\\
		&\leq(p+\theta-1)\|\nabla w_0\|_{L^\infty}\int_\Omega u^{p+\theta-2}|\nabla u|+\|\Delta w_0\|_{L^\infty}\int_\Omega(u^{p+\theta-1}+1) \\
		&\leq\frac{\eta_2}{4}\int_\Omega u^{p-2}|\nabla u|^2+C_{\eta_2}\int_\Omega u^{p+2\theta-2}+\|\Delta w_0\|_{L^\infty}\int_\Omega(u^{p+\theta-1}+1)	\\
		&\leq\frac{\eta_2}{p^2}\int_\Omega |\nabla u^\frac{p}{2}|^2+\eta_3\int_\Omega u^{p+1}+C_{\eta_2\eta_3}.
	\end{split}
	\ee	
Now, substituting \eqref{u-lp-2} and \eqref{u-lp-3} into \eqref{u-lp-1} and fixing sufficiently small $\eta_1, \eta_2, \eta_3>0$, we obtain that
\be\label{hap-taxis-est}
		-\xi\int_\Omega\Phi_3(u)\Delta w\leq\frac{1}{p^2}\int_\Omega |\nabla u^\frac{p}{2}|^2+\frac{\chi}{2p}\int_\Omega u^{p+1}+C_5.
	\ee
Substituting  \eqref{u-lp-v}  and   \eqref{hap-taxis1} in \textbf{Case 1} or \eqref{hap-taxis-est} in \textbf{Case 2} into \eqref{u-lp2} and using   \eqref{4term}, we obtain an important ODI for $\|u|_{L^p}$ as follows:
$$
		\frac{d}{dt}\int_\Omega u^p+\int_\Omega u^p\leq C_6, \quad\forall t\in(0,T_m),
$$
which upon being solved directly yields that
	\be\label{u-lp-bound}
		\int_\Omega u^p(t)\leq C_7(p),\quad\forall t\in(0,T_m).
	\ee
Based on the arbitrary  $L^p$-boundedness of $u$ in \eqref{u-lp-bound}, one can readily  obtain the desired global existence  and boundedness as illustrated in  \textbf{Step 3} in Lemma \ref{global-boundedness-lemma}.
\end{proof}

\subsection{Haptotaxis-only ($\chi=0$) can not induce any blow-up}

	In the previous subsections, we have seen  that the fact $\chi>0$ does have certain positive effect toward boundedness. How about the case  $\chi=0$? In this subsection, we  explore the effect of haptotaxis   on global solvability  for \eqref{main equation} with $\chi=0$. Under certain growth condition on $S$, we shall illustrate  that haptotaxis  can not induce any blow-up by guaranteeing all solutions are globally bounded, showing    Theorem \ref{theorem1} (B2). This extends \cite[Theorem 1.1]{Xiang22} where $S(u)=u$. The following lemma also plays a role in studying the long time dynamics of globally bounded solutions later on.

\begin{lemma}\label{haptotaxis-global-bound-lemma}
	For $n\geq 2, \tau=0$, assume $\theta<1+\frac{2}{n}$ with
	$\nabla\sqrt{w_0}\in L^\infty (\Omega)$ or  $\theta<1+\frac{1}{n}$.  For $n=2, \tau=1$, assume $ \theta<2$ with $\nabla\sqrt{w_0}\in L^\infty (\Omega)$; for $n\geq 2, \tau=1$,  assume
\be\label{theta-con-hap}
\theta< \begin{cases}
\frac{3}{2}, &\text{ if } n=2, \\[0.2cm]
 \frac{4}{n}, &\text{ if } 3\leq n\leq 4, \\[0.2cm]
\frac{6n-n^2-4}{2(n-2)}, &\text{ if } n\geq 5.
\end{cases}
 \ee
  Then the local classical solution $(u,v,w)$ of the system \eqref{main equation} with $\chi=0$ exists globally in time and remains uniformly bounded in the sense of \eqref{global-boundedness}.
\end{lemma}
\begin{proof}
 We only illustrate the proof for the hard case that $\tau=1$  and  $\theta$ satisfies \eqref{theta-con-hap} without  further condition on $w_0$. 	Then, for our later purpose, we first choose
\be\label{p-choice}
\begin{cases}
\max\left\{2-\theta, \   \frac{n}{2}, \  \frac{n}{2}\theta, \    n(\theta-1)\right\}<p,  \ &\text{ if } 2\leq n\leq 4,\\[0.25cm]
\max\left\{2-\theta, \   \frac{n}{2}, \  \frac{n}{2}\theta\right\}<p<\frac{(n-2)(1-\theta)}{n-4},&\text{ if } n\geq  5,
\end{cases}
\ee
then we deduce  from \eqref{u-lp} with $\chi=0$ that
	\be\label{hap1}
		\frac{d}{dt}\int_\Omega u^p+\frac{4(p-1)}{p}\int_\Omega|\nabla u^\frac{p}{2}|^2=-\xi\int_\Omega\Phi_3(u)\Delta w,\quad\forall t\in(0,T_m), 	
	\ee
where $\Phi_3$ is given by \eqref{phi3}.

 In light of the general upper bound of $-\Delta w$ in \eqref{delta-w-estimate0}, the property of $\Phi_3$ in \eqref{phi3}, we derive that
	\be\label{hap2}
	\begin{split}
		-\xi\int_\Omega\Phi_3(u)\Delta w&\leq\frac{2^{p+\theta-2}\tau\xi\|w_0\|_{L^\infty}}{p+\theta-1}\int_\Omega \left(1+u^{p+\theta-1}\right)v		\\
		&\ \ \ -\frac{2^{p+\theta-1}\xi}{p+\theta-1}\int_\Omega \left(1+u^{p+\theta-1}\right)\nabla w_0\cdot\nabla e^{-\int_{0}^{t}v(s)ds}\\
		&\ \ \ +\frac{2^{p+\theta-2}\xi}{p+\theta-1}\left(\frac{\|w_0\|_{L^\infty}}{e}+\|\Delta w_0\|_{L^\infty}\right)\int_\Omega \left(1+u^{p+\theta-1}\right).
	\end{split}
	\ee
 Recalling the mass identity $\|u\|_{L^1}=\|u_0\|_{L^1}$ and the $L^q$-boundedness of $v$ with $q\in[1,\frac{n}{n-2})$ in Lemma \ref{v-lq-norm-lemma}, by the regulation of $\theta$ in  \eqref{theta-con-hap} and the choice of $p$ in \eqref{p-choice}, we can first fix $q$ to fulfill
	\be\label{q-choice}
	 \max\left\{1, \ \frac{n}{n+2-n\theta},  \  \frac{np}{2p-(n-2)(\theta-1)}\right\}<q<\frac{n}{n-2},
	\ee
	then, thanks to the requirement of $\theta$ in \eqref{theta-con-hap} and the choices of $p,q$ in \eqref{p-choice} and \eqref{q-choice}, we  invoke H\"{o}lder's inequality, Gagliardo-Nirenberg inequality and Young's inequality to obtain, for any $\eta_1>0$, that
	\be\label{term1}
	\begin{split}
		&\int_\Omega u^{p+\theta-1}v\\
&\leq\|u^\frac{p}{2}\|_{L^\frac{2(p+\theta-1)q}{p(q-1)}}^\frac{2(p+\theta-1)}{p}\|v\|_{L^q}	\\
		&\leq C_2\left(\|\nabla u^\frac{p}{2}\|_{L^2}^{\frac{2n(p+\theta-2+q^{-1})}{np-n+2}}
\|u^\frac{p}{2}\|_{L^\frac{2}{p}}^{\frac{2[2p-npq^{-1}-(n-2)(\theta-1)]}{(np-n+2)p}}+\|u^\frac{p}{2}\|_{L^\frac{2}{p}}^\frac{2(p+\theta-1)}{p}\right)	\\
		&\leq C_3\|\nabla u^\frac{p}{2}\|_{L^2}^{\frac{2n(p+\theta-2+q^{-1})}{np-n+2}}+C_3	\\
		&\leq\eta_1\int_\Omega|\nabla u^\frac{p}{2}|^2+C_{\eta_1},
	\end{split}
	\ee	
	where the choice of $q$ in \eqref{q-choice} was used  to ensure that
	$$
		\frac{2n(p+\theta-2+q^{-1})}{np-n+2}<2.
	$$
For the second term on the right-hand side of \eqref{hap2}, as in \eqref{u-lp-3}, we employ  integration by parts, the Cauchy-Schwarz inequality, the nonnegativity of $v$ and the fact $\theta<1+\frac{1}{n}$ to deduce, for any $\eta_2, \eta_3>0$, that
	\be\label{term2}
	\begin{split}
		&-\int_\Omega (1+u^{p+\theta-1})\nabla w_0\cdot\nabla e^{-\int_{0}^{t}v(s)ds}\\
		&=(p+\theta-1)\int_\Omega u^{p+\theta-2}e^{-\int_{0}^{t}v(s)ds}\nabla u\cdot\nabla w_0\\
&\ \ +\int_\Omega (1+u^{p+\theta-1})e^{-\int_{0}^{t}v(s)ds}\Delta w_0\\
		&\leq(p+\theta-1)\|\nabla w_0\|_{L^\infty}\int_\Omega  u^{p+\theta-2}|\nabla u|+\|\Delta w_0\|_{L^\infty}\int_\Omega (1+ u^{p+\theta-1})	\\
		&\leq\eta_2\int_\Omega|\nabla u^\frac{p}{2}|^2+C_{\eta_2}\int_\Omega u^{p+2\theta-2}+\|\Delta w_0\|_{L^\infty}\int_\Omega (1+u^{p+\theta-1})\\
&\leq\eta_2\int_\Omega|\nabla u^\frac{p}{2}|^2+ \eta_3\int_\Omega u^{p+\frac{2}{n}}+C_{\eta_2\eta_3}.
	\end{split}
	\ee
	Using the Gagliardo-Nirenberg inequality  and  $\|u\|_{L^1}=\|u_0\|_{L^1}$, we have
	\be\label{u0.5p}
\begin{split}
		\int_\Omega u^{p+\frac{2}{n}}=\|u^\frac{p}{2}\|_{L^\frac{2(2+np)}{np}}^\frac{2(2+np)}{np}&\leq C_4\|\nabla u^\frac{p}{2}\|_{L^2}^2\|u^\frac{p}{2}\|_{L^\frac{2}{p}}^\frac{4}{np}
+C_4\|u^\frac{p}{2}\|_{L^\frac{2}{p}}^\frac{2(2+np)}{np}\\
&\leq C_5\|\nabla u^\frac{p}{2}\|_{L^2}^2+C_5.
\end{split}
	\ee
	Choosing $\eta_i>0 (i=1,2,3)$ sufficiently small and recalling \eqref{4term} implied by \eqref{u0.5p},  we can conclude from \eqref{hap1}, \eqref{hap2}, \eqref{term1}, \eqref{term2} and \eqref{u0.5p}, that
	$$
		\frac{d}{dt}\int_\Omega u^p+\int_\Omega u^p\leq C_6,\quad\forall t\in(0,T_m),
	$$
	which  directly gives rise to
	$$
		\|u(t)\|_{L^p}^p\leq\|u_0\|_{L^p}^p+C_6,\quad\forall t\in(0,T_m).
	$$
When $2\leq n\leq 4$, the exponent $p$ in \eqref{p-choice} can be arbitrarily large; with such large $L^p$-uniform boundedness of $u$, one can readily  derive the desired global existence  and boundedness  as illustrated in  \textbf{Step 3} in Lemma \ref{global-boundedness-lemma}.	

When $n\geq 5$, noticing from \eqref{p-choice} that $p>\frac{n}{2}$, we can choose $q$ and $q^\prime$ to satisfy
$$
n<q^\prime<\frac{np}{(n-p)^+}, \ \ \  q>\max\left\{q^\prime, \ \frac{npq^\prime}{np-(n-p)^+q^\prime}\right\},
$$
which ensures that
\be\label{qqprime}
q>q^\prime>n, \ \ \ \frac{qq^\prime}{q-q^\prime}<\frac{np}{(n-p)^+}.
\ee
Then applying  Lemmas \ref{v-lq-norm-lemma}, \ref{v-lower-bound-lemma} and \ref{w-lemma}, we infer that
\be\label{v-l1q-w-l1q-norm-hap}
\|v(t)\|_{W^{1,\frac{qq^\prime}{q-q^\prime}}}+\|w(t)\|_{W^{1,\frac{qq^\prime}{q-q^\prime}}}\leq C_7,\quad\forall t\in (0, T_m).
\ee
Thus, the growth condition  $0\leq S(u)\leq u(u+1)^{\theta-1}$ implies  that
\be\label{su-bdd}
\|S(u)\|_{L^q}\leq C_8\left(\|u\|_{L^\infty}^\frac{(q\theta-p)^+}{q}\|u\|_{L^p}^\frac{p}{q}+1\right)\leq C_9\left(\|u\|_{L^\infty}^\frac{(q\theta-p)^+}{q}+1\right).
\ee
 Applying Neumann semigroup estimates and H\"{o}lder inequality, we deduce from \eqref{qqprime}, \eqref{v-l1q-w-l1q-norm-hap} and \eqref{su-bdd}  that
 \be\label{hap-est-hard}
\begin{split}
&\xi\int_0^t\left\|e^{(t-s)\Delta }\nabla \cdot((S(u)\nabla w)(s))\right\|_{L^\infty}ds\\
  &\leq  C_{10}\int_0^t\left(1+(t-s)^{-\frac{1}{2}-\frac{n}{2q^\prime}}\right)e^{-\lambda_1(t-s) }\left\|(S(u)\nabla w)(s)\right\|_{L^{q^\prime}}ds\\
 &\leq  C_{10}\int_0^t\left(1+(t-s)^{-\frac{1}{2}-\frac{n}{2q^\prime}}\right)e^{-\lambda_1(t-s) }\left\|S(u)(s)\right\|_{L^q}\left\|\nabla w(s)\right\|_{L^\frac{qq^\prime}{q-q^\prime}}ds\\
  &\leq C_{11}\sup_{s\in(0,t)}\|u(s)\|_{L^\infty}^\frac{(q\theta-p)^+}{q}+C_{11}\\
  &\leq \frac{1}{2}\sup_{s\in(0,t)}\|u(s)\|_{L^\infty}+C_{12},
\end{split}
\ee
where the fact $\theta<1$ due to  \eqref{theta-con-hap} with $n\geq 5$ was used in the last inequality.

Finally, taking $L^\infty$-norm on both sides of  \eqref{u-exp} with $\chi=0$ and using \eqref{hap-est-hard}, we quickly obtain the $L^\infty$-boundedness of $u$, and then,  from Lemmas \ref{v-lq-norm-lemma}, \ref{v-lower-bound-lemma} and \ref{w-lemma}, we derive the $W^{1,\infty}$-boundedness of $v$ and $w$.

Then, as before, we first conclude  global existence and then boundedness  of $(u,v,w)$  in the sense of \eqref{global-boundedness}.  	
\end{proof}

\section{long time behaviors of globally bounded solutions}

	With the global boundedness of classical solutions at hand, we next investigate   large time behaviors of such solutions. We show that any globally bounded solution of the IBVP \eqref{main equation}  converges exponentially to its constant steady state $(\bar{u}_0,\bar{u}_0,0)$,  thus proving  Theorem \ref{theorem1} (B3). To this end, we first establish the   $W^{1,\infty}$-exponential decay of $w$-component via the estimate \eqref{w-lp} and the lower bound  of $v$.

\begin{lemma}\label{w-convergence-lemma}
	Let  $n\geq 2$, $\tau\in\{0,1\}$ and let  $(u,v,w)$ be a globally bounded solution of the IBVP \eqref{main equation} in the topology $(L^\infty,W^{1,\infty},W^{1,\infty})$. Then the $w$-component decays exponentially: there exists $C>0$ such that, for $0<\gamma<v_*$,
	\be\label{w-decay}
		\|w(t)\|_{W^{1,\infty}}\leq Ce^{-\gamma t},\quad\forall t\geq0.
	\ee
\end{lemma}
\begin{proof}
	Noticing from \eqref{v-lower-bound-estimate} that $v\geq v_*>0$ on $\Omega\times(0,T_m)$, we quickly  derive from the expression of $w$ in \eqref{w-expression} that
	\be\label{w-1}
		w(x,t)=w_0e^{-\int_0^tv(x,s)ds}\leq \|w_0\|_{L^\infty}e^{-v_* t}, \ \  (x,t)\in \Omega\times(0,\infty).
	\ee
	Taking $p=\infty$ in \eqref{w-lp},   for any $0<\gamma<v_*$, thanks to the global boundedness of $\|\nabla v(t)\|_{L^\infty}$, we see there exists  $C_1>0$ such that
	\be\label{w-2}
	\begin{split}
		\|\nabla w(t)\|_{L^\infty}&\leq\|\nabla w_0\|_{L^\infty}e^{-v_*t}+\| w_0\|_{L^\infty}\sup_{s\in(0,\infty)}\|\nabla v(s)\|_{L^\infty}te^{-v_*t}	\\
		&\leq\|\nabla w_0\|_{L^\infty}e^{-\gamma t}+\| w_0\|_{L^\infty}\sup_{s\in(0,\infty)}\|\nabla v(s)\|_{L^\infty}te^{-(v_*-\gamma)t}e^{-\gamma t}	\\
		&\leq\left(\|\nabla w_0\|_{L^\infty}+\frac{\| w_0\|_{L^\infty}}{e(v_*-\gamma)}\sup_{s\in(0,\infty)}\|\nabla v(s)\|_{L^\infty}\right)e^{-\gamma t}	\\
		&\leq C_1e^{-\gamma t},\quad\forall t\geq0,
	\end{split}	
	\nonumber
	\ee
	which along with \eqref{w-1} quickly shows the $W^{1,\infty}$-exponential decay of $w$ in \eqref{w-decay}.
\end{proof}

With the establishment of exponential decay of $w$, we proceed to show that the global solution $(u,v,w)$ of the IBVP \eqref{main equation}  converges exponentially to its equilibrium state $(\bar{u}_0,\bar{u}_0,0)$. Our argument relies on  the following functional:
	\be\label{g-def}
		\mathcal{G}(t)=\int_\Omega u\ln\frac{u}{\bar{u}_0}+\frac{\tau\chi}{2}\int_\Omega|\nabla v|^2,
	\ee
which is modified from the functional $\mathcal{F}$ defined by \eqref{F-definition}.
\begin{lemma}\label{g-functional-lemma}
	Under Lemma \ref{w-convergence-lemma}, the time derivative of $\mathcal{G}$ defined by \eqref{g-def} decays exponentially according to
	\be\label{g-decay}
		0\leq\mathcal{G}(t)\leq Ce^{-\delta (t-1)},\quad\forall t\geq1,
	\ee
	where $C>0$ is a time-independent constant and  $\delta>0$  is defined by \eqref{delta-def}.
\end{lemma}	
\begin{proof} Since  $s\ln s\geq s-1$ for any $s>0$, it follows that
		$$
			\int_\Omega u\ln\frac{u}{\bar{u}_0}\geq \int_\Omega \left(u-\bar{u}_0\right)=0,
		$$
which gives the nonnegativity of  $\mathcal{G}$ defined by \eqref{g-def}.

	Multiplying the first equation in \eqref{main equation} by $\ln u-\ln\bar{u}_0$ and then integrating by parts over $\Omega$, we have that
	\be\label{g1}
		\frac{d}{dt}\int_\Omega u\ln\frac{u}{\bar{u}_0}+\int_\Omega\frac{|\nabla u|^2}{u}=-\chi\int_\Omega\nabla u\cdot\nabla v+\xi\int_\Omega
		\frac{S(u)}{u}\nabla u\cdot\nabla w.
	\ee
	We next compute by the second equation in \eqref{main equation} to derive
	\be\label{g2}
		\frac{\tau\chi}{2}\frac{d}{dt}\int_\Omega|\nabla v|^2+\chi\int_\Omega(|\Delta v|^2+|\nabla v|^2)=\chi\int_\Omega\nabla u\cdot\nabla v.
	\ee
	Adding \eqref{g1} and \eqref{g2} gives rise to
	\be\label{g3}
		\mathcal{G}^\prime(t)+\int_\Omega\frac{|\nabla u|^2}{u}+\chi\int_\Omega(|\Delta v|^2+|\nabla v|^2)=\xi\int_\Omega
		\frac{S(u)}{u}\nabla u\cdot\nabla w.
	\ee	
	Noticing that $0\leq S(u)\leq u(u+1)^{\theta-1}$,  the mass identity $\|u\|_{L^1}=\|u_0\|_{L^1}$, the $L^\infty$-boundedness of $u$ and the $W^{1,\infty}$-exponential decays of $w$ in \eqref{w-decay}, we estimate
	\be\label{g4}
	\begin{split}
		\xi\int_\Omega\frac{S(u)}{u}\nabla u\cdot\nabla w&\leq\xi\int_\Omega(u+1)^{\theta-1}|\nabla u\cdot\nabla w|	\\
		&\leq\frac{1}{2}\int_\Omega\frac{|\nabla u|^2}{u}+\frac{\xi^2}{2}\int_\Omega(u+1)^{2\theta-1}|\nabla w|^2	\\
		&\leq\frac{1}{2}\int_\Omega\frac{|\nabla u|^2}{u}+\frac{\xi^2}{2}(\|u\|_{L^\infty}+1)^{(2\theta-1)^+}|\Omega|\|\nabla w\|_{L^\infty}^2	\\
		&\leq\frac{1}{2}\int_\Omega\frac{|\nabla u|^2}{u}+C_1e^{-2\gamma t}.
	\end{split}		
	\ee
	Using the fact  $s\ln s\leq s-1+(s-1)^2$ for all $s>0$ and the Poincar\'{e}'s inequality, upon integrating and noticing the facts $\bar{u}=\bar{u}_0$ and  $u\in L^\infty(0, \infty; L^\infty(\Omega))$, we obtain that
	\be\label{g5}
	\begin{split}
		\int_\Omega \frac{u}{\bar{u}_0}\ln \frac{u}{\bar{u}_0}\leq\int_\Omega\left(\frac{u}{\bar{u}_0}-1\right)
+\int_\Omega\left(\frac{u}{\bar{u}_0}-1\right)^2\leq \frac{C_2}{|\bar{u}_0|^2}\int_\Omega|\nabla u|^2\leq\frac{C_3}{|\bar{u}_0|^2}\int_\Omega\frac{|\nabla u|^2}{u},
	\end{split}	
	\nonumber	
	\ee
	which entails (cf. also \eqref{sharp-ck+} in Lemma \ref{ext-ck})
	\be\label{g6}
		\int_\Omega u\ln\frac{u}{\bar u_0}\leq\frac{C_3}{\bar u_0}\int_\Omega\frac{|\nabla u|^2}{u}.
	\ee
	Substituting \eqref{g4}   into \eqref{g3} and using \eqref{g6} and \eqref{g-def}, we quickly derive  that
$$
\mathcal{G}^\prime(t)+\min\{\frac{\bar u_0}{2C_3}, \  \frac{2}{\tau}\}\mathcal{G}(t)\leq C_1e^{-2\gamma t},\quad\forall t\geq0,
$$
which further implies
\be\label{g7}
		\mathcal{G}^\prime(t)+\delta \mathcal{G}(t)\leq C_1e^{-2\gamma t},\quad\forall t\geq0,
	\ee
where the constant $\delta$ satisfies
	\be\label{delta-def}
		0<\delta<\min\left\{\frac{\bar u_0}{2C_3}, \  \frac{2}{\tau},\  2\gamma\right\}.
	\ee
	Multiplying  \eqref{g7} by $e^{\delta t}$,  integrating it on $[1,t)$ and using \eqref{delta-def}, we conclude that
	 $$
		\mathcal{G}(t)\leq \left(\mathcal{G}(1)+\frac{C_1}{2\gamma-\delta}\right)e^{-\delta (t-1)}, \quad \forall t\geq 1, $$
	trivially yielding the exponential decay of  $\mathcal{G}$ in \eqref{g-decay}.	
\end{proof}
	
	With this crucial exponential decay of $\mathcal{G}$ at hand, we are now in the position to deduce the exponential convergence of $(u,v)$-component  in the topology $(L^\infty,W^{1,\infty})$.
	
\begin{lemma}\label{uv-decay-lemma}
	Under Lemma \ref{w-convergence-lemma}, the $(u,v)$-component of the global bounded classical solution of \eqref{main equation} decays exponentially:  for $p>n$, there exists $C=C(p)>0$,
	\be\label{u-decay}
		\|u(t)-\bar{u}_0\|_{L^\infty}\leq Ce^{-\frac{\delta}{2p}t},\quad t>0,
	\ee
	and
	\be\label{v-decay}
		\|v(t)-\bar{u}_0\|_{W^{1,\infty}}\leq Ce^{-\min\{1,\frac{\delta}{2p}\}t},\quad t>0,
	\ee
	where $0<\gamma<v_*$ and $\delta$ is defined by \eqref{delta-def}.
\end{lemma}
\begin{proof}
	First, the exponential decay of $\mathcal{G}$ in \eqref{g-decay} along  with its definition in \eqref{g-def} quickly shows  that
	$$
		\int_\Omega u\ln\frac{u}{\bar{u}_0}\leq C_1e^{-\delta t},\quad\forall t\geq1.
	$$
	Thanks to the free conservation $\bar{u}=\bar{u}_0$, the widely used  Csisz\'{a}r-Kullback-Pinsker inequality (cf. \cite[Section 4]{Carrillo-et-al01})  implies that
	$$
		\|u-\bar{u}_0\|_{L^1}^2\leq2\bar{u}_0\int_\Omega u\ln\frac{u}{\bar{u}_0}\leq C_2e^{-\delta t},\quad\forall t\geq1,
	$$
	which together with the global $L^\infty$-boundedness of $u$ shows, for any $p>1$,  that
	\be\label{u-lp-decay} \|u-\bar{u}_0\|_{L^p}\leq\|u-\bar{u}_0\|_{L^\infty}^\frac{p-1}{p}\|u-\bar{u}_0\|_{L^1}^\frac{1}{p}\leq C_3e^{-\frac{\delta}{2p}t}, \quad\forall t\geq0.
	\ee
Since $\bar{u}=\bar{u}_0$, taking average of the $v$-equation in \eqref{main equation} gives
	$$
		\tau \bar{v}_t=-\bar{v}+\bar{u}_0,
	$$
	which entails
	\be\label{v-average}
		\bar{v}(t)=\tau (\bar{v}_0-\bar{u}_0)e^{-t}+\bar{u}_0.
	\ee
	Therefore,  the $v$-equation in \eqref{main equation} further shows that
	\be\label{v-bar-exp}
		\tau (v-\bar{v})_t=\Delta(v-\bar{v})-(v-\bar{v})+u-\bar{u}_0.
	\ee
\textbf{Case 1: $\tau=0$}. In this case,  for $p>n$,  by the $W^{2,p}$-elliptic estimate and the continuous embedding $W^{2,p}(\Omega)\hookrightarrow W^{1,\infty}(\Omega)$, it follows from \eqref{u-lp-decay} and \eqref{v-average}  that
\be\label{v-wlinfty-decay} \|v-\bar{u}_0\|_{W^{1,\infty}}\leq C_4 \|v-\bar{u}_0\|_{W^{2,p}}\leq C_5 \|u-\bar{u}_0\|_{L^p}\leq C_6e^{-\frac{\delta}{2p}t}, \quad\forall t\geq0.
	\ee
 \textbf{Case 2: $\tau=1$}. The variation-of-constants formula applied to \eqref{v-bar-exp} entails
	$$
		(v-\bar{v})(t)=e^{t(\Delta-1)}(v_0-\bar{v}_0)+\int_{0}^{t}e^{(t-s)(\Delta-1)}(u(s)-\bar{u}_0)ds.
	$$
	Then for $p>n$,  $\delta$ satisfying \eqref{delta-def} and $\delta<2n\lambda_1$, we use the well-known Neumann heat semigroup estimates and  \eqref{u-lp-decay} to infer
	\be\label{v-vbar-estimate}
	\begin{split}
		&\|(v-\bar{v})(t)\|_{W^{1,\infty}}	\\
		&\leq C_7e^{-(\lambda_1+1)t}+C_7\int_{0}^{t}\left(1+(t-s)^{-\frac{1}{2}-\frac{n}{2p}}\right)e^{-(\lambda_1+1)(t-s)}\|u(s)-\bar{u}_0\|_{L^p}ds	\\
		&\leq C_7e^{-(\lambda_1+1)t}+C_8\int_{0}^{t}\left(1+(t-s)^{-\frac{1}{2}-\frac{n}{2p}}\right)e^{-(\lambda_1+1)(t-s)}e^{-\frac{\delta}{2p}s}ds		\\
		&\leq C_7e^{-(\lambda_1+1)t}+C_8e^{-\frac{\delta}{2p}t}\int_{0}^{\infty}\left(1+z^{-\frac{1}{2}-\frac{n}{2p}}\right)e^{-(\lambda_1+1-\frac{\delta}{2p})z}dz	\\
		&\leq C_9e^{-\frac{\delta}{2p}t},
	\end{split}
	\ee
	where $p>n$ ensures $-\frac{1}{2}-\frac{n}{2p}>-1$ and the finiteness of the gamma function.
	
	Then we easily deduce from \eqref{v-average} and \eqref{v-vbar-estimate} that
	\be\label{v-vbar-decay}
	\begin{split}
		\|v(t)-\bar{u}_0\|_{W^{1,\infty}}&\leq\|v(t)-\bar{v}(t)\|_{W^{1,\infty}}+\|\bar{v}(t)-\bar{u}_0\|_{W^{1,\infty}}	\\
		&\leq C_{10}e^{-\frac{\delta}{2p}t}+|\bar{v}_0-\bar{u}_0|e^{-t}	\\
		&\leq C_{11}e^{-\min\{1,\frac{\delta}{2p}\}t},\quad\forall t>0,
	\end{split}
	\nonumber
	\ee
	which together with \eqref{v-wlinfty-decay} yields  our desired estimate \eqref{v-decay}.

	Finally, we proceed to show the exponential convergence of $u$. Again, the variation-of-constants formula applied to  the $u$-equation in \eqref{main equation} shows that
	\be\label{u-ubar-estimate}
	\begin{split}
		\|u(t)-\bar{u}_0\|_{L^\infty}&\leq\|e^{t\Delta}(u_0-\bar{u}_0)\|_{L^\infty}+\chi\int_{0}^{t}\|e^{(t-s)\Delta}\nabla\cdot(u\nabla(v-\bar{v}))(s)\|_{L^\infty}	\\
		&\ \ \ +\xi\int_{0}^{t}\|e^{(t-s)\Delta}\nabla\cdot(S(u)\nabla w)(s)\|_{L^\infty}ds	\\
		&=:J_1+J_2+J_3.
	\end{split}
	\ee
	Since $\int_\Omega(u_0-\bar{u}_0)=0$, the smoothing $L^p$-$L^q$-estimate of $(e^{t\Delta})_{t\geq0}$ gives
	\be\label{j1}
		J_1\leq C_{12}e^{-\lambda_1t},\quad\forall t>0.
	\ee
	For the term $J_2$, applying the $L^\infty$-boundedness of $u$ and the $W^{1,\infty}$-exponential decay of $v-\bar{v}$ in \eqref{v-wlinfty-decay} for $\tau=0$  or  \eqref{v-vbar-estimate} for $\tau=1$, we obtain  that
	\be\label{j2}
	\begin{split}
		 J_2&\leq C_{13}\int_{0}^{t}\left(1+(t-s)^{-\frac{1}{2}}\right)e^{-\lambda_1(t-s)}\|(u\nabla(v-\bar{v}))(s)\|_{L^\infty}ds	\\
		 &\leq C_{13}\int_{0}^{t}\left(1+(t-s)^{-\frac{1}{2}}\right)e^{-\lambda_1(t-s)}\|u(s)\|_{L^\infty}\|\nabla(v-\bar{v})(s)\|_{L^\infty}ds	\\
		 &\leq C_{14}\int_{0}^{t}\left(1+(t-s)^{-\frac{1}{2}}\right)e^{-\lambda_1(t-s)}e^{-\frac{\delta}{2p}s}ds	\\
		 &\leq C_{15}e^{-\frac{\delta}{2p}t},\quad\forall t>0.
	\end{split}
	\ee
Similarly, by the $W^{1,\infty}$-exponential decay  of $w$ in \eqref{w-decay} with $\gamma<\lambda_1$, it follows  that
	\be\label{j3}
	\begin{split}
		J_3&\leq C_{16}\int_{0}^{t}\left(1+(t-s)^{-\frac{1}{2}}\right)e^{-\lambda_1(t-s)}\|(S(u)\nabla w)(s)\|_{L^\infty}ds	\\
		&\leq C_{16}\int_{0}^{t}\left(1+(t-s)^{-\frac{1}{2}}\right)e^{-\lambda_1(t-s)}\|S(u)(s)\|_{L^\infty}\|\nabla w(s)\|_{L^\infty}ds	\\
		&\leq C_{16}\|u+1\|_{L^\infty}^{\theta^+}\int_{0}^{t}\left(1+(t-s)^{-\frac{1}{2}}\right)e^{-\lambda_1(t-s)}e^{-\gamma s}ds	\\
		&\leq C_{17}e^{-\gamma t},\quad\forall t>0.
	\end{split}
	\ee
	It follows from  $p>n$  and \eqref{delta-def}  that $\frac{\delta}{2p}<\min\{\lambda_1, \ \gamma\}$. Then a simple substitution of \eqref{j1}, \eqref{j2} and \eqref{j3} into \eqref{u-ubar-estimate} immediately show the exponential decay of $u$ to $\bar{u}_0$ in the sense \eqref{u-decay}. Thus we complete the proof.  	
\end{proof}

\section{Small-data global dynamics for $n\geq2$ and $\chi\in\mathbb{R}$}

In the previous section, we showed that any globally bounded solution of the IBVP \eqref{main equation}  converges exponentially to its constant equilibrium. In this section, as a preparatory and    inspired by \cite{Winkler10JDE, Cao15}, we examine the global-in-time existence, boundedness and convergence of classical solutions to \eqref{main equation} for $\theta\leq1$ and $\chi\in\mathbb{R}$ under certain  smallness condition on initial data, and thus proving Theorem \ref{theorem3}. These boundedness and convergence  also improve and extend \cite[Theorem 1.7]{RT24JDE}.

\begin{lemma}
	Let $n\geq 2$ and $\chi\in\mathbb{R}$. If there exists  $\eta_0>0$ such that \eqref{small-initial-data} holds, then for $p\in(n, \infty)$ and $q\in (p,\infty)$, there exist $C>0$ and $\mu\in (0, \min\{\lambda_1, \bar{u}_0\})$ such that the following decay estimates hold for $t>0$:
	\be\label{u-grad v-decay}
	\begin{cases}
		\|u(t)-\bar{u}_0\|_{L^p}\leq C\eta_0\left(1+t^{-1+\frac{n}{2p}}\right)e^{-\mu t}, \\[0.25cm]
		\|\nabla v(t)\|_{L^q}+\|\nabla w(t)\|_{L^q}\leq C\eta_0(1+t^{-\frac{1}{2}+\frac{n}{2q}})e^{-\mu t}.
	\end{cases}
	\ee
	Moreover, the  solution of \eqref{main equation} exists globally and is uniformly bounded according to \eqref{global-boundedness}, and decays exponentially to   $(\bar{u}_0,\bar{u}_0,0)$ in the sense of \eqref{3D-convergence1}.
\end{lemma}
\begin{proof}
	For $\eta_0>0$ to be determined below as in \eqref{epsilon0-def} and assuming \eqref{small-initial-data}, then there exists $\eta\in(0,\eta_0)$ such that
	\be\label{small-data2}
	\|u_0\|_{L^\frac{n}{2}}\leq\eta,\ \tau\|\nabla v_0\|_{L^n}\leq\eta\text{ and }\| \nabla w_0\|_{L^\infty}\leq\eta.
	\ee
	
	\textbf{Step 1:} (Global existence and boundedness) For  any  $p\in(n,\infty)$ and $0<\mu<\min\{\lambda_1,\bar{u}_0\}$, we first define
	\be\label{criterion}
	\begin{split}
		T:=\sup\Bigr\{\tilde{T}\in (0, T_m):  \  & \|u(t)-e^{t\Delta}u_0\|_{L^p} \\
		   & \leq\eta(1+t^{-1+\frac{n}{2p}})e^{-\mu t}\text{ for all } t\in(0,\tilde{T})\Bigr\}.
	\end{split}
	\ee
	 Noticing that $T>0$ is well-defined by the continuity of solutions ensured by the local well-posedness in Lemma \ref{local-existence-lemma}, we shall show that $T=\infty$ via a loop argument. To this end, by  \eqref{small-data2}, the definition of $T$ in \eqref{criterion} and the Neumann heat semigroup estimates, for $ p\in(n,\infty)$, we first derive that
	\be\label{u-u0-lp}
	\begin{split}
		\|u(t)-\bar{u}_0\|_{L^p}&\leq\|u(t)-e^{t\Delta}u_0\|_{L^p}+\|e^{t\Delta}\left(u_0-\bar{u}_0\right)\|_{L^p}	\\
		&\leq\eta(1+t^{-1+\frac{n}{2p}})e^{-\mu t}+C_1(1+t^{-1+\frac{n}{2p}})e^{-\lambda_1 t}\|u_0-\bar{u}_0\|_{L^\frac{n}{2}}	\\
		&\leq C_2\eta(1+t^{-1+\frac{n}{2p}})e^{-\mu t}.
	\end{split}
	\ee
where and below,  we  always assume $t\in(0, T)$,  unless otherwise stated.

	In the case of $\tau=1$, as shown in the preceding section, we have
	\be\label{v-formula}
	v(t)-\bar{v}(t)=e^{t(\Delta-1)}(v_0-\bar{v}_0)+\int_{0}^{t}e^{(t-s)(\Delta-1)}(u(s)-\bar{u}_0)ds.
	\ee
	Thus, applying $\nabla$ to both sides of \eqref{v-formula} and using the Neumann heat semigroup estimate, the smallness condition \eqref{small-data2}, \eqref{u-u0-lp} and the estimate for integral in Lemma \ref{integral-lemma}, we derive, for $p=\max\{p,n\}\leq q<\frac{np}{(n-p)^+}=\infty$, that
	\be\label{v-lq}
	\begin{split}
		\|\nabla v(t)\|_{L^q}&\leq\|\nabla e^{t(\Delta-1)}(v_0-\bar{v}_0)\|_{L^q}+\int_{0}^{t}\|\nabla e^{(t-s)(\Delta-1)}(u(s)-\bar{u}_0)\|_{L^q}ds\\
		&\leq C_3(1+t^{-\frac{1}{2}+\frac{n}{2q}})e^{-(\lambda_1+1) t}\|\nabla v_0\|_{L^n}\\
		&\   +C_3\int_{0}^{t}\left(1+(t-s)^{-\frac{1}{2}-\frac{n}{2}(\frac{1}{p}-\frac{1}{q})}\right)e^{-(\lambda_1+1)(t-s)}\|u(s)-\bar{u}_0\|_{L^p}ds\\
		&\leq C_4\eta(1+t^{-\frac{1}{2}+\frac{n}{2q}})e^{-\mu t}.
	\end{split}		
	\ee
When $\tau=0$, one can easily get \eqref{v-lq} by the elliptic estimate as shown in \eqref{v-wlinfty-decay}.

	Again, with no loss of generality by shifting initial times, we proceed to assume that $v\geq v_*$ on $\Omega \times (0, T_m)$ as in \eqref{v-lower-bound-estimate}. Then in light of \eqref{w-lp}, \eqref{v-lq} and \eqref{small-data2},  we obtain, for $p\leq q<\infty$, that
	\be\label{grad w-lp}
	\begin{split}
		\|\nabla w(t)\|_{L^q}&\leq\|\nabla w_0\|_{L^q} e^{-v_* t}+\|w_0\|_{L^\infty}e^{-v_*t}\int_{0}^{t}\|\nabla v(s)\|_{L^q}ds\\
		&\leq  C_5 \eta e^{-v_* t}+C_5\eta e^{-v_*t}\int_{0}^{t}(1+s^{-\frac{1}{2}+\frac{n}{2q}})e^{-\mu s}ds\\
		&\leq C_6\eta e^{-v_* t}.
	\end{split}		
	\ee
	Now, applying the well-known  semigroup estimates to the semigroup representation of $u$-equation in \eqref{main equation} and noticing $0\leq S(u)\leq u$ due to $\theta\leq1$, we infer that
	\be\label{u-loop}
	\begin{split}
		\|u(t)-e^{t\Delta}u_0\|_{L^p}&\leq|\chi|\int_{0}^{t}\|e^{(t-s)\Delta}\nabla\cdot(u\nabla v)(s)\|_{L^p}ds\\
		&\ \ \ +|\xi|\int_{0}^{t}\|e^{(t-s)\Delta}\nabla\cdot(S(u)\nabla w)(s)\|_{L^p}ds\\
		&\leq C_7\int_{0}^{t}\left(1+(t-s)^{-\frac{1}{2}}\right)e^{-\lambda_1(t-s)}\|(u\nabla v)(s)\|_{L^p}ds\\
		&\ \ \ +C_7\int_{0}^{t}\left(1+(t-s)^{-\frac{1}{2}}\right)e^{-\lambda_1(t-s)}\|(u\nabla w)(s)\|_{L^p}ds\\
		&=:C_7  I_1+C_7I_2.
	\end{split}
	\ee
	To estimate $I_1$ and $I_2$, for any $p\in(n, \infty)$, we choose $ p_1\in (p,\  2p]$ and
 $$
 q_1=\frac{p_1p}{p_1-p}\Longleftrightarrow \frac{1}{p_1}+\frac{1}{q_1}=\frac{1}{p}.
 $$
 Then we find  that  $p_1\leq q_1<\infty$.  Thus, by  H\"{o}lder's interpolation inequality, \eqref{small-data2}, \eqref{u-u0-lp}, \eqref{v-lq} and the integral estimate in Lemma \ref{integral-lemma}, we conclude  that
	\be\label{i1}
	\begin{split}
		I_1&\leq\int_{0}^{t}\left(1+(t-s)^{-\frac{1}{2}}\right)e^{-\lambda_1(t-s)}\|((u-\bar{u}_0)\nabla v)(s)\|_{L^p}ds\\
		&\ \ \ +\int_{0}^{t}\left(1+(t-s)^{-\frac{1}{2}}\right)e^{-\lambda_1(t-s)}\|\bar{u}_0\nabla v(s)\|_{L^p}ds	\\
		&\leq \int_{0}^{t}\left(1+(t-s)^{-\frac{1}{2}}\right)e^{-\lambda_1(t-s)}\|u(s)-\bar{u}_0\|_{L^{p_1}}\|\nabla v(s)\|_{L^{q_1}}ds\\
		&\ \ \ +\int_{0}^{t}\left(1+(t-s)^{-\frac{1}{2}}\right)e^{-\lambda_1(t-s)}\|\bar{u}_0\|_{L^{p_1}}\|\nabla v(s)\|_{L^{q_1}}ds	\\
		&\leq C_8\eta^2\int_{0}^{t}\left(1+(t-s)^{-\frac{1}{2}}\right)e^{-\lambda_1(t-s)}(1+s^{-1+\frac{n}{2p_1}})(1+s^{-\frac{1}{2}+\frac{n}{2q_1}})e^{-2\mu s}ds\\
		&\ \ \ +C_9\eta^2\int_{0}^{t}\left(1+(t-s)^{-\frac{1}{2}}\right)e^{-\lambda_1(t-s)}(1+s^{-\frac{1}{2}+\frac{n}{2q_1}})e^{-\mu s}ds	\\
		&\leq C_{10}\eta^2(1+t^{-1+\frac{n}{2p}})e^{-\mu t}.
	\end{split}
	\ee
	In a similar fashion,  we use \eqref{small-data2}, \eqref{u-u0-lp}, \eqref{grad w-lp} to estimate $I_2$ as follows:
	\be\label{i2}
	I_2\leq C_{11}\eta^2(1+t^{-1+\frac{n}{2p}})e^{-\mu t}.
	\ee
	An easy substitution of  \eqref{i1} and \eqref{i2} into \eqref{u-loop} entails that
	$$
	\|u(t)-e^{t\Delta}u_0\|_{L^p}\leq C_{12}\eta^2(1+t^{-1+\frac{n}{2p}})e^{-\mu t}, \ \ \forall t\in (0, T).
	$$
Then, noticing $C_{12}$ is independent of $\eta$, upon choosing
 \be\label{epsilon0-def}
 \eta_0=\frac{1}{C_{12}}
 \ee
  and $\eta_1\in(\eta, \eta_0)$, we see, for all $t\in (0, T)$ and $p\in (n,  \infty)$,  that
$$		
			\|u(t)-e^{t\Delta}u_0\|_{L^p}\leq C_{12}\eta_1\eta\left(1+t^{-\frac{1}{2}+\frac{n}{2p}}\right)
e^{-\mu t},
$$
which coupled  with the fact $ C_{12}\eta_1<1$ and the definition of $T$ in \eqref{criterion} implies that $T=\infty$. Consequently, based on \eqref{criterion} and repeating the argument leading to \eqref{u-u0-lp}, and then using \eqref{v-lq} and  \eqref{grad w-lp}, we quickly obtain the desired decay estimate  \eqref{u-grad v-decay}.

Thus, trivially $T_m=\infty$, then, with the aid of \eqref{u-grad v-decay} and  the local smoothness of solutions  in Lemma \ref{local-existence-lemma}, one can easily see that  $u\in L^\infty((0,\infty); L^p(\Omega))$ for $p>n$. The latter enables one to adapt the argument in  \textbf{Step 3} in Lemma \ref{global-boundedness-lemma} to conclude the  desired global boundedness as in \eqref{global-boundedness}.
	
	\textbf{Step 2 :} (Exponential convergence)  Given the global boundedness, from Lemma \ref{w-convergence-lemma}, we first obtain the $W^{1,\infty}$-exponential decay of $w$. Then, on the basis of $L^p$-exponential decay in \eqref{u-grad v-decay}, modifying the arguments from \eqref{u-lp-decay} with $\frac{\delta}{2p}$ replaced by $\mu$  to the end proof of Lemma \ref{uv-decay-lemma}, one can readily derive the exponential decay of  $(u(t)-\bar{u}_0, v(t)-\bar{u}_0,w(t))$ in the sense of \eqref{3D-convergence1}, as desired.
\end{proof}

\section{Global existence, eventual regularity and exponential convergence of  weak solutions in 3D}
	
	In this section, we aim to prove Theorem \ref{theorem4} by deriving the global existence, eventual regularity and exponential convergence of  weak solutions to the IBVP \eqref{main equation} with $\tau=1$ in 3D. To this end, we first state a natural notion of weak solutions based on integration by parts and introduce an approximation problem of \eqref{main equation}, and then derive some $\epsilon$-independent estimates of approximating solutions, which are crucial to derive global existence of weak solutions. Finally, to establish the eventual regularity of global weak solutions, we adopt the method based on an analysis of the corresponding linearized system and the stability of constant solutions in scaling-invariant spaces,  developed by Jiang in \cite{Jiang19, Jiang22} for chemo-repulsion systems. To avoid more repetitions, unless otherwise stated, we shall tacitly assume that $\tau=1$ and $n=3$ in this section.
	
\subsection{Existence of global weak solutions to the IBVP \eqref{main equation}}

	Based on integration by parts, we first introduce the following definition of weak solutions.

\begin{definition}\label{weak-definition}
	A global weak solution of \eqref{main equation} is a pair of nonnegative functions
	$$
		(u,v)\in C([0,\infty);weak-L^1(\Omega;\mathbb{R}^2)), \ \  \ w\in C([0,\infty);weak-L^\infty(\Omega))
	$$
  such that
 	$$
 		\nabla u, \nabla v, u\nabla v \text{ and } S(u)\nabla w    \text{ belong to } L_{loc}^1([0,\infty); \left(L^1(\Omega)\right)^3)
 	$$
 	and that the following identities from integration by parts
 	\be\label{weak-cond1}
 	\begin{split}
 		\int_\Omega(u(t)-u_0)\varphi=&-\int_{0}^{t}\int_\Omega\nabla u\cdot\nabla\varphi-\chi\int_{0}^{t}\int_\Omega u\nabla v\cdot\nabla\varphi+\xi\int_{0}^{t}\int_\Omega S(u)\nabla w\cdot\nabla\varphi,
 	\end{split}		
 	\nonumber
 	\ee
 	\be\label{weak-cond2}
 		 \int_\Omega (v(t)-v_0)\varphi=-\int_{0}^{t}\int_\Omega\nabla v\cdot\nabla\varphi-\int_{0}^{t}\int_\Omega v\varphi+\int_{0}^{t}\int_\Omega u\varphi
 	\nonumber
 	\ee
 	and
 	\be\label{weak-cond3}
 		\int_\Omega (w(t)-w_0)\varphi=-\int_{0}^{t}\int_\Omega vw\varphi
 	\nonumber
 	\ee
 	hold for each $t\geq0$ and $\varphi\in W^{1,\infty}(\Omega)$.
\end{definition}

	To derive global weak solution  of  \eqref{main equation} in the sense of Definition \ref{weak-definition}, for any $\epsilon >0$, we investigate  the following  approximating system corresponding to \eqref{main equation}:
	\be\label{regularized-equation}
	\begin{cases}
		u_{\epsilon t} =  \Delta u_\epsilon+\chi \nabla\cdot(u_\epsilon (1-\epsilon u_\epsilon)\nabla v_\epsilon) \\	\ \ \ \ \ \ \ \
-\xi\nabla\cdot(S(u_\epsilon) (1-\epsilon u_\epsilon)\nabla w_\epsilon)  &\text{in } \Omega\times(0,\infty), \\[0.2cm]
		 v_{\epsilon t} = \Delta v_\epsilon-v_\epsilon+u_\epsilon  &\text{in } \Omega\times(0,\infty),\\[0.2cm]
		w_{\epsilon t}=-v_\epsilon w_\epsilon 	 &\text{in } \Omega\times(0,\infty),\\[0.2cm]
		\frac{\partial u_\epsilon}{\partial\nu}=\frac{\partial v_\epsilon}{\partial\nu}=\frac{\partial w_\epsilon}{\partial\nu}=0  &\text{on } \partial\Omega\times(0,\infty), \\[0.2cm]
		(u_\epsilon,v_\epsilon,w_\epsilon)=(u_0,v_0,w_0)  &\text{on } \Omega\times\{0\}.
	\end{cases}
	\ee
Here, we adopt  the equivalent homogeneous boundary conditions in \eqref{nbc}.  A formal limit of $\epsilon\searrow 0$ in \eqref{regularized-equation} yields \eqref{main equation}. In the sequel, we shall demonstrate that such limit process is possible in the sense of Definition \ref{weak-definition}. To this purpose, we first show the  global  well-posedness for the approximating system \eqref{regularized-equation}.
\begin{lemma}\label{regular-system-local-lemma}
	 Let $\chi,\xi >0$, $S$ satisfy \eqref{S-con},  $\Omega\subset \mathbb{R}^n(n\geq 1)$ be a bounded and smooth domain, and let the initial data $(u_0,v_0,w_0)$ fulfill  \eqref{initial-data} and $0<\epsilon < \|u_0\|_{L^\infty}^{-1}$. Then the IBVP  \eqref{regularized-equation} possesses a unique global  classical solution pair
	$$
		(u_\epsilon,v_\epsilon,w_\epsilon)\in[C^0(\bar{\Omega}\times[0,\infty)\cap C^{2,1}(\bar{\Omega}\times(0,\infty))]^3,
	$$
which satisfy, for  $(x,t)\in \bar\Omega\times (0, \infty)$,
	\be\label{u-l1-estimate2}
		0<u_\epsilon (x,t)<\frac{1}{\epsilon},\ \ 0<v_\epsilon(x,t), \ \   w_\epsilon(x,t)=w_0(x)e^{-\int_{0}^{t}v_\epsilon(x,s)ds},
	\ee
	and are uniformly bounded in $(L^1, L^1, L^\infty)$:
	\be\label{v-l1-estimate2}
		\|u_\epsilon(t)\|_{L^1}=\|u_0\|_{L^1},  \   \|v_\epsilon(t)\|_{L^1}\leq \max\left\{ \|v_0\|_{L^1},    \|u_0\|_{L^1}\right\}, \   0\leq w_\epsilon\leq  \|w_0\|_{L^\infty}.
	\ee
\end{lemma}
\begin{proof}
	Using the well-established arguments in chemotaxis systems as illustrated   in Lemma \ref{local-existence-lemma}, within its lifespan  $(0,T_{m\epsilon})$, one can easily derive the local well-posedness,   estimates   \eqref{v-l1-estimate2} and \eqref{u-l1-estimate2} except $u_\epsilon<\frac{1}{\epsilon}$,  and if $T_{m\epsilon}<\infty$, then
\be\label{extensibility-criterion+}
		\limsup_{t\nearrow T_m}\left(\|(u_\epsilon(t)\|_{L^{\infty}}+\|v_\epsilon(t)\|_{W^{1,\infty}}+\|w_\epsilon(t)\|_{W^{1,\infty}}\right)=\infty.
	\ee
We next invoke arguments based on maximum principle  for parabolic equations as done in \cite[Theorem 2.5]{Xiang15} to  show that $T_{m\epsilon}=\infty$ and $u_\epsilon<\frac{1}{\epsilon}$ on $\bar\Omega \times [0,\infty)$. Indeed,  noticing $0<\epsilon < \|u_0\|_{L^\infty}^{-1}$,  we expand  out the $u_\epsilon$-equation in \eqref{regularized-equation} to find that
	\be\label{u-equation}
	\begin{cases}
		u_{\epsilon t}-\Delta u_\epsilon+b(x,t)\nabla u_\epsilon+c(x,t)(1-\epsilon u_\epsilon)u_\epsilon=0, & x\in\Omega,\ t>0,	\\[0.2cm]
		\frac{\partial u_\epsilon}{\partial\nu}=0,  &x\in\partial\Omega,\ t>0,	\\[0.2cm]
		u_\epsilon(x,0)=u_0(x)<\frac{1}{\epsilon}, & x\in\Omega,
	\end{cases}
	\nonumber
	\ee
	where
	$$
		b(x,t)=-\chi(1-2\epsilon u_\epsilon)\nabla v_\epsilon+\xi [ S^\prime(u_\epsilon)(1-\epsilon u_\epsilon) -\epsilon S(u_\epsilon)]\nabla w_\epsilon
	$$
	and
	$$
		c(x,t)=-\chi\Delta v_\epsilon+\xi\frac{S(u_\epsilon)}{u_\epsilon}\Delta w_\epsilon.
	$$
Thus, applying the comparison principle, the strong maximum principle and the Hopf boundary point lemma to the above system, we immediately obtain that $u_\epsilon(x,t)<\frac{1}{\epsilon}$ for  $(x,t)\in \bar{\Omega}\times[0,T_{m\epsilon})$. Then we deduce  from Lemma \ref{v-lq-norm-lemma} that $v_\epsilon\in L^\infty(0, T_{m\epsilon}; W^{1,\infty}(\Omega))$, which along with applications of Lemmas \ref{v-lower-bound-lemma} and \ref{w-lemma} to the second and third equation in \eqref{regularized-equation} further  shows that $w_\epsilon \in L^\infty(0, T_{m\epsilon}; W^{1,\infty}(\Omega))$. Thus, the extensibility criterion \eqref{extensibility-criterion+} forces $T_{m\epsilon}=\infty$.
\end{proof}
	
To show the regularity of approximating solutions, we establish an energy-dissipation identity for the approximating model \eqref{regularized-equation}, which is crucial in deriving the global existence of weak solutions.
	
\begin{lemma}\label{funtional-lemma}
	Under Lemma \ref{regular-system-local-lemma}, the global classical solution $(u_\epsilon,v_\epsilon,w_\epsilon)$ of the IBVP \eqref{regularized-equation} fulfills, for $t>0$, that
	\be\label{h-functional}
		\mathcal{H_\epsilon}^\prime(t)+\int_\Omega\frac{|\nabla u_\epsilon|^2}{u_\epsilon(1-\epsilon u_\epsilon)}+\chi\int_\Omega(|\Delta v_\epsilon|^2+|\nabla v_\epsilon|^2)=\xi\int_\Omega\frac{S(u_\epsilon)}{u_\epsilon}\nabla u_\epsilon\cdot\nabla w_\epsilon,
	\ee
	where $\mathcal{H_\epsilon}$ is given by
	\be\label{H-def}
		\mathcal{H_\epsilon}(t)=\int_\Omega\left(u_\epsilon\ln u_\epsilon+\frac{1}{\epsilon}(1-\epsilon u_\epsilon)\ln(1-\epsilon u_\epsilon)+\frac{\chi}{2}|\nabla v_\epsilon|^2\right).
	\ee
\end{lemma}
\begin{proof}
Noticing $1-\epsilon u_\epsilon>0$ on $\bar{\Omega}\times[0,\infty)$, then  multiplying the $u_\epsilon$-equation in \eqref{regularized-equation} by $\ln u_\epsilon-\ln(1-\epsilon u_\epsilon)$ and   integrating by parts over $\Omega$, we get  that
	\be\label{h1}
	\begin{split}
		&\frac{d}{dt}\int_\Omega\left(u_\epsilon\ln u_\epsilon+\frac{1}{\epsilon}(1-\epsilon u_\epsilon)\ln(1-\epsilon u_\epsilon)\right)	\\
		&=-\int_\Omega\frac{|\nabla u_\epsilon|^2}{u_\epsilon(1-\epsilon u_\epsilon)}-\chi\int_\Omega\nabla u_\epsilon\cdot\nabla v_\epsilon+\xi\int_\Omega\frac{S(u_\epsilon)}{u_\epsilon}\nabla u_\epsilon\cdot\nabla w_\epsilon.
	\end{split}		
	\ee
	On the other hand, taking $\nabla$ on  the $v_\epsilon$-equation in \eqref{regularized-equation} and  then multiplying by $\nabla v_\epsilon$ on both sides and integrating by parts over $\Omega$, we derive that
	$$
		\frac{1}{2}\frac{d}{dt}\int_\Omega|\nabla v_\epsilon|^2=-\int_\Omega|\Delta v_\epsilon|^2-\int_\Omega|\nabla v_\epsilon|^2+\int_\Omega\nabla u_\epsilon\cdot\nabla v_\epsilon,
	$$
	which along with \eqref{h1} quickly entails \eqref{h-functional} with $\mathcal{H_\epsilon}(t)$ given by \eqref{H-def}.
\end{proof}

	Based on the evolution identity established by Lemma \ref{funtional-lemma}, we obtain the uniform-in-time boundedness of $\|u\ln u\|_{L^1}$ under a restriction on the growth of $S$.

\begin{lemma}\label{ulnu-epsilon-lemma}
	Under  Lemma \ref{regular-system-local-lemma} and  $\theta<\frac{4}{3}$,   the global classical solution $(u_\epsilon,v_\epsilon,w_\epsilon)$ of the IBVP \eqref{regularized-equation} fulfills, for $t>0$, that
	\be\label{ulnu-epsilon}
	\int_\Omega\left(|u_\epsilon\ln u_\epsilon|+|\nabla v_\epsilon|^2\right)\leq C
	\ee
	and
	\be\label{delta v-grad v-estimate-2}
	\int_{t}^{t+1}\int_\Omega\left(\frac{|\nabla u_\epsilon|^2}{u_\epsilon}+|\Delta v_\epsilon|^2+|\nabla v_\epsilon|^2\right)\leq C,
	\ee
where the constant  $C>0$ is independent of $\epsilon$.
\end{lemma}
\begin{proof}
	Arguing similarly as done in Lemma \ref{grad u-grad w-lemma},  we proceed with $0\leq S(u_\epsilon)\leq u_\epsilon(1+u_\epsilon)^{\theta-1}$ with $0<\theta<\frac{4}{3}$.  Then, by the general one-sided upper bound of $-\Delta w_\epsilon$ in \eqref{delta-w-estimate0}, we estimate that
	\be\label{S-u-w}
	\begin{split}
		&\xi\int_\Omega\frac{S(u_\epsilon)}{u_\epsilon}\nabla u_\epsilon\cdot\nabla w_\epsilon	\\
		&=-\xi\int_\Omega\Phi_1(u_\epsilon)\Delta w_\epsilon	\\
		&\leq\frac{\xi\|w_0\|_{L^\infty}}{\theta}\int_\Omega(u_\epsilon+1)^\theta v_\epsilon-\frac{2\xi}{\theta}\int_\Omega(u_\epsilon+1)^\theta\nabla w_0\cdot\nabla e^{-\int_{0}^{t}\nabla v_\epsilon(s)ds}	\\
		&\ \ \ +\frac{\xi}{\theta}\left(\frac{\|w_0\|_{L^\infty}}{e}+\|\Delta w_0\|_{L^\infty}\right)\int_\Omega(u_\epsilon+1)^\theta,	
	\end{split}
	\ee
	where $\Phi_1$ is defined by \eqref{phi1} and satisfies $0\leq \Phi_1(u_\epsilon)\leq (1+u_\epsilon)^\theta/\theta$.
	
	Thanks to $\theta<\frac{4}{3}$,  we first fix $q$ to satisfy
\be\label{q-theta-esp}
		\max\left\{1,\  \frac{3}{3-\theta}, \   \frac{3}{5-3\theta}\right\}<q<\min\left\{3, \ \frac{1}{(1-\theta)^+}\right\}.
	\ee
	Next, by H\"{o}lder's inequality, the 3D Gagliardo-Nirenberg inequality, the $L^q$-bound of $v_\epsilon$ with $q\in[1,3)$ in \eqref{v-lq-norm-estimate1} and the mass identity $\|u_\epsilon\|_{L^1}=\|u_0\|_{L^1}$, we deduce,  for any $\eta_1>0$, that
	\be\label{term-utheta-v}
	\begin{split}
		\int_\Omega (u_\epsilon+1)^\theta v_\epsilon&\leq\|v_\epsilon\|_{L^q}\left(\int_\Omega(u_\epsilon+1)^{\frac{q\theta }{q-1}}\right)^{\frac{q-1}{q}}	\\
		&\leq C_1\|(u_\epsilon+1)^\frac{1}{2}\|_{L^\frac{2q\theta}{q-1}}^{2\theta}\\
		&\leq C_2\left(\|\nabla (u_\epsilon+1)^\frac{1}{2}\|_{L^2}^{3\theta-\frac{3(q-1)}{q}}
\|(u_\epsilon+1)^\frac{1}{2}\|_{L^2}^{\frac{3(q-1)}{q}-\theta}+1\right)	\\
		&\leq\eta_1\int_\Omega\frac{|\nabla u_\epsilon|^2}{u_\epsilon}+C_{\eta_1},	
	\end{split}
	\ee
	where, due to  the    choice of $q$ in \eqref{q-theta-esp}, the fact that $3\theta-\frac{3(q-1)}{q}<2$ was used.
	
	Similarly, using the 3D GN inequality or simply applying the estimate  \eqref{u0.5p} with $p=1,n=3$,  for  $\theta<\frac{5}{3}$ and for  any $\eta_2>0$, one can see  that
	\be\label{term-utheta}
	\int_\Omega (u_\epsilon+1)^\theta\leq\eta_2\int_\Omega\frac{|\nabla u_\epsilon|^2}{u_\epsilon}+C_{\eta_2}.
	\ee
	As for the second term in \eqref{S-u-w}, we use integration by parts and Young's inequality with epsilon  to derive, for any $\eta_3>0$, that
	\be\label{term22}
	\begin{split}
		&-\int_\Omega(u_\epsilon+1)^\theta\nabla w_0\cdot\nabla e^{-\int_{0}^{t}v_\epsilon(s)ds}	\\
		&=\theta\int_\Omega(u_\epsilon+1)^{\theta-1}e^{-\int_{0}^{t} v_\epsilon(s)ds}\nabla u_\epsilon\cdot\nabla w_0+\int_\Omega(u_\epsilon+1)^\theta e^{-\int_{0}^{t} v_\epsilon(s)ds}\Delta w_0	\\
		&\leq\theta\|\nabla w_0\|_{L^\infty}\int_\Omega(u_\epsilon+1)^{\theta-1}|\nabla u_\epsilon|+\|\Delta w_0\|_{L^\infty}\int_\Omega(u_\epsilon+1)^\theta	\\
		&\leq\theta\|\nabla w_0\|_{L^\infty}\left(\eta_3\int_\Omega\frac{|\nabla u_\epsilon|^2}{u_\epsilon}+\frac{1}{4\eta_3}\int_\Omega(u_\epsilon+1)^{2\theta-1}\right)	\\
		&\ \ \ +\|\Delta w_0\|_{L^\infty}\int_\Omega(u_\epsilon+1)^\theta.
	\end{split}
	\ee
	Since $2\theta-1<\frac{5}{3}$, applying \eqref{u0.5p} again shows, for any $\eta_4>0$, that
	\be\label{term23}
		\int_\Omega(u_\epsilon+1)^{2\theta-1}\leq\eta_4\int_\Omega\frac{|\nabla u_\epsilon|^2}{u_\epsilon}+C_{\eta_4}.
	\ee
	Substituting \eqref{term-utheta-v},  \eqref{term-utheta}, \eqref{term22} and \eqref{term23} into \eqref{S-u-w} and choosing  $\eta_i>0 (i=1,2,3,4)$ sufficiently small, we deduce,  for $\theta<\min\{\frac{4}{3},\frac{5}{3}\}=\frac{4}{3}$, that
	\be\label{Suw-estimate}
		\xi\int_\Omega\frac{S(u_\epsilon)}{u_\epsilon}\nabla u_\epsilon\cdot\nabla w_\epsilon\leq\frac{1}{4}\int_\Omega\frac{|\nabla u_\epsilon|^2}{u_\epsilon}+C_3.
	\ee
	For any $\eta_5>0$, since  $z\ln z\leq \eta_5 z^\frac{4}{3}+C_{\eta_5}$ for $z>0$, one can readily deduce from  Young's inequality and \eqref{u0.5p} that
	\be\label{ulnu-epsilon1}
		\int_\Omega u_\epsilon\ln u_\epsilon\leq\frac{1}{4}\int_\Omega\frac{|\nabla u_\epsilon|^2}{u_\epsilon}+C_4,
	\ee
Observing that
	$$
		\frac{|\nabla u_\epsilon|^2}{u_\epsilon(1-\epsilon u_\epsilon)}\geq\frac{|\nabla u_\epsilon|^2}{u_\epsilon},\ \ \  \frac{1}{\epsilon}(1-\epsilon  u_\epsilon)\ln(1-\epsilon  u_\epsilon)<0,
	$$
 plugging \eqref{Suw-estimate} into the evolution identity \eqref{h-functional}, combining the definition of $\mathcal{H_\epsilon}$ in \eqref{H-def} and \eqref{ulnu-epsilon1}, we deduce that
	\be\label{h2}
	\mathcal{H_\epsilon}^\prime(t)+\mathcal{H_\epsilon}(t)+\frac{1}{2}\int_\Omega\frac{|\nabla u_\epsilon|^2}{u_\epsilon}+\chi\int_\Omega|\Delta v_\epsilon|^2+\frac{\chi}{2}\int_\Omega|\nabla v_\epsilon|^2\leq C_5.
	\ee
Multiplying this  Gronwall inequality  by $e^t$, integrating on $(0,t)$ and noticing $\|u_\epsilon\|_{L^1}=\|u_0\|_{L^1}$, $0\leq \epsilon u_\epsilon<1$ and $0\leq \epsilon u_0<1$,  we infer, for $t>0$, that
\begin{equation*}
\begin{split}
&\int_\Omega \left(u_\epsilon\ln u_\epsilon+ u_\epsilon+\frac{1}{\epsilon}(1-\epsilon  u_\epsilon)\ln(1-\epsilon  u_\epsilon)+\frac{\chi}{2}|\nabla v_\epsilon|^2\right)(t)\\
&\leq \int_\Omega \left( u_0\ln u_0+ u_0+\frac{1}{\epsilon}(1-\epsilon  u_0)\ln(1-\epsilon  u_0)+\frac{\chi}{2} |\nabla v_0|^2\right)+C_5\\
&\leq \int_\Omega \left( u_0\ln u_0+ u_0+\frac{\chi}{2} |\nabla v_0|^2\right)+C_5,
\end{split}
\end{equation*}
which along with the facts
\be\label{ulnu-bdd-key}
u_\epsilon \ln  u_\epsilon\geq- e^{-1}, \ \  \   u_\epsilon+\frac{1}{\epsilon}(1-\epsilon  u_\epsilon)\ln(1-\epsilon  u_\epsilon)\geq0,
\ee
   enables us to achieve the desired  estimate \eqref{ulnu-epsilon}. Then integrating \eqref{h2} on $(t,t+1)$, we achieve the desired space-time estimate \eqref{delta v-grad v-estimate-2}.	
\end{proof}

	Before deriving  another uniform boundedness, it is a good place for us to  show a crucial $L^2$-exponential decay of  $\nabla w_\epsilon$, which is essential in proving the eventual regularity of weak solutions in Subsection 6.2 below.
	
\begin{lemma}\label{gradw-l2-decay}
	Under Lemma \ref{regular-system-local-lemma} and  $\inf_\Omega v_0>0$, there exists $v_*>0$ independent of $\epsilon$ such that  the solution component $v_\epsilon$ of \eqref{regularized-equation} has a positive lower bound:
 \be\label{v-low-bdd-eps}
 v_\epsilon\geq v_* \text{ on } \bar\Omega \times [0,\infty),
 \ee
 and  for  $0<\tilde{\gamma}<2v_*$,   the  solution component $w_\epsilon$  fulfills, for some $C>0$, that
	\be\label{grad w-l2}
		\int_\Omega|\nabla w_\epsilon|^2\leq Ce^{-\tilde{\gamma}t},\quad\forall t>0.
	\ee
\end{lemma}
\begin{proof}
Applying Lemma \ref{v-lower-bound-lemma} to the second equation in \eqref{regularized-equation}, we see there exists an increasing function $\zeta: [0,\infty)\mapsto \mathbb{R}^+$ such that
\be\label{v-low-bdd-pf-ep}
v_\epsilon(x,t)\geq \|u_0\|_{L^1} \zeta(t), \ \ \ \forall (x,t)\in \Omega\times [0,\infty),
\ee
On the other hand, the comparison principle  gives rise to
$$
v_\epsilon(x,t)\geq e^{-t} \inf_{y\in\Omega} v_0(y), \ \ \ \forall (x,t)\in \Omega\times [0,\infty),
$$
which along with \eqref{v-low-bdd-pf-ep} and the fact $\inf_\Omega v_0>0$  quickly implies
$$
v_\epsilon(x,t)\geq \sup_{s\in (0,\infty)}\min\left\{\|u_0\|_{L^1} \zeta(s), \ \  e^{-s} \inf_{y\in\Omega} v_0(y)\right\}>0, \ \  \forall (x,t)\in \Omega\times [0,\infty),
$$
yielding, upon an evident choice of $v_*$,  our desired uniform bound \eqref{v-low-bdd-eps}.

Next, by the $w_\epsilon$-equation in \eqref{regularized-equation}, we  compute directly that
	$$
		\nabla w_\epsilon(t)=\nabla w_0e^{-\int_{0}^{t}v_\epsilon(s)ds}-w_0e^{-\int_{0}^{t}v_\epsilon(s)ds}\int_{0}^{t}\nabla v_\epsilon(s)ds,
	$$
and so, we see
	$$
		|\nabla w_\epsilon(t)|^2\leq2|\nabla w_0|^2e^{-2\int_{0}^{t}v_\epsilon(s)ds}+2\|w_0\|_{L^\infty}^2e^{-2\int_{0}^{t}v_\epsilon(s)ds}\left(\int_{0}^{t}\nabla v_\epsilon(s)ds\right)^2.
$$
	Integrating the inequality above over $\Omega$, applying the uniform lower bound of $v_\epsilon$ in \eqref{v-low-bdd-eps} and using integration by parts, we deduce, for $0<\tilde{\gamma}<2v_*$ and $t>0$, that
	\be
	\begin{split}
		\int_\Omega|\nabla w_\epsilon(t)|^2&\leq2\|\nabla w_0\|_{L^2}^2e^{-2v_*t}+2\|w_0\|_{L^\infty}^2\int_\Omega e^{-2\int_{0}^{t}v_\epsilon(s)ds}\left(\int_{0}^{t}\nabla v_\epsilon(s)ds\right)^2	\\
		&=2\|\nabla w_0\|_{L^2}^2e^{-2v_*t}-\|w_0\|_{L^\infty}^2\int_\Omega\nabla\cdot\left(e^{-2\int_{0}^{t}v_\epsilon(s)ds}\int_{0}^{t}\nabla v_\epsilon(s)ds\right)	\\
		&\ \ \ +\|w_0\|_{L^\infty}^2\int_\Omega e^{-2\int_{0}^{t}v_\epsilon(s)ds}\int_{0}^{t}\Delta v_\epsilon(s)ds	\\
		&=2\|\nabla w_0\|_{L^2}^2e^{-2v_*t}+\|w_0\|_{L^\infty}^2\int_\Omega e^{-2\int_{0}^{t}v_\epsilon(s)ds}\int_{0}^{t}(v_{\epsilon t}+v_\epsilon-u_\epsilon)(s)ds	\\
		&=2\|\nabla w_0\|_{L^2}^2e^{-2v_*t}+\|w_0\|_{L^\infty}^2\int_\Omega e^{-2\int_{0}^{t}v_\epsilon(s)ds}(v_{\epsilon}(t)-v_0)	\\
		&\ \ \ +\|w_0\|_{L^\infty}^2\int_\Omega e^{-2\int_{0}^{t}v_\epsilon(s)ds}\int_{0}^{t}(v_\epsilon(s)-u_\epsilon(s))ds	\\
		&\leq2\|\nabla w_0\|_{L^2}^2e^{-2v_*t}+\|w_0\|_{L^\infty}^2e^{-2v_*t}\left(\int_\Omega v_\epsilon+\int_0^t\int_\Omega v_\epsilon\right)\\
&\leq 2\|\nabla w_0\|_{L^2}^2e^{-2v_*t}+\max\{\|v_0\|_{L^1},\|u_0\|_{L^1}\}(1+t)e^{-2v_*t}\\
&\leq \left(2\|\nabla w_0\|_{L^2}^2+\max\{\|v_0\|_{L^1},\|u_0\|_{L^1}\}(1+\frac{1}{2v_*e})\right)e^{-\tilde{\gamma}t},
	\end{split}
	\nonumber
	\ee
	where we used the facts  $\Delta v_\epsilon= v_{\epsilon t}+v_\epsilon-u_\epsilon$ and the nonnegativity of $v_0$ and $u_\epsilon$.
\end{proof}

In order to lift the estimates  \eqref{ulnu-epsilon} and  \eqref{delta v-grad v-estimate-2}, thanks to the logarithmic-type Sobolev  inequality in Lemma \ref{ext-ck} and  the crucial $L^2$-decay estimate \eqref{grad w-l2} for   $\nabla w_\epsilon$, we improve  the estimates  \eqref{ulnu-epsilon} and  \eqref{delta v-grad v-estimate-2} under a further restriction on $\theta$.

\begin{lemma}\label{exponential-decay-lemma}
	Under  Lemma \ref{gradw-l2-decay} and  $\theta\leq \frac{1}{2}$, there exists some $C>0$ such that
	\be\label{ulnu-epsilon2}
		\int_\Omega\left(u_\epsilon\ln\frac{u_\epsilon}{\bar{u}_0}+\frac{1}{\epsilon}(1-\epsilon u_\epsilon)\ln(1-\epsilon u_\epsilon)+\frac{\chi}{2}|\nabla v_\epsilon|^2\right)\leq Ce^{-\tilde\delta t}, \ \ \forall t>0,
	\ee
	and
	\be\label{delta v-grad v-estimate-1}
		\int_{0}^{\infty}\int_\Omega\left(\frac{|\nabla u_\epsilon|^2}{u_\epsilon}+|\Delta v_\epsilon|^2+|\nabla v_\epsilon|^2\right)\leq C,
	\ee
	where $\tilde{\delta}>0$ is defined by \eqref{delta2-def}.
\end{lemma}	
\begin{proof}	
	Proceeding similarly as Lemma \ref{funtional-lemma}, we use a refined functional $\tilde{\mathcal{H}}_\epsilon$  given by
	\be\label{h2-def} \mathcal{\tilde{H}}_\epsilon(t)=\int_\Omega\left(u_\epsilon\ln\frac{u_\epsilon}{\bar{u}_0}+\frac{1}{\epsilon}(1-\epsilon u_\epsilon)\ln(1-\epsilon u_\epsilon)+\frac{\chi}{2}|\nabla v_\epsilon|^2\right),
	\ee
which, due to $\int_\Omega u_\epsilon=\int_\Omega u_0$ and Lemma \ref{funtional-lemma},  fulfills
	\be\label{h2-functional}
	\mathcal{\tilde{H}}_\epsilon^\prime(t)+\int_\Omega\frac{|\nabla u_\epsilon|^2}{u_\epsilon(1-\epsilon u_\epsilon)}+\chi\int_\Omega(|\Delta v_\epsilon|^2+|\nabla v_\epsilon|^2)=\xi\int_\Omega\frac{S(u_\epsilon)}{u_\epsilon}\nabla u_\epsilon\cdot\nabla w_\epsilon.
	\ee
Recalling that $0\leq S(u)\leq u(u+1)^{\theta-1}$ with $\theta\leq\frac{1}{2}$,  the mass identity $\|u_\epsilon\|_{L^1}=\|u_0\|_{L^1}$, the $L^2$-exponential decay of $\nabla w_\epsilon$ in \eqref{grad w-l2}, we infer, for $0<\tilde\gamma<2v_*$, that
	 \be\label{temp1}
	 \begin{split}
	 	\xi\int_\Omega\frac{S(u_\epsilon)}{u_\epsilon}\nabla u_\epsilon\cdot\nabla w_\epsilon&\leq\frac{1}{2}\int_\Omega\frac{|\nabla u_\epsilon|^2}{u_\epsilon}+\frac{\xi^2}{2}\int_\Omega\frac{S^2(u_\epsilon)}{u_\epsilon}|\nabla w_\epsilon|^2\\
	 	&\leq\frac{1}{2}\int_\Omega\frac{|\nabla u_\epsilon|^2}{u_\epsilon}+\frac{\xi^2}{2}\int_\Omega\frac{u_\epsilon}{(u_\epsilon+1)^{2-2\theta}}|\nabla w_\epsilon|^2	\\
	 	&\leq\frac{1}{2}\int_\Omega\frac{|\nabla u_\epsilon|^2}{u_\epsilon}+\frac{\xi^2}{2}\int_\Omega|\nabla w_\epsilon|^2	\\
	 	&\leq\frac{1}{2}\int_\Omega\frac{|\nabla u_\epsilon|^2}{u_\epsilon}+C_1e^{-\tilde{\gamma}t}.
	 \end{split}
	 \ee
In view of \eqref{sharp-ck+}  in Lemma \ref{ext-ck}, a logarithmic-type Sobolev inequality follows:
	 \be\label{temp2}
	 	\int_\Omega u_\epsilon\ln\frac{u_\epsilon}{\bar{u}_0}\leq C_2\int_\Omega\frac{|\nabla u_\epsilon|^2}{u_\epsilon}.
	 \ee
Inserting  \eqref{temp1} into \eqref{h2-functional} and combining  \eqref{temp2} and  the definition of $\mathcal{\tilde{H}}_\epsilon$ in \eqref{h2-def}, we readily infer, with $\tilde A=\min\{\frac{1}{4C_2},1\}$,  that
	 \be\label{ulnu-decay}	\mathcal{\tilde{H}}_\epsilon^\prime(t)+\tilde A\mathcal{\tilde{H}}_\epsilon(t)+\frac{1}{4}\int_\Omega\frac{|\nabla u_\epsilon|^2}{u_\epsilon}+\frac{\chi}{2}\int_\Omega(|\Delta v_\epsilon|^2+|\nabla v_\epsilon|^2)\leq C_1e^{-\tilde{\gamma}t}.	 	
	\ee
Fixing
	\be\label{delta2-def}
		0<\tilde{\delta}<\min\{\tilde{A}, \ \  \tilde{\gamma}\}=\min\{\frac{1}{4C_2},\ \ 1, \ \ \tilde{\gamma}\},
	\ee
 then solving  the differential inequality \eqref{ulnu-decay} on $[0,t)$, we infer
	\be\label{decay1}
	\mathcal{\tilde{H}}_\epsilon(t)\leq \mathcal{\tilde{H}}_\epsilon(0)e^{-\tilde{A} t}+C_1
	\begin{cases} \frac{1}{\tilde  A-\tilde \gamma }\left(e^{-\tilde A t}-e^{-\tilde \gamma t}\right), &\text{if } \tilde{A}\neq\tilde\gamma,\\[0.2cm]
		te^{-\tilde A t}, &\text{if } \tilde{A}=\tilde\gamma,
	\end{cases}
	\nonumber
	\ee
which along with   the definition of $\mathcal{\tilde{H}}_\epsilon$ in \eqref{h2-def} implies
	\begin{equation*}
	\mathcal{\tilde{H}}_\epsilon(t)\leq e^{-\tilde \delta t}\int_\Omega\left(u_0\ln\frac{u_0}{\bar{u}_0}+\frac{\chi}{2}|\nabla v_0|^2\right)+C_1e^{-\tilde \delta t}
	\begin{cases} \frac{1}{|\tilde  A-\tilde \gamma| }, &\text{if } \tilde{A}\neq\tilde\gamma,\\[0.2cm]
		\frac{1}{e(\tilde A-\tilde\delta)}, &\text{if } \tilde{A}=\tilde\gamma,
	\end{cases}
	\end{equation*}
yielding directly \eqref{ulnu-epsilon2}.  Next, a simple substitution of  \eqref{temp1} into \eqref{h2-functional} gives
$$
\mathcal{\tilde{H}}_\epsilon^\prime(t)+\frac{1}{2}\int_\Omega\frac{|\nabla u_\epsilon|^2}{u_\epsilon}+\chi\int_\Omega(|\Delta v_\epsilon|^2+|\nabla v_\epsilon|^2)\leq C_1e^{-\tilde{\gamma}t}. 	 	
$$
Then we  integrate this inequality on $[0,t)$ and use \eqref{ulnu-bdd-key} to deduce that
	\begin{equation*}
\begin{split}
		&\int_{0}^t\int_\Omega\left(\frac{|\nabla u_\epsilon|^2}{2u_\epsilon}+\chi|\Delta v_\epsilon|^2+\chi|\nabla v_\epsilon|^2\right)\\
&\leq \frac{C_1}{\tilde \gamma}+\mathcal{\tilde{H}}_\epsilon(0)-\mathcal{\tilde{H}}_\epsilon(t)\\
&\leq \frac{C_1}{\tilde \gamma}+\int_\Omega\left(u_0\ln\frac{u_0}{\bar{u}_0}+\frac{\chi}{2}|\nabla v_0|^2\right)+(1+e^{-1} )\int_\Omega u_0,
	\end{split}
\end{equation*}
	which  upon sending $t\rightarrow \infty$ entails the desired space-time bound  \eqref{delta v-grad v-estimate-1}.	
\end{proof}
	
	With the help of  Lemma \ref{ulnu-epsilon-lemma}, we now prove Theorem \ref{theorem4} on the 3D global existence of   weak solutions by means of  compactness arguments via  Aubin-Lions lemma.

\begin{lemma}\label{weak-solution-proposition}
	Assume that the initial datum $(u_0,v_0,w_0)$ satisfies \eqref{initial-data} and $\theta<\frac{7}{6}$. Then the IBVP \eqref{main equation} possesses a global weak solution in the sense of Definition \ref{weak-definition} and it fulfills the regularities stated in \eqref{reg-weak1}.
\end{lemma}
\begin{proof}
	 With the preparations provided so far, we next  innovate the framework in  \cite[Lemma 4.1]{Cieslak-Laurencot-MM08} and \cite[Section 3.2]{Jiang22} to show  the global existence of  3D weak solution.  To this end, for any $\epsilon\in (0,   \|u_0\|_{L^\infty}^{-1})$, we see from Lemma \ref{regular-system-local-lemma} that the IBVP \eqref{regularized-equation} has a unique global classical solution.  For any $T>0$, since $\nabla u_\epsilon^\frac{1}{2}=\frac{\nabla u_\epsilon}{2\sqrt{u_\epsilon}}$ and $\|u_\epsilon(t)\|_{L^1}=\|u_0\|_{L^1}$ by \eqref{v-l1-estimate2}, we deduce from the space-time bound \eqref{delta v-grad v-estimate-2} that
	$$
		\int_{0}^{T}\|u_\epsilon^\frac{1}{2}\|_{W^{1,2}}^2=\int_0^T\int_\Omega \left(u_\epsilon+\frac{1}{4}\frac{|\nabla u_\epsilon|^2}{u_\epsilon}\right)\leq C_1(T),
	$$
	which along with the continuous embedding $W^{1,2}(\Omega)\hookrightarrow L^6(\Omega)$ with $\Omega\subset\mathbb{R}^3$ implies
	\be\label{e2}
	\int_{0}^{T}\|u_\epsilon\|_{L^3}=\int_{0}^{T}\|u_\epsilon^\frac{1}{2}\|_{L^6}^2\leq C_2\int_{0}^{T}\|u_\epsilon^\frac{1}{2}\|_{W^{1,2}}^2\leq  C_3(T).
    \ee
	In light of the H\"{o}lder interpolation inequality, it follows    that
	\be\label{e3}
		\int_{0}^{T}\|u_\epsilon\|_{L^\frac{3p}{3p-2}}^p\leq \int_0^T\|u_\epsilon\|_{L^1}^{p-1}\|u_\epsilon\|_{L^3}=\|u_0\|_{L^1}^{p-1}\int_0^T\|u_\epsilon\|_{L^3}\leq C_4(T),\ \ \forall p\in[1,\infty),
	\nonumber
	\ee
	and in particular, it follows upon setting $p=\frac{5}{3}$ that
	\be\label{e4}
		\int_{0}^{T}\int_\Omega u_\epsilon^\frac{5}{3}\leq C_5(T).
	\ee
	Therefore, applying H\"{o}lder's inequality, \eqref{delta v-grad v-estimate-2} and \eqref{e4}, we derive that
	\be\label{e5}
		\int_{0}^{T}\int_\Omega|\nabla u_\epsilon|^\frac{5}{4}\leq\left(\int_{0}^{T}\int_\Omega\frac{|\nabla u_\epsilon|^2}{u_\epsilon}\right)^\frac{5}{8}\left(\int_{0}^{T}\int_\Omega u_\epsilon^\frac{5}{3}\right)^\frac{3}{8}\leq C_6(T).
	\ee
	A combination of  \eqref{e4} and \eqref{e5} then shows that
	\be\label{bdd1}
		\{u_\epsilon\}_{\epsilon>0} \text{ are uniformly bounded in } L^\frac{5}{4}(0,T;W^{1,\frac{5}{4}}(\Omega))\cap L^\frac{5}{3}(\Omega\times(0,T)).
 	\ee
	Hence, by  the second equation in  \eqref{regularized-equation} and the  space-time bound of $u_\epsilon$ in \eqref{e4}, one can conclude from the maximal Sobolev regularity first
$$
		v_\epsilon\in W^{1,\frac{5}{3}}(0,T;L^\frac{5}{3}(\Omega))\cap L^\frac{5}{3}(0,T;W^{2,\frac{5}{3}}(\Omega)),
	$$
	and then
$$
		\int_{0}^{T}\int_\Omega|v_{\epsilon t}|^\frac{5}{3}+\int_{0}^{T}\int_\Omega|(-\Delta+1)v_\epsilon|^\frac{5}{3}\leq C_7(T).
	$$
	The above two estimates  together with elliptic estimate theory indicate that
	\be\label{bdd2}
		\{v_\epsilon\}_{\epsilon>0} \text{ are uniformly bounded in } L^\frac{5}{3}(0,T;W^{2,\frac{5}{3}}(\Omega))
	\nonumber
	\ee
	and
	\be\label{bdd3}
		\{v_{\epsilon t}\}_{\epsilon>0} \text{ are uniformly bounded in } L^\frac{5}{3}(\Omega\times(0,T)).
	\ee	
	On the other hand,  by \eqref{v-l1-estimate2}, \eqref{ulnu-epsilon} and \eqref{delta v-grad v-estimate-2}, one can deduce further  regularity estimate for $v_\epsilon$:
	\be\label{bdd4}
		\{v_\epsilon\}_{\epsilon>0} \text{ are uniformly bounded in } L^\infty(0,T;W^{1,2}(\Omega))\cap L^2(0,T;W^{2,2}(\Omega)).
	\ee
	Then \eqref{bdd4} along  with the continuous  embedding  $W^{2,2}(\Omega)\hookrightarrow  W^{1,6}(\Omega)$ with $\Omega \subset \mathbb{R}^3$ shows that $\nabla v_\epsilon$ are uniformly bounded in both $L^\infty(0,T;L^2(\Omega))$ and $L^2(0,T;L^6(\Omega))$, whence again by interpolation,
	\begin{equation*}
\begin{split}
		\int_{0}^{T}\|\nabla v_\epsilon\|_{L^\frac{6p}{3p-4}}^p&\leq 	\|\nabla v_\epsilon\|_{L^\infty(0,T; L^2(\Omega))}^{p-2}\int_0^T\|\nabla v_\epsilon\|_{L^6}^2\\
&\leq C_8\|\nabla v_\epsilon\|_{L^\infty(0,T; L^2(\Omega))}^{p-2}\int_0^T\| v_\epsilon\|_{W^{2,2}}^2 \leq C_9(T), \  \forall p\in[2,\infty),
\end{split}
\end{equation*}
	and in particular,   it follows upon setting respectively  $p=2$ and $p=\frac{10}{3}$ that
	\be\label{e8}
		\int_{0}^{T}\|\nabla v_\epsilon\|_{L^6}^2+\int_{0}^{T}\int_\Omega|\nabla v_\epsilon|^\frac{10}{3}\leq C_{10}(T),
	\ee
	which along with \eqref{e4} implies
	\be\label{bdd5}
		\{u_\epsilon\nabla v_\epsilon\}_{\epsilon>0} \text{ are uniformly bounded in } L^\frac{10}{9}(\Omega\times(0,T)).
	\ee
	Using   the nonnegativity of $v_\epsilon$ and the formula for  $w_\epsilon$  in \eqref{u-l1-estimate2}, we see  that
\be\label{gradw-eps}
\|\nabla w_\epsilon(t)\|_{L^q}\leq \|\nabla w_0\|_{L^q}+\|w_0\|_{L^\infty}\int_0^t\|\nabla v_\epsilon(s)\|_{L^q}, \ \  \  \forall q\in[1,\infty),
\ee
which, upon setting $q=6$ and   the use of  \eqref{v-l1-estimate2}, \eqref{ulnu-epsilon} and \eqref{e8},  shows that
	\be\label{w1}
\begin{cases}
		&\{w_\epsilon\}_{\epsilon>0} \text{ are uniformly bounded } \\[0.2cm]
 &  \text{ in  } L^\infty(0,T;W^{1,2}(\Omega))\cap L^2(0,T;W^{1,6}(\Omega)) \cap  L^\infty(\Omega\times(0,T)).
\end{cases}
	\ee
	Then applying \eqref{v-lq-norm-estimate1} of Lemma \ref{v-lq-norm-lemma} to the IBVP \eqref{regularized-equation} with $p=1,n=3$, we  get, for $1\leq r<3$, that
	\be\label{w2}
		\{w_{\epsilon t}\}_{\epsilon>0}, \{v_\epsilon w_\epsilon\}_{\epsilon>0} \text{ are uniformly bounded in } L^r(\Omega\times(0,T)).
	\ee
	In light of  \eqref{e8} and \eqref{gradw-eps} with $q=\frac{10}{3}$, one can derive that
	\be\label{w4}
		\int_{0}^{T}\int_\Omega|\nabla w_\epsilon|^\frac{10}{3}\leq C_{11}(T).
	\ee
	Noticing that $0\leq S(u_\epsilon)\leq  u_\epsilon(1+u_\epsilon)^{\theta-1}$, we utilize  \eqref{e4} and \eqref{w4} to obtain that
	\be\label{w3}
		\{S(u_\epsilon)\nabla w_\epsilon\}_{\epsilon>0} \text{ are uniformly bounded in } L^{\frac{10}{3+6\theta^+}}(\Omega\times(0,T)).
	\ee	
	Next, for  $\phi\in W^{1,\infty}(\Omega)$, it follows from the $u_\epsilon$-equation in  \eqref{regularized-equation} and \eqref{ulnu-epsilon} that
	\begin{equation*}
	\begin{split}
		\bigg|\int_\Omega u_{\epsilon t}\phi\bigg|&\leq\int_\Omega|\nabla u_\epsilon||\nabla\phi|+\chi\int_\Omega u_\epsilon(1-\epsilon u_\epsilon)|\nabla v_\epsilon||\nabla\phi|	\\
		&\ \ \ +\xi\int_\Omega S(u_\epsilon)(1-\epsilon u_\epsilon)|\nabla w_\epsilon||\nabla\phi|	\\
		&\leq \|\nabla\phi\|_{L^\infty}\|\nabla u_\epsilon\|_{L^\frac{5}{4}}|\Omega|^\frac{1}{5}+\chi\|u_\epsilon\|_{L^2}\|\nabla v_\epsilon\|_{L^2}\|\nabla\phi\|_{L^\infty}	\\
		&\ \ \ +\xi\|\nabla\phi\|_{L^\infty}\| S(u_\epsilon)\nabla w_\epsilon\|_{L^1}	\\
		&\leq C_{12}\left(\|\nabla u_\epsilon\|_{L^\frac{5}{4}}+\|u_\epsilon\|_{L^2}+\| S(u_\epsilon)\nabla w_\epsilon\|_{L^1}\right)\|\nabla\phi\|_{L^\infty},
	\end{split}
	\end{equation*}
	which coupled with  \eqref{e2}, \eqref{e5} and \eqref{w3} with $\theta<\frac{7}{6}$ gives rise to
	\be\label{bdd6}
		\{u_{\epsilon t}\}_{\epsilon>0} \text{ are uniformly bounded in } L^1(0,T;(W^{1,\infty}(\Omega))^*).
	\ee
	Thanks to  \eqref{bdd1} and \eqref{bdd6}, by the compact embedding $W^{1,\frac{5}{4}}(\Omega)\hookrightarrow L^1(\Omega)$ and  the continuous embedding $L^1(\Omega) \hookrightarrow (W^{1,\infty}(\Omega))^*$, we conclude from the well-known Aubin-Lions lemma (c.f  \cite{Simon87}) that $u_\epsilon$ are relatively compact in $L^1(\Omega\times(0,T))$, which can be strengthen by \eqref{e4} and interpolation that
	\be\label{e10}
		\{u_\epsilon\}_{\epsilon>0} \text{ are relatively compact in } L^p(\Omega\times(0,T)),\quad\forall p\in[1,\frac{5}{3}).
	\ee
	Notice that $\|u_\epsilon\|_{L^1}=\|u_0\|_{L^1}$, based on \eqref{e10} and  \eqref{bdd6}, an application of the Arzel\`{a}-Ascoli compactness theorem shows, as $\epsilon\searrow 0$ (without relabeling for convenience), there exists some nonnegative function $u$ such that, for any $p\in[1,\frac{5}{3})$,
	\be\label{converge1}
		u_\epsilon\rightarrow u \text{ in } C([0,T];(W^{1,\infty}(\Omega))^*)\cap L^p(\Omega\times(0,T)) \text{ and a.e. in } \Omega\times (0, T).
	\ee

 Similarly,  by the compact embedding $W^{1,2}(\Omega)\hookrightarrow L^2(\Omega)$ and  the continuous embedding $L^2(\Omega) \hookrightarrow L^\frac{5}{3}(\Omega)$, from  \eqref{bdd3}, \eqref{bdd4} and \cite[Corollary 4]{Simon87},  we infer, for some nonnegative function  $v$,  that
	\be\label{converge2}
		v_\epsilon\rightarrow v \text{ in } C([0,T];L^2(\Omega))\cap L^2(0,T;W^{1,2}(\Omega)) \text{ and a.e. in } \Omega\times (0, T),
	\ee
	and from \eqref{w1}, \eqref{w2} and again \cite[Corollary 4]{Simon87}, for some nonnegative function  $w$,  that
	\be\label{converge3}
		w_\epsilon\rightarrow w \text{  in } C([0,T];L^2(\Omega))\cap L^2(0,T;W^{1,2}(\Omega)) \text{ and a.e. in } \Omega\times (0, T).
	\ee	
	Then  for any $\varphi\in W^{1,\infty}(\Omega)$, multiplying  by  $\varphi$ the $ v_\epsilon$-equation and $w_\epsilon$-equation in  \eqref{regularized-equation} respectively and integrating by parts  over $\Omega\times (0,t)$, by the convergence properties \eqref{converge1}  \eqref{converge2} and  \eqref{converge3}, up to a subsequence of $\epsilon\rightarrow0$, we  obtain that
	\be\label{weak-v}
		\int_\Omega(v(t)-v_0)\varphi+\int_{0}^{t}\int_\Omega(\nabla v\cdot\nabla\varphi+(v-u)\varphi)=0,\quad\forall t\in[0,T]
	\ee
	and
	\be\label{weak-w}
		\int_\Omega(w(t)-w_0)\varphi+\int_{0}^{t}\int_\Omega vw\varphi=0,\quad\forall t\in[0,T].
	\ee	
Note that the uniform boundedness in \eqref{bdd1} implies the weak convergence:
 \be\label{bdd1+}
		\nabla u_\epsilon \rightharpoonup \nabla u\text{ in } \left(L^\frac{5}{4}(\Omega \times(0,T))\right)^3.
 	\ee
	In view of \eqref{bdd5} and the fact $0<1-\epsilon u_\epsilon<1$, we know that $\{u_\epsilon(1-\epsilon u_\epsilon)\nabla v_\epsilon\}_{\epsilon>0}$ are also bounded in $L^\frac{10}{9}(\Omega\times(0,T))$ and thus weakly compact and  uniformly integrable in $L^1(\Omega\times(0,T))$. Noting from the facts  \eqref{converge1} and \eqref{converge2} (up to  subsequence) that
	$$
		u_\epsilon(1-\epsilon u_\epsilon)\nabla v_\epsilon\longrightarrow u\nabla v \text{ a.e. in } \Omega\times(0,T),
	$$
	one can readily infer  from  the Vitali convergence theorem (roughly, uniform integrability plus a.e. convergence imply convergence) that
	\be\label{converge4}
		u_\epsilon(1-\epsilon u_\epsilon)\nabla v_\epsilon\longrightarrow u\nabla v \text{ in } \left(L^1(\Omega\times(0,T))\right)^3.
	\ee
	On the other hand, it follows from \eqref{w3} that  $\{S(u_\epsilon)(1-\epsilon u_\epsilon)\nabla w_\epsilon\}_{\epsilon>0}$ are bounded in $L^\frac{10}{3+6\theta^+}(\Omega\times(0,T))$, and thus weakly compact and uniformly integrable in $L^1(\Omega\times(0,T))$ (since  $\frac{10}{3+6\theta^+}>1$ due to $\theta<\frac{7}{6}$). Again by \eqref{converge1}, \eqref{converge3} and the Vitali theorem, we obtain
	\be\label{converge5}
		S(u_\epsilon)(1-\epsilon u_\epsilon)\nabla w_\epsilon\longrightarrow S(u)\nabla w \text{ in } \left(L^1(\Omega\times(0,T))\right)^3.
	\ee	
	Now,  we pass to the limit as $\epsilon\rightarrow0$ in the integral version of the  $u_\epsilon$-equation in \eqref{regularized-equation} to derive from \eqref{bdd1+}, \eqref{converge4} and \eqref{converge5}, for any $\varphi\in W^{1,\infty}(\Omega)$, that
	\be\label{weak-u}
		\int_\Omega(u(t)-u_0)\varphi+\int_{0}^{t}\int_\Omega(\nabla u+\chi u\nabla v-\xi S(u)\nabla w)\cdot\nabla\varphi=0,\quad\forall t\in[0,T],		
	\nonumber
	\ee
	which along with \eqref{weak-v} and \eqref{weak-w} readily shows that $(u,v,w)$ is a global weak solution of the IBVP \eqref{regularized-equation} in the sense of Definition \ref{weak-definition}.

The stated intermediate regularities  in \eqref{reg-weak1} follow from the lower   semi-continuity of norms (Fatou's lemma) and the established uniform estimates for the approximating solutions $(u_\epsilon, v_\epsilon, w_\epsilon)$ in \eqref{bdd1}, \eqref{bdd4} with the embedding  $W^{2,2}(\Omega)\hookrightarrow W^{1,6}(\Omega)$, and \eqref{w1}. We thus  complete the proof.	
\end{proof}	

\subsection{Eventual regularity and exponential convergence of weak solutions}
	
	In this subsection, we examine the eventual regularity and exponential stability of weak  solutions based on a linearized system of the approximating model \eqref{regularized-equation} and the stability of the homogeneous solution $(\bar{u}_0,\bar{u}_0,0)$ in scaling-invariant spaces.
	
\subsubsection{The linearized system}
	
To analyze  decay properties of the linearized system around $(\bar{u}_0,\bar{u}_0,0)$, we  introduce $\rho_\epsilon:=u_\epsilon-\bar{u}_0,c_\epsilon:=v_\epsilon-\bar{u}_0$, upon some manipulations,  we find that the  approximating system \eqref{regularized-equation} can be rewritten as
	\be\label{differ-system}
	\begin{cases}
	 \rho_{\epsilon t}-\Delta \rho_\epsilon-\alpha\Delta c_\epsilon+\beta\Delta
  w_\epsilon=\chi\nabla\cdot J_{1\epsilon}-\xi\nabla\cdot J_{2\epsilon} &\text{in } \Omega\times(0,\infty), \\[0.2cm]
		  c_{\epsilon t}-\Delta  c_\epsilon+ c_\epsilon- \rho_\epsilon=0  &\text{in } \Omega\times(0,\infty),\\[0.2cm]
		  w_{\epsilon t}+\bar u_0 w_\epsilon=-c_\epsilon w_\epsilon	 &\text{in } \Omega\times(0,\infty),\\[0.2cm]
		\frac{\partial  \rho_\epsilon}{\partial\nu}=\frac{\partial c_\epsilon}{\partial\nu}=\frac{\partial  w_\epsilon}{\partial\nu}=0  &\text{on } \partial\Omega\times(0,\infty),\\[0.2cm]
		( \rho_\epsilon, \   c_\epsilon, \   w_\epsilon)=(u_0-\bar u_0, \ v_0-\bar u_0, w_0)=:(\rho_I,\  c_I, \  w_I) &\text{on } \Omega\times\{0\},
	\end{cases}
	\ee
where the parameters $\alpha$ and $\beta$ and fluxes  $J_{1\epsilon}$ and $J_{2\epsilon}$ are defined by
\be\label{tranfs-J1J2-etc}
\begin{cases}
\alpha=\chi\bar{u}_0 (1-\epsilon\bar{u}_0), \ \ \  \beta=\xi S(\bar{u}_0) (1-\epsilon\bar{u}_0), \\[0.2cm]
J_{1\epsilon}=\left[1-\epsilon(\rho_\epsilon+\bar{u}_0)-\epsilon\bar{u}_0\right]\rho_\epsilon\nabla c_\epsilon, \\[0.2cm]
J_{2\epsilon}=\left[(S(\rho_\epsilon+\bar{u}_0)-S(\bar{u}_0))(1-\epsilon\bar{u}_0)-\epsilon \rho_\epsilon S(\rho_\epsilon+\bar{u}_0)\right]\nabla w_\epsilon.
\end{cases}
\ee
Then its  corresponding linearized system around the origin reads as follows:
	\be\label{linear-system}\begin{cases}
		\tilde\rho_{\epsilon t}-\Delta \tilde\rho_\epsilon-\alpha\Delta \tilde c_\epsilon+\beta\Delta \tilde w_\epsilon=0  &\text{in } \Omega\times(0,\infty), \\[0.2cm]
		\tilde c_{\epsilon t}-\Delta\tilde c_\epsilon+\tilde c_\epsilon-\tilde\rho_\epsilon=0  &\text{in } \Omega\times(0,\infty),\\[0.2cm]
		\tilde w_{\epsilon t}+\bar{u}_0\tilde w_\epsilon=0 	 &\text{in } \Omega\times(0,\infty),\\[0.2cm]
		\frac{\partial \tilde\rho_\epsilon}{\partial\nu}=\frac{\partial\tilde c_\epsilon}{\partial\nu}=\frac{\partial\tilde w_\epsilon}{\partial\nu}=0  &\text{on } \partial\Omega\times(0,\infty), \\[0.2cm]
		(\tilde\rho_\epsilon, \tilde c_\epsilon, \tilde w_\epsilon)=(\rho_I,c_I,w_I)  &\text{on } \Omega\times\{0\},
	\end{cases}
	\ee
and, by introducing $(\hat\rho_\epsilon, \hat c_\epsilon, \hat w_\epsilon):=(\rho_\epsilon-\tilde\rho_{\epsilon }, c_\epsilon-\tilde c_\epsilon, w_\epsilon-\tilde w_\epsilon)$, and then subtracting \eqref{linear-system}  from \eqref{differ-system}, we find that the nonlinear part of \eqref{differ-system} reads as follows:
\be\label{difference-system}
	\begin{cases}
		\hat \rho_{\epsilon t}-\Delta \hat\rho_\epsilon-\alpha\Delta \hat c_\epsilon+\beta\Delta
\hat w_\epsilon =\chi\nabla\cdot J_{1\epsilon}-\xi\nabla\cdot J_{2\epsilon} &\text{in } \Omega\times(0,\infty), \\[0.2cm]
		\hat c_{\epsilon t}-\Delta\hat c_\epsilon+\hat c_\epsilon-\hat\rho_\epsilon=0  &\text{in } \Omega\times(0,\infty),\\[0.2cm]
		\hat w_{\epsilon t}+\bar u_0\hat w_\epsilon=-c_\epsilon w_\epsilon	 &\text{in } \Omega\times(0,\infty),\\[0.2cm]
		\frac{\partial \hat\rho_\epsilon}{\partial\nu}=\frac{\partial\hat c_\epsilon}{\partial\nu}=\frac{\partial\hat w_\epsilon}{\partial\nu}=0  &\text{on } \partial\Omega\times(0,\infty),\\[0.2cm]
		(\hat\rho_\epsilon, \  \hat c_\epsilon, \ \hat w_\epsilon)=(0,\  0, \  0) &\text{on } \Omega\times\{0\}.
	\end{cases}
	\ee
Notice   that $w_\epsilon$, $\tilde w_\epsilon$  and thus $\hat w_\epsilon$ can be expressed explicitly: indeed,
\be\label{tilde w3}
		 \left(w_\epsilon, \ \  \tilde{w}, \ \  \hat w_\epsilon)\right (t)=w_Ie^{-\bar{u}_0t}\left(e^{-\int_0^t c_\epsilon(s)ds}, \  \ 1, \ \  e^{-\int_0^t c_\epsilon(s)ds}-1\right).
	\ee
 In the sequel, we shall combine energy estimates and operator semigroup estimates to derive solution estimates of \eqref{differ-system} via \eqref{linear-system} and \eqref{difference-system}. For this purpose, for $p\in (1,\infty)$, we define the realization   $\mathcal{A}:  D(\mathcal{A})\subset (L^p(\Omega))^2\mapsto (L^p(\Omega))^2$ by
\be\label{Amatrix-def}
		\mathcal{A}(\rho, c)=\begin{pmatrix}
		\Delta & \alpha\Delta \\
		1 & \Delta-1
		\end{pmatrix}\begin{pmatrix}\rho\\
c
\end{pmatrix},\ee
where the domain  $D(\mathcal{A})=W^{2,p}_{N, 0}(\Omega)\times W^{2,p}_N(\Omega)$  with $W^{2,p}_N(\Omega)=\{z\in W^{2,p}(\Omega): \ \nabla z\cdot \nu=0 \text{ on }\partial \Omega\}$ and $W^{2,p}_{N, 0}(\Omega)=\{z\in W_N^{2,p}(\Omega): \ \bar z=0\}$.

Now, if   $\mathcal{A}(\rho, c)=(0, 0)$  in $D(\mathcal{A})$, then we quickly see  that   $(\rho,c)= (0, 0)$, which shows that zero is not the  eigenvalue of $\mathcal{A}$. Moreover,  by  eigen-expansion (cf. \cite[Proposition 2.2]{XT15JDE}), we find that  the spectrum  of  $-\mathcal{A}$ is  given by
 $$
 \sigma(-\mathcal{A})=  \left\{\frac{1+2\lambda_k+\sqrt{1-4\alpha\lambda_k}}{2}\right\}_{k=1}^\infty\bigcup
 \left\{\frac{1+2\lambda_k-\sqrt{1-4\alpha\lambda_k}}{2}\right\}_{k=1}^\infty,
 $$
 where $\{\lambda_k\}_{k=1}^\infty$ is the complete increasing sequence of eigenvalues of the realization of the Neumann $-\Delta$   in $L_0^p(\Omega)$ with $0<\lambda_1\leq \lambda_2\leq\cdots$ and $\lambda_k\rightarrow \infty$ as $k\rightarrow\infty$. Here   $L^p_0(\Omega)=\{z\in L^p(\Omega): \ \bar z=0\}$.

Then, noting from \eqref{tranfs-J1J2-etc} that  $\alpha=\chi\bar{u}_0 (1-\epsilon\bar{u}_0)$, we apparently deduce that
\be\label{A-leg}
\inf\left\{\text{Re}\lambda: \lambda\in \sigma(-\mathcal{A})\right\}\geq \lambda_1>0, \ \  \forall \epsilon\in [0, \   \frac{1}{\bar u_0}].
\ee
Consequently, $\mathcal{A}$ is sectorial (see also   \cite[Section 4.1]{Jiang22}), and generates an analytical semigroup on $L_0^p(\Omega)\times L^p(\Omega)$.  Hence, using variation-of-constant-formula to the first two equations in \eqref{differ-system}, the $(\rho_\epsilon, c_\epsilon)$-solution component of \eqref{differ-system} can be written through \eqref{linear-system} and \eqref{difference-system} as
 \be\label{rho-epsilon-rep+}
 \begin{split}
	\begin{pmatrix}
		\rho_\epsilon(t) \\ c_\epsilon(t)
	\end{pmatrix}&=e^{t\mathcal{A}}\begin{pmatrix}
		\rho_I \\ c_I
	\end{pmatrix}+\int_{0}^{t} e^{(t-s)\mathcal{A}}
	\begin{pmatrix}
	-\beta\Delta w_\epsilon+\chi\nabla\cdot J_{1\epsilon}-\xi\nabla\cdot J_{2\epsilon} \\ 0
	\end{pmatrix}(s)ds\\[0.25cm]
&\ \ =\begin{pmatrix}
		\tilde\rho_\epsilon(t) \\ \tilde c_\epsilon(t)
	\end{pmatrix}+\begin{pmatrix}
		\hat\rho_\epsilon(t)  \\  \hat c_\epsilon(t)
	\end{pmatrix}.
\end{split} 		\ee
 Henceforth, for convenience, due to \eqref{A-leg},   we consider $\epsilon\in [0, \frac{1}{2\bar u_0}]$ so that $\alpha$ and $\beta$ are bounded from below and above  uniformly in $\epsilon$,  and thus, we shall  omit the subscript $\epsilon$ within this subsection since the involving estimates are independent of such $\epsilon$. In the following, we shall also use the notation like $C_N^2(\bar{\Omega})=\{z\in C^2(\bar{\Omega}): \nabla w_I\cdot\nu=0 \text{ on }   \partial\Omega\}$. With theses preparations, we next  derive some main  $\epsilon$-independent estimates for  solutions to   \eqref{linear-system} and \eqref{difference-system}.

\begin{lemma}\label{rho-l2-grad c-l2-lemma}
	 Let the initial datum $(\rho_I,c_I,w_I)\in L_0^2(\Omega)\times H^1(\Omega)\times C_N^2(\bar{\Omega})$ and let $0<\mu<\min\{\lambda_1,2\bar{u}_0\}$. Then the solution of  the linearized system  \eqref{linear-system} decays  exponentially: there exists $C>0$ such that, for all $t\geq0$,
	 \be\label{rho-l2-grad c-l2-norm}
	 	\|\tilde{\rho}(t)\|_{L^2}^2+\alpha\|\nabla\tilde{c}(t)\|_{L^2}^2\leq Ce^{-\mu t}\left(\|\rho_I\|_{L^2}^2+\alpha\|\nabla c_I\|_{L^2}^2+\|\nabla w_I\|_{L^2}^2\right).
	 \ee
\end{lemma}	
\begin{proof}
	Integrating by parts over $\Omega$ simply in the $\tilde{\rho}$-equation of  \eqref{linear-system} shows that $\int_\Omega \tilde{\rho}=\int_\Omega \rho_I=0$,  and observing  the expression for  $\tilde{w}$ in \eqref{tilde w3}  that
	\be\label{tilde w}
		\tilde{w}(t)=w_Ie^{-\bar{u}_0t},\quad\forall t\geq0.
	\ee
	Multiplying the $\tilde{\rho}$-equation of \eqref{linear-system} by $\tilde{\rho}$ and the $\tilde{c}$-equation by $-\alpha\Delta\tilde{c}$ respectively and then integrating by parts over $\Omega$, we derive that
	\be\label{energy-function}
		\frac{1}{2}\frac{d}{dt}\int_\Omega\left(\tilde{\rho}^2+\alpha|\nabla\tilde{c}|^2\right)+\int_\Omega|\nabla\tilde{\rho}|^2+\alpha\int_\Omega\left(|\Delta\tilde{c}|^2+|\nabla\tilde{c}|^2\right)=\beta\int_\Omega\nabla\tilde{\rho}\cdot\nabla\tilde{w}.
	\ee
	Applying Young's inequality and \eqref{tilde w}, we estimate that
	\be\label{term-1}
		\beta\int_\Omega\nabla\tilde{\rho}\cdot\nabla\tilde{w}=\beta e^{-\bar{u}_0t}\int_\Omega\nabla\tilde{\rho}\cdot\nabla w_I\leq\frac{1}{2}\int_\Omega|\nabla\tilde{\rho}|^2+\frac{\beta^2}{2}e^{-2\bar{u}_0t}\int_\Omega|\nabla w_I|^2.
	\ee
Since $\bar{\tilde{\rho}}=0$, we employ  H\"{o}lder's inequality and  Poincar\'{e}'s inequality to infer that
	\be\label{term-2}
	\begin{split}
		\|\nabla\tilde{c}\|_{L^2}^2=-\int_\Omega(\tilde{c}-\bar{\tilde{c}})\Delta\tilde{c}\leq\|\tilde{c}-\bar{\tilde{c}}\|_{L^2}\|\Delta\tilde{c}\|_{L^2}\leq\frac{1}{\sqrt{\lambda_1}}\|\nabla\tilde{c}\|_{L^2}\|\Delta\tilde{c}\|_{L^2}
	\end{split}
	\ee
	and that
	\be\label{term-3}
		\lambda_1\|\tilde{\rho}\|_{L^2}^2\leq\|\nabla\tilde{\rho}\|_{L^2}^2.
	\ee
	Substituting \eqref{term-1}, \eqref{term-2} and \eqref{term-3} into \eqref{energy-function} simply gives rise to
	$$
		\frac{d}{dt}\int_\Omega\left(\tilde{\rho}^2+\alpha|\nabla\tilde{c}|^2\right)+
\lambda_1\int_\Omega\tilde{\rho}^2+2\alpha(1+\lambda_1)\int_\Omega|\nabla\tilde{c}|^2\leq
\beta^2e^{-2\bar{u}_0t}\int_\Omega|\nabla w_I|^2.	
	$$
	Solving this differential inequality, we finally arrive at the estimate \eqref{rho-l2-grad c-l2-norm}.	
\end{proof}
	
	Based on the well-known smoothing $L^p$-$L^q$ estimates for semigroup $(e^{t\Delta})_{t\geq0}$ and Lemma \ref{rho-l2-grad c-l2-lemma}, we next establish some further $L^p$-bounds of $\tilde{\rho}$ and $\nabla\tilde{c}$.
	
\begin{lemma}\label{rho-grad c-lp-lemma}
	 Let the initial datum $(\rho_I,c_I,w_I)\in (C(\bar{\Omega})\cap L_0^1(\Omega))\times W_N^{1,\infty}(\Omega)\times C_N^2(\bar{\Omega})$  and let $0<\mu<\min\{\lambda_1,\bar{u}_0\}$, then for some $C>0$ and for all $t>0$, there holds, for any $\frac{n}{2}\leq p<\infty$, that
	\be\label{rho-lp}
		\|\tilde{\rho}(t)\|_{L^p}\leq Ce^{-\mu t}\left(1+t^{-1+\frac{n}{2p}}\right)\left(\|\rho_I\|_{L^\frac{n}{2}}+\|\nabla c_I\|_{L^n}+\|\nabla w_I\|_{L^n}\right),
	\ee
	and for any $n\leq p<\infty$, that
	\be\label{grad c-lp}
		\|\nabla\tilde{c}(t)\|_{L^p}\leq Ce^{-\mu t}\left(1+t^{-\frac{1}{2}+\frac{n}{2p}}\right)\left(\|\rho_I\|_{L^\frac{n}{2}}+\|\nabla c_I\|_{L^n}+\|\nabla w_I\|_{L^n}\right).
	\ee
\end{lemma}
\begin{proof}
	For all $t>0$, the variation-of-constants formula applied to the $\tilde{\rho}$- and $\tilde{c}$-equation in \eqref{linear-system} show  that
	\be\label{rho-representation}
		\tilde\rho(t)=e^{t\Delta}\rho_I+\alpha\int_{0}^{t}e^{(t-s)\Delta}\Delta\tilde{c}(s)ds-\beta\int_{0}^{t}e^{(t-s)\Delta}\Delta\tilde{w}(s)ds
	\ee
	and that
	\be\label{c-representation}
		\tilde c(t)=e^{t(\Delta-1)}c_I+\int_{0}^{t}e^{(t-r)(\Delta-1)}\tilde{\rho}(r)dr.
	\ee
	We plug  \eqref{tilde w} and \eqref{c-representation} into \eqref{rho-representation} to see that
	\be\label{rho2}
	\begin{split}
		\tilde{\rho}(t)=&e^{t\Delta}\rho_I+\alpha\int_{0}^{t}e^{(t-s)\Delta}\Delta e^{s(\Delta-1)}c_Ids	\\
		&+\alpha\int_{0}^{t}e^{(t-s)\Delta}\Delta \int_{0}^{s}e^{(s-r)(\Delta-1)}\tilde{\rho}(r)drds-\beta \int_{0}^{t}e^{-\bar{u}_0s}e^{(t-s)\Delta}\Delta w_Ids.
	\end{split}	
	\nonumber	
	\ee
 With the exponential decay in \eqref{rho-l2-grad c-l2-norm}, the $L^p$-bound of the first three terms on the right-hand side of the above identity can be treated in the same  way as in \cite[Lemma 3.5]{Jiang20} without smallness condition on $\alpha$, thus we omit the details here. As for the last term, applying the heat semigroup estimates, we will show, for any $1<p<\infty$, $n\geq2$ and $0<\mu<\min\{\lambda_1,\bar{u}_0\}$, that
	\be\label{term-0}
		\bigg\|\beta\int_{0}^{t} e^{-\bar{u}_0s}e^{(t-s)\Delta}\Delta w_Ids\bigg\|_{L^p}\leq C_1e^{-\mu t}\|\nabla w_I\|_{L^n}.
	\ee
	Indeed, if $n\leq p<\infty$, we infer by the semigroup estimates that
	\be\label{term01}
	\begin{split}
		&\bigg\|\beta\int_{0}^{t} e^{-\bar{u}_0s}e^{(t-s)\Delta}\Delta w_Ids\bigg\|_{L^p}	\\
		&\leq C_2\int_{0}^{t}e^{-\bar{u}_0s}\left(1+(t-s)^{-\frac{1}{2}-\frac{n}{2}(\frac{1}{n}-\frac{1}{p})}\right)e^{-\lambda_1(t-s)}\|\nabla w_I\|_{L^n}ds	\\
		&\leq C_3e^{-\mu t}\|\nabla w_I\|_{L^n}\int_{0}^{t}\left(1+(t-s)^{-1+\frac{n}{2p}}\right)e^{-(\lambda_1-\mu)(t-s)}ds 	\\
		&\leq C_4e^{-\mu t}\|\nabla w_I\|_{L^n}.
	\end{split}
	\ee
	On the other hand, if $n>p$, we apply H\"{o}lder's inequality and  the well-known  semigroup estimate again to deduce that
	\be\label{term02}
	\begin{split}
		&\bigg\|\beta\int_{0}^{t} e^{-\bar{u}_0s}e^{(t-s)\Delta}\Delta w_Ids\bigg\|_{L^p}	\\
		&\leq \beta|\Omega|^{\frac{1}{p}-\frac{1}{n}}\int_{0}^{t}e^{-\bar{u}_0s}\|e^{(t-s)\Delta}\Delta w_I\|_{L^n}ds	\\
		&\leq C_5e^{-\mu t}\int_{0}^{t}\left(1+(t-s)^{-\frac{1}{2}}\right)e^{-(\lambda_1-\mu)(t-s)}\|\nabla w_I\|_{L^n}ds 	\\
		&\leq C_6e^{-\mu t}\|\nabla w_I\|_{L^n},
	\end{split}
	\nonumber
	\ee
	which together with \eqref{term01} gives rise to  \eqref{term-0}.
	
	Then one can proceed in the same manner as done in \cite[Lemma 3.5]{Jiang20} via   semigroup estimates to obtain the desired estimates based on  \eqref{term-0} and Lemma \ref{rho-l2-grad c-l2-lemma}.
\end{proof}

	In light of \eqref{A-leg}, the following lemma can also be  similarly derived as in \cite[Lemma 3.6]{Jiang20} (cf. also \cite[Lemma 4.4]{Jiang19}).
	
\begin{lemma}\label{tilde rho-lp-lemma}
	Assume that $\rho_I=\nabla\cdot z_I$, $c_I=0$ and  $w_I=0$.  Then there holds, for some $C>0$ and $\frac{n}{2}<q\leq p<\infty$, that
	$$
		\|\tilde{\rho}(t)\|_{L^p}\leq Ce^{-\lambda_1 t}\left(1+t^{-\frac{1}{2}-\frac{n}{2}(\frac{1}{q}-\frac{1}{p})}\right)\|z_I\|_{L^q}.
	$$
\end{lemma}

\subsubsection{Stabilization in scaling-invariant spaces.}

	With the $\epsilon$-independent estimates provided in Subsection 6.1, we next show a stability result of the approximating system \eqref{differ-system} based on a loop-type argument illustrated in \cite[Proposition 4.1]{Jiang22}.

\begin{lemma}\label{proposition}
	Let $\chi,\xi >0$, $S$ satisfy \eqref{S-con} with the derivative $S^\prime$ being bounded on $\mathbb{R}^+$,  and let the initial datum $(\rho_I,c_I,w_I)\in(C(\bar{\Omega})\cap L_0^1(\Omega))\times W_N^{1,\infty}(\Omega)\times C_N^2(\bar{\Omega})$, $\rho_I\geq-\bar{u}_0$, $ c_I>-\bar{u}_0$ in $\bar{\Omega}$. Then for any fixed $n<q_0<2n$ and $q_0\leq p_0<\frac{nq_0}{q_0-n}$,  there exists $\eta_0>0$ such that, for $0<\eta<\eta_0$, whenever
\be\label{initial-eps-small}
\|\rho_I\|_{L^\frac{n}{2}}+\|\nabla c_I\|_{L^n}+\|c_I\|_{L^\infty}+\|\nabla w_I\|_{L^{p_0}}+\|w_I\|_{L^\infty}\leq\eta,
 \ee
 then the solution $(\rho_\epsilon,c_\epsilon, w_\epsilon)$  of the approximating system \eqref{differ-system} with $\epsilon\in(0,\epsilon_0)$  is globally bounded  and satisfies, for  $0<\mu<\min\{\lambda_1, \bar{u}_0, v_*\}$ with some $v_*>0$ independent of $\epsilon$, being the largest lower bound of $c_\epsilon+\bar u_0$, that
	\be\label{rho-rhoepsilon}
		\|\rho_\epsilon(t)-\tilde{\rho}_\epsilon(t)\|_{L^{q_0}}\leq\eta e^{-\mu t}\left(1+t^{-1+\frac{n}{2q_0}}\right),\quad\forall t>0,
	\ee
	where $\epsilon_0=\min\{\frac{1}{2\bar{u}_0},\frac{1}{\|\rho_I+\bar{u}_0\|_{L^\infty}}\}$ and $(\tilde{\rho}_\epsilon,\tilde{c}_\epsilon,\tilde{w}_\epsilon)$ is the solution of the linearized system \eqref{linear-system}.
\end{lemma}
\begin{proof}
	Recall that $(u_\epsilon, v_\epsilon, w_\epsilon)=(\rho_\epsilon+\bar u_0, c_\epsilon+\bar u_0,w_\epsilon)$ and that  $0<\epsilon<\epsilon_0\leq \|\rho_I+\bar u_0\|_{L^\infty}^{-1}$,  the  global existence in Lemma \ref{regular-system-local-lemma} shows that the IBVP \eqref{regularized-equation} admits a unique global classical solution on $[0,\infty)$. Thus,   the the approximating system \eqref{differ-system} admits a unique global classical solution $(\rho_\epsilon,c_\epsilon,w_\epsilon)=(u_\epsilon-\bar{u}_0,v_\epsilon-\bar{u}_0,w_\epsilon)$ with   $\|u_\epsilon(t)\|_{L^\infty}<\frac{1}{\epsilon}$ on $\bar\Omega \times  [0,\infty)$, and then,    the $(W^{1,\infty}, W^{1,\infty})$-boundedness of $(c_\epsilon, w_\epsilon)$ follows readily.

 Notice that $v_\epsilon(0)=c_I+\bar u_0>0$, Lemma \ref{gradw-l2-decay} implies, for some $v_*>0$ independent of $\epsilon$, that
	\be\label{c-lower-bound}
		c_\epsilon(x,t)+\bar{u}_0=v_\epsilon(x,t)\geq v_*>0,\ \ \ \forall (x,t)\in \bar\Omega\times[0, \infty).
	\ee
 In order to conduct on a loop argument to show \eqref{rho-rhoepsilon},  with $q_0, \eta$ given as in the lemma,  for any given $\epsilon\in(0,\epsilon_0)$ and $0<\mu<\min\{\lambda_1,\bar{u}_0,v_*\}$, we define
	\be\label{T-def}
		T_\epsilon:=\sup\left\{T>0: \|\rho_\epsilon(t)-\tilde{\rho}_\epsilon (t)\|_{L^{q_0}}\leq\eta e^{-\mu t}(1+t^{-1+\frac{n}{2q_0}}), \ \forall t\in(0,T)\right\}.
	\ee
	First, the fact $\rho_\epsilon(0)-\tilde{\rho}_\epsilon (0)=0$ and the continuity of
	$\rho_\epsilon(t)-\tilde{\rho}_\epsilon(t)$ implies that $T_\epsilon$ is well-defined and $T_\epsilon>0$. In the sequel, we use connectedness (openness and closedness)  argument to show that if $\eta_0$ is sufficiently small, then $T_\epsilon=\infty$.
	
	To begin with, we deduce from \eqref{rho-lp} along with \eqref{initial-eps-small} and \eqref{T-def} that
	\be\label{rho-epsilon-q0}
		\|\rho_\epsilon(t)\|_{L^{q_0}}\leq C_1\eta e^{-\mu t}\left(1+t^{-1+\frac{n}{2q_0}}\right)
	\ee
	for all  $t\in (0, T_\epsilon)$, which will be assumed till the end proof to avoid unnecessary repetitions. By the $c_\epsilon$-equation in \eqref{differ-system} and  $\tilde c_\epsilon$-equation in \eqref{linear-system}, one gets that
	\be\label{c-rep}
	\begin{split}
		c_\epsilon(t)&=e^{t(\Delta-1)}c_I+\int_{0}^{t}e^{(t-s)(\Delta-1)}\rho_\epsilon(s)ds	\\
		&=e^{t(\Delta-1)}c_I+\int_{0}^{t}e^{(t-s)(\Delta-1)}(\rho_\epsilon(s)-\tilde{\rho}_\epsilon(s)+\tilde{\rho}_\epsilon(s))ds	\\
		&=\tilde{c}_\epsilon(t)+\int_{0}^{t}e^{(t-s)(\Delta-1)}(\rho_\epsilon(s)-\tilde{\rho}_\epsilon(s))ds.
	\end{split}
	\ee
Now,  by the maximum principle and semigroup estimates,  from    \eqref{initial-eps-small},  \eqref{rho-epsilon-q0} and the condition $\|c_I\|_{L^\infty}\leq \eta$, we use the integral estimates in  Lemma \ref{integral-lemma} to infer that
	\be\label{c-epsilon-infty}
	\begin{split} \|c_\epsilon(t)\|_{L^\infty}&\leq\|e^{t(\Delta-1)}c_I\|_{L^\infty}+\int_{0}^{t}\|e^{(t-s)(\Delta-1)}\rho_\epsilon(s)\|_{L^\infty}ds	\\
&\leq \|c_I\|_{L^\infty}e^{-t}+C_2\int_{0}^{t}\left(1+(t-s)^{-\frac{n}{2q_0}}\right)e^{-(\lambda_1+1)(t-s)}\|\rho_\epsilon(s)\|_{L^{q_0}}ds\\
		&\leq \|c_I\|_{L^\infty}e^{-t}  +C_2\eta\int_{0}^{t}\left(1+(t-s)^{-\frac{n}{2q_0}}\right)e^{-(\lambda_1+1)(t-s)} e^{-\mu s}\left(1+s^{-1+\frac{n}{2q_0}}\right)ds	\\
		&\leq C_3\eta e^{-\min\{\mu,1\}t}.
	\end{split}
	\ee
Then, we notice from the inequality $|e^s-1|\leq |s|e^{|s|}$ for all $s\in\mathbb{R}$  that
\be\label{int-c-eps}
\left|e^{-\int_{0}^{t}c_\epsilon(s)ds}-1\right|\leq \left|\int_{0}^{t}c_\epsilon(s)ds\right|e^{\left|\int_{0}^{t}c_\epsilon(s)ds\right|}\leq \frac{C_3\eta}{ \min\{\mu,1\}}e^{\frac{C_3\eta}{ \min\{\mu,1\}}}.
\ee
Next, taking the gradient $\nabla$ on both sides of \eqref{c-rep} and applying the smoothing semigroup estimates, we infer from the definition of $T_\epsilon$ in \eqref{T-def} and the integral estimate in Lemma \ref{integral-lemma},  for $ n<q_0\leq p_0<\frac{nq_0}{q_0-n}$, that
	\be\label{grad c-differ}
	\begin{split}
		&\|\nabla c_\epsilon(t)-\nabla\tilde{c}_\epsilon(t)\|_{L^{p_0}}	\\
		&\leq\int_{0}^{t}\|\nabla e^{(t-s)(\Delta-1)}(\rho_\epsilon(s)-\tilde{\rho}_\epsilon(s))\|_{L^{p_0}}ds	\\
		&\leq C_4\int_{0}^{t}\left(1+(t-s)^{-\frac{1}{2}-\frac{n}{2}(\frac{1}{q_0}-\frac{1}{p_0})}\right)e^{-(\lambda_1+1)(t-s)}\|\rho_\epsilon(s)-\tilde{\rho}_\epsilon(s)\|_{L^{q_0}}ds	\\
		&\leq  C_5\eta\int_{0}^{t}\left(1+(t-s)^{-\frac{1}{2}-\frac{n}{2}(\frac{1}{q_0}-\frac{1}{p_0})}\right)e^{-(\lambda_1+1)(t-s)} e^{-\mu s}\left(1+s^{-1+\frac{n}{2q_0}}\right)ds	\\
		&\leq C_6\eta\left(1+t^{-\frac{1}{2}+\frac{n}{2p_0}}\right)e^{-\mu t},
	\end{split}
	\nonumber
	\ee
	which along with \eqref{grad c-lp} and the smallness condition \eqref{initial-eps-small} entails that
	\be\label{grad c-epsilon-lp0}
	\begin{split}
		\|\nabla c_\epsilon(t)\|_{L^{p_0}}&\leq\|\nabla c_\epsilon-\nabla\tilde{c}_\epsilon(t)\|_{L^{p_0}}+\|\nabla\tilde{c}_\epsilon(t)\|_{L^{p_0}}	\\
		&\leq C_7\eta e^{-\mu t}\left(1+t^{-\frac{1}{2}+\frac{n}{2p_0}}\right).
	\end{split}
	\ee
On the other hand, extracting  the $w_\epsilon$-expression  in \eqref{tilde w3}, we find that
	\be\label{w-epsilon}
		w_\epsilon=w_I  e^{-\int_0^t(c_\epsilon(s)+\bar u_0ds},
	\ee
	which along  with $c_\epsilon+\bar{u}_0\geq v_*$  as in \eqref{c-lower-bound}, \eqref{grad c-epsilon-lp0} and $\|\nabla w_I\|_{L^{p_0}}<\eta$ implies  that
	\be\label{grad w-epsilon-lp0}
	\begin{split}
		\|\nabla w_\epsilon(t)\|_{L^{p_0}}&\leq\|\nabla w_I\|_{L^{p_0}}e^{-v_*t}+\|w_I\|_{L^\infty}e^{-v_*t}\int_{0}^{t}\|\nabla c_\epsilon(s)\|_{L^{p_0}}	\\
		&\leq\eta e^{-v_*t}+C_7\eta e^{-v_*t}\|w_I\|_{L^\infty}\int_{0}^{\infty}\left(1+s^{-\frac{1}{2}+\frac{n}{2p_0}}\right)e^{-\mu s}ds	\\
		&\leq C_8 \eta e^{-v_*t}.
	\end{split}
	\ee
To proceed further, using  the transformations in \eqref{tranfs-J1J2-etc}, we rewrite \eqref{difference-system} as follows:
\be\label{difference-system+}
	\begin{cases}
		\hat \rho_{\epsilon t}-\Delta \hat\rho_\epsilon-\alpha\Delta \hat c_\epsilon=-\beta\Delta
\hat w_\epsilon +\chi\nabla\cdot J_{1\epsilon}-\xi\nabla\cdot J_{2\epsilon} &\text{in } \Omega\times(0,\infty), \\[0.2cm]
		\hat c_{\epsilon t}-\Delta\hat c_\epsilon+\hat c_\epsilon-\hat\rho_\epsilon=0  &\text{in } \Omega\times(0,\infty),\\[0.2cm]
		\hat w_\epsilon=w_Ie^{-\bar{u}_0t}\left(e^{-\int_0^tc_\epsilon(s)ds}-1\right)	 &\text{in } \Omega\times(0,\infty),\\[0.2cm]
		\frac{\partial \hat\rho_\epsilon}{\partial\nu}=\frac{\partial\hat c_\epsilon}{\partial\nu}=\frac{\partial\hat w_\epsilon}{\partial\nu}=0  &\text{on } \partial\Omega\times(0,\infty),\\[0.2cm]
		(\hat\rho_\epsilon, \  \hat c_\epsilon, \ \hat w_\epsilon)=(0,\  0, \  0) &\text{on } \Omega\times\{0\}.
	\end{cases}
	\ee
Taking $\nabla$ to the third equation in \eqref{difference-system+}, noting that $c_\epsilon+\bar{u}_0\geq v_*$  in \eqref{c-lower-bound} and $\|w_I\|_{L^\infty}+\|\nabla w\|_{L^{p_0}}\leq \eta$ in \eqref{initial-eps-small},  we infer from   \eqref{int-c-eps} and \eqref{grad c-epsilon-lp0} that
\be\label{gradw-hat-est}
\begin{split}
&\|\nabla (w_\epsilon-\tilde w_\epsilon)\|_{L^{p_0}}=\|\nabla \hat w_\epsilon\|_{L^{p_0}}\\
&\leq \|\nabla w_I\|_{L^{p_0}}e^{-\bar u_0 t}|e^{-\int_0^tc_\epsilon(s)ds}-1| +\|w_I\|_{L^\infty}e^{-v_*t}\int_0^t\|\nabla c_\epsilon(s)\|_{L^{p_0}}ds\\
&\leq C_9\eta\|\nabla w_I\|_{L^{p_0}} e^{-\bar u_0 t}+\eta\|w_I\|_{L^\infty} e^{-v_*t}\int_{0}^t\left(1+s^{-\frac{1}{2}+\frac{n}{2p_0}}\right)
e^{-\mu s}ds	\\
		&\leq C_{10}\eta^2e^{-\min\{\bar u_0,  v_*\}t}.
\end{split}
\ee
Next, applying  the semigroup $(e^{t\mathcal{A}})_{t\geq 0}$ to the $(\tilde \rho_\epsilon, \tilde c_\epsilon)$ equations in \eqref{linear-system}, and combining the semigroup representation of $(\rho_\epsilon, c_\epsilon)$ in \eqref{rho-epsilon-rep+}, we see that
	\be\label{rho-epsilon-rep}
\begin{split}
	&\begin{pmatrix}
		\rho_\epsilon(t)-\tilde{\rho}_\epsilon(t) \\ c_\epsilon(t)-\tilde{c}_\epsilon(t)
	\end{pmatrix}\\[0.2cm]
&=\int_{0}^{t} e^{(t-s)\mathcal{A}}
	\begin{pmatrix}
		-\beta\Delta(w_\epsilon-\tilde{w}_\epsilon)+\chi\nabla\cdot J_{1\epsilon}-\xi\nabla\cdot J_{2\epsilon} \\ 0
	\end{pmatrix}(s)ds\\[0.2cm]
&=:\int_0^t\begin{pmatrix}
	 e_{11}^{(t-s)\mathcal{A}}(-\beta\Delta(w_\epsilon-\tilde{w}_\epsilon),0 ) \\ e_{21}^{(t-s)\mathcal{A}}(-\beta\Delta(w_\epsilon-\tilde{w}_\epsilon),0)
	\end{pmatrix}ds \\[0.2cm]
&\ \ +\int_0^t\begin{pmatrix}
		 e_{12}^{(t-s)\mathcal{A}}(\chi\nabla\cdot J_{1\epsilon}, 0) \\  e_{22}^{(t-s)\mathcal{A}}(\chi\nabla\cdot J_{1\epsilon}, 0)
	\end{pmatrix} ds+\int_0^t\begin{pmatrix}
		 e_{13}^{(t-s)\mathcal{A}}(-\xi\nabla\cdot J_{2\epsilon}, 0)\\  e_{23}^{(t-s)\mathcal{A}}(-\xi\nabla\cdot J_{2\epsilon}, 0)
	\end{pmatrix}ds.
\end{split}
	\ee
In the sequel, we aim to control the terms on the right-hand side of \eqref{rho-epsilon-rep}. For the first term, using the semigroup  estimate for $(e^{t\mathcal{A}})_{t\geq 0}$  in Lemma \ref{tilde rho-lp-lemma},  H\"{o}lder's inequality  and Lemma \ref{integral-lemma}, we use \eqref{gradw-hat-est} and the fact $0<\beta\leq \xi S(\bar u_0)$ to bound
\be\label{sem-est1}
\begin{split}
&\int_{0}^{t} \|e_{11}^{(t-s)\mathcal{A}} \left(-\beta\Delta(w_\epsilon-\tilde{w}_\epsilon), \ 0\right)\|_{L^{q_0}}ds\\
&\leq C_{11}\beta\int_{0}^{t}\left(1+(t-s)^{-\frac{1}{2}}\right)e^{-\lambda_1(t-s)}
\|\nabla(\tilde{w}_\epsilon-w_\epsilon)(s)\|_{L^{q_0}}ds\\
&\leq C_{11}\xi S(\bar u_0)|\Omega|^{\frac{1}{q_0}-\frac{1}{p_0}}\int_{0}^{t}\left(1+(t-s)^{-\frac{1}{2}}\right)
e^{-\lambda_1(t-s)}\|\nabla(\tilde{w}_\epsilon-w_\epsilon)(s)\|_{L^{p_0}}ds\\
&\leq C_{12}\eta^2\int_{0}^{t}\left(1+(t-s)^{-\frac{1}{2}}\right)
e^{-\lambda_1(t-s)}e^{-\min\{\bar u_0,  v_*\}s}ds\\
&\leq C_{13}\eta^2(1+\sqrt{t})e^{-\min\{\lambda_1, \bar u_0,  v_*\}t}.
\end{split}
\ee
Similarly, by \eqref{rho-epsilon-q0},  \eqref{grad c-epsilon-lp0},  \eqref{tranfs-J1J2-etc} and the facts $0\leq 1-\epsilon(\rho_\epsilon+\bar u_0)\leq 1$ and $0<\epsilon \bar u_0<1$, for $\frac{1}{r_0}=\frac{1}{p_0}+\frac{1}{q_0}$, we derive that
\begin{equation*}
\begin{split}
\|J_{1\epsilon}\|_{L^{r_0}}&=\left\|[1-\epsilon(\rho_\epsilon+\bar{u}_0)-\epsilon\bar{u}_0]\rho_\epsilon\nabla c_\epsilon\right\|_{L^{r_0}}\\
&\leq 2\|\rho_\epsilon\nabla c_\epsilon\|_{L^{r_0}}\\
&\leq  2\|\rho_\epsilon\|_{L^{q_0}}\|\nabla c_\epsilon\|_{L^{p_0}}\\
&\leq C_{14}\eta^2
e^{-2\mu t} (1+t^{-1+\frac{n}{2q_0}})(1+t^{-\frac{1}{2}+\frac{n}{2p_0}})\\
&\leq C_{15}\eta^2
e^{-2\mu t} (1+t^{-\frac{3}{2}+\frac{n}{2r_0}}),
\end{split}
\end{equation*}
so that we obtain that
\be\label{sem-est2}
\begin{split}
&\int_{0}^{t} \|e_{12}^{(t-s)\mathcal{A}}\left(\chi \nabla\cdot J_{1\epsilon}, 0\right)\|_{L^{q_0}}ds\\
&\leq C_{16}\int_{0}^{t}\left(1+(t-s)^{-\frac{1}{2}-\frac{n}{2}(\frac{1}{r_0}-\frac{1}{q_0})}\right)e^{-\lambda_1(t-s)}
\|J_{1\epsilon}(s)\|_{L^{r_0}}ds\\
&\leq C_{17}\eta^2\int_{0}^{t}\left(1+(t-s)^{-\frac{1}{2}-\frac{n}{2p_0}}\right)e^{-\lambda_1(t-s)}
e^{-2\mu s} (1+s^{-\frac{3}{2}+\frac{n}{2r_0}})ds\\
&\leq C_{18}\eta^2\left(1+t^{-1+\frac{n}{2q_0}}\right) e^{-\min\{\lambda_1, 2\mu\}t}.
\end{split}
\ee
Recalling  that $S(0)=0$, $0\leq 1-\epsilon\bar{u}_0\leq 1$, $0\leq \epsilon (\rho_\epsilon+\bar u_0)\leq 1$ and that the derivative $S^\prime$ is bounded on $\mathbb{R}^+$,  from  \eqref{rho-epsilon-q0},   \eqref{grad w-epsilon-lp0} and \eqref{tranfs-J1J2-etc}, we use the intermediate  mean value theorem to bound, for some $\sigma_\epsilon$ and $\delta_\epsilon$,  that
\begin{equation*}
\begin{split}
\|J_{2\epsilon}\|_{L^{r_0}}&=\left\|\left[(S(\rho_\epsilon+\bar{u}_0)-S(\bar{u}_0))(1-\epsilon\bar{u}_0)-\epsilon \rho_\epsilon S(\rho_\epsilon+\bar{u}_0)\right]\nabla w_\epsilon\right\|_{L^{r_0}}\\
&=\left\|\left[S^\prime(\sigma_\epsilon)(1-\epsilon\bar{u}_0)-\epsilon  S^\prime(\delta_\epsilon)(\rho_\epsilon+\bar{u}_0)\right]\rho_\epsilon\nabla w_\epsilon\right\|_{L^{r_0}}\\
&\leq C_{19}\left\| \rho_\epsilon \nabla w_\epsilon\right\|_{L^{r_0}}\\
&\leq C_{19} \|\rho_\epsilon\|_{L^{q_0}}\|\nabla w_\epsilon\|_{L^{p_0}}\\
&\leq C_{20}\eta^2
e^{-(\mu+v_*) t} (1+t^{-1+\frac{n}{2q_0}}),
\end{split}
\end{equation*}
and therefore, we have
\be\label{sem-est3}
\begin{split}
&\int_{0}^{t} \|e_{13}^{(t-s)\mathcal{A}} \left(-\xi\nabla\cdot J_2, \ 0\right)\|_{L^{q_0}}ds\\
&\leq C_{21}\int_{0}^{t}\left(1+(t-s)^{-\frac{1}{2}-\frac{n}{2}(\frac{1}{r_0}-\frac{1}{q_0})}\right)e^{-\lambda_1(t-s)}
\|J_2(s)\|_{L^{r_0}}ds\\
&\leq C_{22}\eta^2\int_{0}^{t}\left(1+(t-s)^{-\frac{1}{2}-\frac{n}{2p_0}}\right)e^{-\lambda_1(t-s)}
e^{-(\mu+v_*) s} (1+s^{-1+\frac{n}{2q_0}})ds\\
&\leq C_{23}\eta^2\left(1+t^{-\frac{1}{2}+\frac{n}{2}(\frac{1}{q_0}-\frac{1}{p_0})}\right) e^{-\min\{\lambda_1, \mu+v_*\}t}.
\end{split}
\ee
Finally, substituting \eqref{sem-est1}, \eqref{sem-est2} and \eqref{sem-est3} into \eqref{rho-epsilon-rep} and using the linearity of $e^{t\mathcal{A}}$, we  end up with
	\be\label{rho-rho-differ}
	\begin{split}
		&\|\rho_\epsilon(t)-\tilde{\rho}_\epsilon(t)\|_{L^{q_0}}\\
&\leq C_{13}\beta\eta^2(1+\sqrt{t})e^{-\min\{\lambda_1, \bar u_0,  v_*\}t}+C_{18}\eta^2\left(1+t^{-1+\frac{n}{2q_0}}\right) e^{-\min\{\lambda_1, 2\mu\}t}\\
&\ \ +C_{23}\eta^2\left(1+t^{-\frac{1}{2}+\frac{n}{2}(\frac{1}{q_0}-\frac{1}{p_0})}\right) e^{-\min\{\lambda_1, \mu+v_*\}t}\\
&\leq C_{24}\eta^2\left(1+t^{-1+\frac{n}{2q_0}}\right) e^{-\mu t},
	\end{split}
	\ee
where the facts  $p_0>n, q_0>n$ and $0<\mu<\min\{\lambda_1, \bar u_0, v_*\}$ were used.

	Upon fixing $0<\eta_0<\frac{1}{C_{24}}$ in \eqref{rho-rho-differ}, we infer from the definition of $T_\epsilon$ in \eqref{T-def} that $T_\epsilon=\infty$, and thus, we finally achieve  the desired estimate \eqref{rho-rhoepsilon}.	
\end{proof}
	
 \begin{remark} One may also view the nonlinear system \eqref{differ-system} as
 \begin{equation*}
	\begin{pmatrix}
		\rho_\epsilon(t) \\ c_\epsilon(t) \\ w_\epsilon
	\end{pmatrix}=e^{t\mathcal{B}}\begin{pmatrix}
		\rho_I \\ c_I\\ w_I
	\end{pmatrix}+\int_{0}^{t} e^{(t-s)\mathcal{B}}
	\begin{pmatrix}
	\chi\nabla\cdot J_1-\xi\nabla\cdot J_2 \\ 0 \\
    -c_\epsilon w_\epsilon
	\end{pmatrix}(s)ds, 	\end{equation*}
where the realization $\mathcal{B}:D(\mathcal{B})\subset (L^p(\Omega))^3\mapsto (L^p(\Omega))^3$ is defined as follows:
\begin{equation*}\label{Bmatrix-def}
		\mathcal{B}=\begin{pmatrix}
		\Delta & \alpha\Delta	& -\beta \Delta \\
		1 & \Delta-1 & 0 \\
        0  & 0        & -\bar u_0
		\end{pmatrix} \end{equation*}
with domain $D(\mathcal{B})=W^{2,p}_{N, 0}(\Omega)\times W^{2,p}_N(\Omega)\times W^{2,p}_N(\Omega)$.  Then one can see that   $\mathcal{B}(\rho, c, w)=(0,0,0)$ in $D(\mathcal{B})$ implies  $(\rho, c, w)= (0,0, 0)$, and when $\epsilon\leq \frac{1}{2\bar u_0}$, one can readily see that $-\mathcal{B}$ has a complete sequence of eigenvalues $\{\tilde\lambda_k\}_{k=1}^\infty$ with $0<\text{Re}(\tilde\lambda_k)\leq \text{Re}(\tilde\lambda_{k+1})$ for all $k\geq 1$ and that the infimum $\hat\lambda_1$ of the real part of the spectrum $\sigma(-\mathcal{B})$ is positive. Upon noticing that
$$
	\begin{pmatrix}
		\rho_\epsilon(t)-\tilde{\rho}_\epsilon(t) \\ c_\epsilon(t)-\tilde{c}_\epsilon(t)\\ w_\epsilon(t)-\tilde w_\epsilon(t)
	\end{pmatrix}=\int_{0}^{t} e^{(t-s)\mathcal{B}}
	\begin{pmatrix}
	\chi\nabla\cdot J_{1\epsilon}-\xi\nabla\cdot J_{2\epsilon}\\ 0 \\
    -c_\epsilon w_\epsilon
	\end{pmatrix}(s)ds,
$$
applying the usual semigroup estimates for $e^{t\mathcal{B}}$,  one may derive Lemma \ref{proposition} as well.
 \end{remark}
 	Based on Lemmas \ref{exponential-decay-lemma} and \ref{proposition}, for  initial datum  $(u_0,v_0,w_0)$ fulfilling \eqref{initial-data},  we set  $\epsilon_0=\min\{\frac{1}{2\bar u_0}, \ \ \|u_0\|_{L^\infty}^{-1}\}$ and then fix $\eta_0>0$  as provided by Lemma \ref{proposition}, and then use the relations  $(u_\epsilon, v_\epsilon, w_\epsilon)=(\rho_\epsilon+\bar u_0, c_\epsilon+\bar u_0,w_\epsilon)$  to demonstrate that the smallness condition \eqref{initial-eps-small} can be fulfilled at certain time for $n=3$.

\begin{lemma}\label{rho-l1.5-lemma}
	Assume the basic  Lemma \ref{regular-system-local-lemma} holds with $n=3$,  $\inf_\Omega v_0>0$,  $S$ satisfies \eqref{S-con} with $\theta\leq\frac{1}{2}$ and the derivative $S^\prime$ being bounded on $\mathbb{R}^+$. Then,  for  $\epsilon_0, \eta_0$ be as given by Lemma \ref{proposition} and for  $\eta\in (0, \eta_0)$, there exists $t_\eta>0$ such that the  solution of the approximating system  \eqref{regularized-equation} fulfills, for all $\epsilon\in(0,\epsilon_0)$,
	\be\label{rho-l1.5}
\begin{split}
		&\|u_\epsilon(t_\eta)-\bar{u}_0\|_{L^\frac{3}{2}}+\|\nabla v_\epsilon(t_\eta)\|_{L^6}+\|v_\epsilon(t_\eta)-\bar{u}_0\|_{L^\infty}\\[0.2cm]
&+\|w_\epsilon(t_\eta)\|_{L^\infty}+\|\nabla w_\epsilon(t_\eta)\|_{L^6}\leq\eta.
\end{split}	\ee
\end{lemma}
\begin{proof}
	By the  Sobolev embedding $W^{1,1}(\Omega)\hookrightarrow L^\frac{3}{2}(\Omega)$ with $\Omega\subset\mathbb{R}^3$,   Poincar\'{e}'s inequality and H\"{o}lder's inequality enable us to see that
	\be\label{u-l1.5}
		\|u_\epsilon-\bar{u}_0\|_{L^\frac{3}{2}}^2\leq C_1\|\nabla u_\epsilon\|_{L^1}^2\leq C_1\int_\Omega u_\epsilon\int_\Omega\frac{|\nabla u_\epsilon|^2}{u_\epsilon}=C_1\|u_0\|_{L^1}\int_\Omega\frac{|\nabla u_\epsilon|^2}{u_\epsilon}.
	\ee
	Again by the Sobolev embedding $W^{2,2}(\Omega)\hookrightarrow W^{1,6}(\Omega)$ in 3D and the $W^{2,2}$-elliptic estimate, we obtain that
	\be\label{grad v-6}
		\|\nabla v_\epsilon\|_{L^6}^2=\|\nabla (v_\epsilon-\bar v_\epsilon)\|_{L^6}^2\leq C_2\|v_\epsilon-\bar v_\epsilon\|_{W^{2,2}}^2\leq C_3\|\Delta v_\epsilon\|_{L^2}^2.
	\ee
Hence, we employ the  3D Gagliardo-Nirenberg inequality as in \eqref{GN-intro}, Poincar\'{e}'s inequality and \eqref{grad v-6} to deduce that
	$$
		\|v_\epsilon-\bar{v}_\epsilon\|_{L^\infty}^2\leq C_4\|\nabla v_\epsilon\|_{L^6}\|v_\epsilon-\bar{v}_\epsilon\|_{L^6}\leq C_5\|\nabla v_\epsilon\|_{L^6}^2\leq C_6\|\Delta v_\epsilon\|_{L^2}^2,
	$$
	which together with the $v_\epsilon$-equation in \eqref{regularized-equation} implies that
	\be\label{v-infty}
	\begin{split}
		\|v_\epsilon-\bar{u}_0\|_{L^\infty}^2&\leq2\|v_\epsilon-\bar{v}_\epsilon\|_{L^\infty}^2+2\|\bar{v}_\epsilon-\bar{u}_0\|_{L^\infty}^2	\\
		&\leq 2C_6\|\Delta v_\epsilon\|_{L^2}^2+2|\bar{v}_0-\bar{u}_0|e^{-2t}.
	\end{split}
	\ee
Therefore, we combine  \eqref{u-l1.5}, \eqref{grad v-6} and \eqref{v-infty} together, and then, we apply \eqref{delta v-grad v-estimate-1} of  Lemma \ref{exponential-decay-lemma} thanks to $\theta\leq \frac{1}{2}$ to infer that
	\be\label{bridge1}
	\begin{split}
		&\int_{0}^{\infty}\left(\|u_\epsilon-\bar{u}_0\|_{L^\frac{3}{2}}^2+\|\nabla v_\epsilon\|_{L^6}^2+\|v_\epsilon-\bar{u}_0\|_{L^\infty}^2\right)	\\
		&\leq C_7\int_{0}^{\infty}\int_\Omega\left(\frac{|\nabla u_\epsilon|^2}{u_\epsilon}+|\Delta v_\epsilon|^2\right)+C_7\int_{0}^{\infty}e^{-2t}dt	\leq C_8.
	\end{split}
	\nonumber
	\ee
 In view of the fact in Lemma \ref{limit-lemma} with any $a>1$, for any $\eta\in(0,\eta_0)$, one can find some  $\tilde{t}_\eta>a$ independent of $\epsilon$, and a   time $\tilde{t}_\epsilon(\eta)\in(a,\tilde{t}_\eta)$ such that
\be\label{small-need1}
		\|u_\epsilon(\tilde{t}_\epsilon(\eta))-\bar{u}_0\|_{L^\frac{3}{2}}+\|\nabla v_\epsilon(\tilde{t}_\epsilon(\eta))\|_{L^6}+\|v_\epsilon(\tilde{t}_\epsilon(\eta))-\bar{u}_0\|_{L^\infty}
\leq\frac{\eta}{2}.
	\ee
Then we use  the fact $v_\epsilon\geq v_*$ in \eqref{v-low-bdd-eps} of Lemma \ref{gradw-l2-decay} and  the counterparts for $w_\epsilon$ in  \eqref{w-expression} and \eqref{w-lp} of Lemma \ref{w-lemma}  to deduce that
\begin{equation*}
\begin{split}
		&\|w_\epsilon(\tilde{t}_\epsilon(\eta))\|_{L^\infty}+\|\nabla w_\epsilon(\tilde{t}_\epsilon(\eta))\|_{L^6}-\left(\|w_0\|_{L^\infty}+\|\nabla w_0\|_{L^6}\right)e^{-v_*\tilde{t}_\epsilon(\eta)}\\
&\leq\|w_0\|_{L^\infty} e^{-v_*\tilde{t}_\epsilon(\eta)}\int_0^{\tilde{t}_\epsilon(\eta)}\|\nabla v_\epsilon(s)\|_{L^6}ds\\
&\leq \|w_0\|_{L^\infty} e^{-v_*\tilde{t}_\epsilon(\eta)}\sqrt{\tilde{t}_\epsilon(\eta)}\left(\int_0^{\tilde{t}_\epsilon(\eta)}
\|\nabla v_\epsilon(s)\|_{L^6}^2ds\right)^\frac{1}{2}\\
&\leq \|w_0\|_{L^\infty} e^{-\frac{v_*}{2}\tilde{t}_\epsilon(\eta)}\sup_{t>0}\left(t^\frac{1}{2}e^{-\frac{v_*}{2}t}\right)\left(\int_0^\infty\|\nabla v_\epsilon(s)\|_{L^6}^2ds\right)^\frac{1}{2}\\
&\leq \left(\frac{C_8}{ev_*}\right)^\frac{1}{2}\|w_0\|_{L^\infty} e^{-\frac{v_*}{2}\tilde{t}_\epsilon(\eta)},
\end{split}
	\end{equation*}
which along with the fact $\tilde{t}_\epsilon(\eta)>a$ implies
\be\label{w-need-es}
\begin{split}
		&\|w_\epsilon(\tilde{t}_\epsilon(\eta))\|_{L^\infty}+\|\nabla w_\epsilon(\tilde{t}_\epsilon(\eta))\|_{L^6}\\
&\leq \left(\|w_0\|_{L^\infty}+\|\nabla w_0\|_{L^6}\right)e^{-v_*a}+\left(\frac{C_8}{ev_*}\right)^\frac{1}{2}\|w_0\|_{L^\infty} e^{-\frac{v_*}{2}a},
\end{split}
\ee
which, upon choosing $a>1$ suitably large, along with \eqref{small-need1} quickly entails that
\be\label{small-need2}
\begin{split}
&\|u_\epsilon(\tilde{t}_\epsilon(\eta))-\bar{u}_0\|_{L^\frac{3}{2}}+\|\nabla v_\epsilon(\tilde{t}_\epsilon(\eta))\|_{L^6}+\|v_\epsilon(\tilde{t}_\epsilon(\eta))
-\bar{u}_0\|_{L^\infty}\\[0.2cm]
&\ \ +\|w_\epsilon(\tilde{t}_\epsilon(\eta))\|_{L^\infty}+\|\nabla w_\epsilon(\tilde{t}_\epsilon(\eta))\|_{L^6}
\leq \eta.
\end{split}
	\ee
 Now, viewing the solution  for  $t\geq\tilde{t}_\epsilon(\eta)$ as a trajectory that initially starts from the datum $(u_\epsilon(\tilde{t}_\epsilon(\eta)),v_\epsilon(\tilde{t}_\epsilon(\eta), w_\epsilon(\tilde{t}_\epsilon(\eta))$.  Recalling that $(u_\epsilon, v_\epsilon, w_\epsilon)=(\rho_\epsilon+\bar u_0, c_\epsilon+\bar u_0,w_\epsilon)$,  by the facts that $u_\epsilon>0$ and $v_\epsilon\geq v_*>0$ on $\Omega \times (0, \infty)$, we apparently see  that $\rho_\epsilon(\tilde{t}_\epsilon(\eta))=u_\epsilon(\tilde{t}_\epsilon(\eta))-\bar u_0>-\bar u_0$ and $c_\epsilon(\tilde{t}_\epsilon(\eta))=v_\epsilon(\tilde{t}_\epsilon(\eta))-\bar u_0>-\bar u_0$; also, it follows from   \eqref{small-need2}  that the combined smallness condition \eqref{initial-eps-small} is satisfied for $n=3$,  $q_0\in(3, 6)$ and $p_0\in [q_0, 6]$. Therefore,  for any $\epsilon\in(0, \epsilon_0)$, $\eta\in (0, \eta_0)$, $q_0\in (3, 6)$ and $\mu\in (0,\min\{1, \lambda_1, \bar u_0, v_*\})$, we utilize  Lemma \ref{proposition} to derive that
	$$
		\|\rho_\epsilon(t)-\tilde{\rho}_\epsilon(t)\|_{L^{q_0}}\leq\eta e^{-\mu(t-\tilde{t}_\epsilon(\eta))}\left(1+(t-\tilde{t}_\epsilon(\eta))^{-1+\frac{3}{2q_0}}\right),\quad\forall t>\tilde{t}_\epsilon(\eta).
	$$
	This along with  application of  \eqref{rho-lp} in  Lemma \ref{rho-grad c-lp-lemma}   shows that
	\be\label{u-u0-q0}
	\begin{split}
		\|u_\epsilon(t)-\bar{u}_0\|_{L^{q_0}}&\leq\|\rho_\epsilon(t)-\tilde{\rho}_\epsilon(t)\|_{L^{q_0}}+\|\tilde{\rho}_\epsilon(t)\|_{L^{q_0}}	\\
		&\leq\eta e^{-\mu(t-\tilde{t}_\epsilon(\eta))}\left(1+(t-\tilde{t}_\epsilon(\eta))^{-1+\frac{3}{2q_0}}\right)	\\
		&\ \ \ +C_9e^{-\mu (t-\tilde{t}_\epsilon(\eta))}\left(1+(t-\tilde{t}_\epsilon(\eta))^{-1+\frac{3}{2q_0}}\right)	\\
		&\ \ \ \times\left(\|\rho_\epsilon(\tilde{t}_\epsilon(\eta))\|_{L^\frac{3}{2}}+\|\nabla c_\epsilon(\tilde{t}_\epsilon(\eta))\|_{L^3}+\|\nabla w_\epsilon(\tilde{t}_\epsilon(\eta))\|_{L^3}\right)	\\
		&\leq C_{10}\eta e^{-\mu(t-\tilde{t}_\epsilon(\eta))}\left(1+(t-\tilde{t}_\epsilon(\eta))^{-1+\frac{3}{2q_0}}\right),\ \ \forall t>\tilde{t}_\epsilon(\eta).
	\end{split}
	\nonumber
	\ee
On the other hand,  the  estimate for $c_\epsilon=v_\epsilon-\bar u_0$ in  \eqref{c-epsilon-infty} entails that
	$$
	  \|v_\epsilon(t)-\bar u_0\|_{L^\infty}\leq  C_{11}\eta e^{-\mu(t-\tilde{t}_\epsilon(\eta))},\ \ \forall t>\tilde{t}_\epsilon(\eta).
$$
Furthermore, for $p_0\in [q_0, 6]$, we infer from \eqref{grad c-lp} and \eqref{grad c-epsilon-lp0} that
	\be\label{grad v-p0}
	\begin{split}
		\|\nabla v_\epsilon(t)\|_{L^{p_0}}&\leq\|\nabla c_\epsilon(t)-\nabla\tilde{c}_\epsilon(t)\|_{L^{p_0}}+\|\nabla\tilde{c}_\epsilon(t)\|_{L^{p_0}}	\\
		&\leq C_{12}\eta e^{-\mu(t-\tilde{t}_\epsilon(\eta))}\left(1+(t-\tilde{t}_\epsilon(\eta))^{-\frac{1}{2}+\frac{3}{2p_0}}\right),\quad\forall t>\tilde{t}_\epsilon(\eta).
	\end{split}
	\nonumber
	\ee
For $t>\tilde{t}_\epsilon(\eta)$, using  the expressions
$$
w_\epsilon(t)=w_\epsilon(\tilde{t}_\epsilon(\eta))e^{-\int_{\tilde{t}_\epsilon(\eta)}^tv_\epsilon(s)ds}
$$
and
$$
\nabla   w_\epsilon(t)=\nabla w_\epsilon(\tilde{t}_\epsilon(\eta))e^{-\int_{\tilde{t}_\epsilon(\eta)}^tv_\epsilon(s)ds}
-w_\epsilon(\tilde{t}_\epsilon(\eta))e^{-\int_{\tilde{t}_\epsilon(\eta)}^t
v_\epsilon(s)ds}\int_{\tilde{t}_\epsilon(\eta)}^t\nabla v_\epsilon(s)ds,
$$
by \eqref{small-need2} and the fact  $v_\epsilon>v_*>0$, similar to \eqref{w-need-es}, we obtain that
\begin{equation*}
\begin{split}
		&\|w_\epsilon(t)\|_{L^\infty}+\|\nabla w_\epsilon(t)\|_{L^6}\\
&\leq \left(\|w(\tilde{t}_\epsilon(\eta))\|_{L^\infty}+\|\nabla w(\tilde{t}_\epsilon(\eta))\|_{L^6}\right)e^{-v_*(t-\tilde{t}_\epsilon(\eta))}\\
&\ \
+C_{13}\|w(\tilde{t}_\epsilon(\eta))\|_{L^\infty} e^{-\mu(t-\tilde{t}_\epsilon(\eta))}\\
&\leq C_{14}\eta e^{-\mu(t-\tilde{t}_\epsilon(\eta))}, \ \ \forall t>\tilde{t}_\epsilon(\eta).
\end{split}
\end{equation*}
	Recalling that $\tilde{t}_\epsilon(\eta)<\tilde{t}_\eta$ with $\tilde{t}_\eta$ independent of $\epsilon$, collecting the relative estimates above, we conclude   that
\begin{equation*} \begin{split}
		&\|u_\epsilon(t)-\bar{u}_0\|_{L^{q_0}}+\|\nabla v_\epsilon(t)\|_{L^{p_0}}+\|v_\epsilon(t)-\bar{u}_0\|_{L^\infty}\\[0.2cm]
&+\|w_\epsilon(t)\|_{L^\infty}+\|\nabla w_\epsilon(t)\|_{L^6}\leq C_{15}\eta  e^{-\mu(t-1-\tilde{t}_\eta)}, \ \ \forall t\geq 1+\tilde{t}(\eta).
\end{split}	
\end{equation*}
Finally, our asserted estimate \eqref{rho-l1.5} follows upon    choosing $t_\eta>1+\tilde{t}_\eta$ sufficiently large and noticing $q_0\in (3,6)$ and $p_0\in[q_0,6]$.
\end{proof}

	Before showing the eventual regularity of the limiting functions, we must derive certain higher-order estimates for $(u_\epsilon,v_\epsilon,w_\epsilon)$ that are independent of $\epsilon$, cf. \cite{Jiang22, Lankeit15}.
	
\begin{lemma}\label{u-epsilon-infty-lemma}
	Under Lemma \ref{rho-l1.5-lemma}, there exists $C>0$ independent of $\epsilon$ such that the  solution of the approximating system  \eqref{regularized-equation} fulfills, for all $\epsilon\in(0,\epsilon_0)$, $\eta\in (0, \eta_0)$ and $\mu\in (0,\min\{1, \lambda_1, \bar u_0, v_*\})$, that
	\be\label{uvw-com-con-es}
\begin{split}
		&\|u_\epsilon(t)-\bar{u}_0\|_{L^\infty}+\|\nabla v_\epsilon(t)\|_{L^\infty}+\|v_\epsilon(t)-\bar{u}_0\|_{L^\infty}\\[0.2cm]
&\ \ +\|w_\epsilon(t)\|_{L^\infty}+\|\nabla w_\epsilon(t)\|_{L^6}\leq Ce^{-\mu t}, \ \ \forall t\geq 3+t_\eta.
\end{split}
	\ee
	Moreover, there is $\sigma\in(0,1)$ such that for any $T\geq 3+t_\eta$, there exists $C=C(T)>0$ such that
	\be\label{u-epsilon-regularity}
		\|u_\epsilon\|_{C^{\sigma,\frac{\sigma}{2}}(\bar{\Omega}\times[3+t_\eta, T])}
+\|v_\epsilon\|_{C^{2+\sigma,1+\frac{\sigma}{2}}(\bar{\Omega}\times[3+t_\eta, T])}\leq C(T).
	\ee
\end{lemma}
\begin{proof}
Due to the uniqueness of classical solutions, we can again regard the family of solutions $(u_\epsilon,v_\epsilon,w_\epsilon)$ for $t\geq t_\eta$ as a set of trajectories uniformly starting from $(u_\epsilon(t_\eta),v_\epsilon(t_\eta),w_\epsilon(t_\eta))$. Notice then that \eqref{rho-l1.5} implies \eqref{initial-eps-small} for  $n=3$,  $q_0\in(3, 6)$ and $p_0\in [q_0, 6]$.  Applying and repeating the  arguments in Lemma \ref{proposition}, we infer, for   $\epsilon\in(0,\epsilon_0)$, $\eta\in (0, \eta_0)$, $\mu\in (0,\min\{1, \lambda_1, \bar u_0, v_*\})$ and for $t>t_\eta$,  that
	$$
		\|\rho_\epsilon(t)-\tilde{\rho}_\epsilon(t)\|_{L^{q_0}}\leq\eta e^{-\mu(t-t_\eta)}\left(1+(t-t_\eta)^{-1+\frac{3}{2q_0}}\right),
$$
	and thus
	\be\label{u-q0}
		\|\rho_\epsilon(t)\|_{L^{q_0}}=\|u_\epsilon(t)-\bar{u}_0\|_{L^{q_0}}\leq C_1\eta e^{-\mu(t-t_\eta)}\left(1+(t-t_\eta)^{-1+\frac{3}{2q_0}}\right)
	\ee
	and
	\be\label{grad vw-p0}
		\|\nabla c_\epsilon(t)\|_{L^{p_0}}=\|\nabla v_\epsilon(t)\|_{L^{p_0}}\leq C_2\eta e^{-\mu(t-t_\eta)}\left(1+(t-t_\eta)^{-\frac{1}{2}+\frac{3}{2p_0}}\right)
	\ee
as well as
\be\label{w-gradw-vlinfty}
		\|v_\epsilon(t)-\bar{u}_0\|_{L^\infty}+\|w_\epsilon(t)\|_{L^\infty}+\|\nabla w_\epsilon(t)\|_{L^6}\leq C_3\eta  e^{-\mu(t-t_\eta)}.
\ee
Applying the variation-of-constants formula \eqref{c-rep}  for $c_\epsilon$-equation  and using  \eqref{rho-l1.5}  and \eqref{u-q0}, we obtain that
	\be\label{grad v-epsilon-infty}
	\begin{split}
		&\|\nabla v_\epsilon(t)\|_{L^\infty}=\|\nabla c_\epsilon(t)\|_{L^\infty}	\\
		&\leq\|\nabla e^{(t-t_\eta)(\Delta-1)}c_\epsilon(t_\eta)\|_{L^\infty}+\int_{t_\eta}^{t}\|\nabla e^{(t-s)(\Delta-1)}\rho_\epsilon(s)\|_{L^\infty}ds	\\
		&\leq C_4e^{-(\lambda_1+1)(t-t_\eta)}\left(1+(t-\eta)^{-\frac{1}{2}}\right)\|c_\epsilon(t_\eta)\|_{L^\infty}\\
		&\ \ \ +C_4\int_{t_\eta}^{t}\left(1+(t-s)^{-\frac{1}{2}-\frac{3}{2q_0}}\right)
e^{-\lambda_1(t-s)}\|\rho_\epsilon(s)\|_{L^{q_0}}ds	\\
		&\leq C_5e^{-(\lambda_1+1)(t-t_\eta)}\left(1+(t-t_\eta)^{-\frac{1}{2}}\right)\\
		& +C_5\int_{t_\eta}^{t}\left(1+(t-s)^{-\frac{1}{2}-\frac{3}{2q_0}}\right)e^{-\lambda_1(t-s)}
\left(1+(s-t_\eta)^{-1+\frac{3}{2q_0}}\right)e^{-\mu(s-t_\eta)}ds	\\
		&\leq C_6(1+(t-t_\eta)^{-\frac{1}{2}})e^{-\mu (t-t_\eta)}\leq 2C_6e^{-\mu (t-t_\eta)},\quad\forall t\geq1+t_\eta.
	\end{split}
	\ee
Hence, for $q_0\in(3,6)$ and $s\geq 1+t_\eta$, the estimates \eqref{u-q0} and \eqref{grad v-epsilon-infty} imply
\begin{equation*}
\begin{split}
\|u_\epsilon(s)\nabla c_\epsilon(s)\|_{L^{q_0}}&\leq \|u_\epsilon(s)\|_{L^{q_0}}\|\nabla c_\epsilon(s)\|_{L^\infty}\\
&\leq \left(\|u_\epsilon(s)-\bar u_0\|_{L^{q_0}}+\|\bar u_0\|_{L^{q_0}}\right)\|\nabla c_\epsilon(s)\|_{L^\infty}\leq C_7 e^{-\mu s}.
\end{split}
\end{equation*}
  On the other hand, since $\theta\leq \frac{1}{2}$ implies $0\leq S(u_\epsilon)\leq u_\epsilon$, and then, for $s\geq 1+t_\eta$, the estimates \eqref{u-q0} and \eqref{w-gradw-vlinfty} also imply
\begin{equation*}
\begin{split}
\|(S(u_\epsilon)\nabla w_\epsilon)(s)\|_{L^\frac{6q_0}{6+q_0}}&\leq \|S(u_\epsilon(s))\|_{L^{q_0}}\|\nabla w_\epsilon(s)\|_{L^6}\leq \|u_\epsilon(s)\|_{L^{q_0}}\|\nabla w_\epsilon(s)\|_{L^6}\\
&\leq \left(\|u_\epsilon(s)-\bar u_0\|_{L^{q_0}}+\|\bar u_0\|_{L^{q_0}}\right)\|\nabla w_\epsilon(s)\|_{L^6}\leq C_8 e^{-\mu s}.
\end{split}
\end{equation*}
By  the two estimates above and the facts $0\leq1-\epsilon(\rho_\epsilon+\bar{u}_0)\leq1$, for $q_0\in (3,6)$, we employ the semigroup representation of $\rho_\epsilon=u_\epsilon-\bar u_0$ in \eqref{differ-system} to deduce that
	\be\label{u-epsilon}
	\begin{split}
		&\|u_\epsilon(t)-\bar{u}_0\|_{L^\infty}=\|\rho_\epsilon(t)\|_{L^\infty}\\
		&\leq\|e^{(t-1-t_\eta)\Delta}\rho_\epsilon(1+t_\eta)\|_{L^\infty}+
\chi\int_{1+t_\eta}^{t}\|e^{(t-s)\Delta}\nabla\cdot(u_\epsilon(1-\epsilon u_\epsilon))\nabla c_\epsilon\|_{L^\infty}ds	\\
		&\ \ \ +\xi\int_{1+t_\eta}^{t}\|e^{(t-s)\Delta}\nabla\cdot(S(u_\epsilon)(1-\epsilon u_\epsilon))\nabla w_\epsilon\|_{L^\infty}ds	\\
		&\leq C_9e^{-\lambda_1(t-1-t_\eta)}\left(1+(t-1-t_\eta)^{-\frac{3}{2q_0}}\right)\|\rho_\epsilon(1+t_\eta)\|_{L^{q_0}}	\\
		&\ \ \ +C_9\int_{1+t_\eta}^{t}\left(1+(t-s)^{-\frac{1}{2}-\frac{3}{2q_0}}\right)e^{-\lambda_1(t-s)}
\|u_\epsilon\nabla c_\epsilon\|_{L^{q_0}}ds	\\
		&\ \ \ +C_9\int_{1+t_\eta}^{t}
\left(1+(t-s)^{-\frac{3}{4}-\frac{3}{2q_0}}\right)e^{-\lambda_1(t-s)}
\|S(u_\epsilon)\nabla w_\epsilon\|_{L^\frac{6q_0}{6+q_0}}ds\\
		&\leq C_{10}e^{-\lambda_1(t-1-t_\eta)}+C_{10}\int_{1+t_\eta}^{t}\left(1+(t-s)^{-\frac{1}{2}-\frac{3}{2q_0}}\right)
e^{-\lambda_1(t-s)}e^{-\mu s}ds\\
&\ \ \ +C_{10}\int_{1+t_\eta}^{t}\left(1+(t-s)^{-\frac{3}{4}-\frac{3}{2q_0}}\right)
e^{-\lambda_1(t-s)}e^{-\mu s}ds\\
		&\leq C_{11}(1+(t-1-t_\eta)^{\frac{1}{4}-\frac{3}{2q_0}})e^{-\mu t}\leq C_{12}e^{-\mu t},\quad\forall t\geq2+t_\eta.
	\end{split}
	\nonumber
	\ee
	This in conjunction with \eqref{w-gradw-vlinfty} and \eqref{grad v-epsilon-infty} gives rise to \eqref{uvw-com-con-es}.

Now, rewriting the $u_\epsilon$ as
$$
u_{\epsilon t}=\Delta u_\epsilon+\nabla\cdot (u_\epsilon f_\epsilon(x,t)),
\  \ \   f_\epsilon(x,t):=(1-\epsilon u_\epsilon)\left(\chi \nabla v_\epsilon-\xi\frac{S(u_\epsilon)}{u_\epsilon}\nabla w_\epsilon\right)(x,t),
$$
and noting from \eqref{w-gradw-vlinfty} and \eqref{grad v-epsilon-infty} that
$$
\|f_\epsilon(t)\|_{ L^6}\leq \chi\|\nabla v_\epsilon(t)\|_{L^6}+\xi\|\nabla w_\epsilon(t)\|_{L^6}\leq C_{13}, \ \ \forall t\geq 1+t_\eta,
$$
we therefore apply the standard parabolic H\"{o}lder theory (\cite[Theorem 1.3]{PV93}) to see there exists  $\sigma_1\in(0,1)$ with the property that, for any $T>2+t_\eta$, we can find $C_{14}(T)>0$ such that
\be\label{holder-es}
		\|u_\epsilon\|_{C^{\sigma_1,\frac{\sigma_1}{2}}(\bar{\Omega}\times[2+t_\eta,T])}\leq C_{14}(T).
\ee
Then, applying the standard parabolic Schauder theory (cf. \cite{Friedman69}) to the $v_\epsilon$-equation in \eqref{regularized-equation} shows directly, for any $T>3+t_\eta$,
$$
\|v_\epsilon\|_{C^{2+\sigma_1,1+\frac{\sigma_1}{2}}(\bar{\Omega}\times[3+t_\eta,T])}\leq C_{15}(T),
$$
which along with   \eqref{holder-es} or \eqref{holder-es+} below implies \eqref{u-epsilon-regularity}, completing the proof.
\end{proof}
\begin{remark}In fact, using smoothing effect of the Neumann semigroup,  the H\"{o}lder estimate \eqref{holder-es} can be shown as follows: for some $\sigma_2\in (0,1)$,
\be\label{holder-es+}
		\|u_\epsilon\|_{C^{\sigma_2,\frac{\sigma_2}{2}}(\bar{\Omega}\times[t, t+1])}\leq C_{15}, \ \ \forall t\geq 2+t_\eta.
\ee
As a matter of fact, by the semigroup representation of $u_\epsilon$: with $t\geq t_0:=1+t_\eta$
$$
u_\epsilon(t)-\bar u_0=e^{(t-t_0)\Delta }(u_\epsilon(t_0)-\bar u_0)+\int_{t_0}^te^{(t-s)\Delta }\nabla\cdot(u_\epsilon(s)f_\epsilon(s))ds,
$$
with the uniform $(L^\infty, L^6)$-boundedness of $(u_\epsilon, f)$ on $[t_0, \infty)$, we find (cf. \cite{Fujie-Senba18,HW05}), for any $\gamma\in (0, 1)$ and $q>1$ such that  $2\gamma-\frac{3}{2q}>0$ and $\frac{1}{4}+\frac{3}{2q}-\gamma>0$,
\be\label{holder1-es}
\begin{split}
&\|(-\Delta)^\gamma (u_\epsilon(t)-\bar u_0)\|_{L^q}\\
&\leq \|(-\Delta)^\gamma e^{\Delta} (u_\epsilon(t-1)-\bar u_0)\|_{L^q}\\
&\ \ \ +\int_{t-1}^t\|(-\Delta)^\gamma e^{(t-\tau)\Delta}\nabla\cdot(u(\tau)f_\epsilon(\tau))\|_{L^q}d\tau\\
&\leq C_{16}\|u_\epsilon(t-1)-\bar u_0\|_{L^q}\\
&\ \ \ +C_{16}\int_{t-1}^t (t-\tau)^{-\gamma-\frac{3}{4}+\frac{3}{2q}}e^{-\mu (t-\tau)}\|u(\tau)f_\epsilon(\tau)\|_{L^6}d\tau\\
&\leq C_{17}+C_{17}\int_0^1 \tau^{-\gamma-\frac{3}{4}+\frac{3}{2q}}d\tau=:C_{18},  \  \forall t>2+t_\eta,
\end{split}
\ee
and, for any $1+t_0\leq s\leq t\leq s+1$,
\be\label{holder2-es}
\begin{split}
&\|(-\Delta)^\gamma (u_\epsilon(t)-u_\epsilon(s))\|_{L^q}\\
&\leq \|(-\Delta)^\gamma (e^{(t-t_0)\Delta} -e^{(s-t_0)\Delta})(u_\epsilon(t_0)-\bar u_0)\|_{L^q}\\
&\ \ \ +\int_{t_0}^s\|(-\Delta)^\gamma \left(e^{(t-\tau)\Delta}-e^{(s-\tau)\Delta}\right)\nabla\cdot(u(\tau)f_\epsilon(\tau))\|_{L^q}d\tau\\
&\ \ \ +\int_s^t\|(-\Delta)^\gamma e^{(t-\tau)\Delta}\nabla\cdot(u(\tau)f_\epsilon(\tau))\|_{L^q}d\tau\\
&\leq C_{19}|t-s|\|u_\epsilon(t_0)-\bar u_0\|_{L^q}\\
&\ \ \ +C_{19}|t-s|\int_{t_0}^s\tau^{-\frac{3}{2}-\gamma-\frac{3}{2}(\frac{1}{6}-\frac{1}{q})}\|u(\tau)f_\epsilon(\tau)\|_{L^6}d\tau\\
&\ \ \ +C_{19}\int_s^t (t-\tau)^{-\gamma-\frac{1}{2}-\frac{3}{2}(\frac{1}{6}-\frac{1}{q})}e^{-\mu (t-\tau)}\|u(\tau)f_\epsilon(\tau)\|_{L^6}d\tau\\
&\leq C_{20}|t-s|+C_{20}\int_0^{t-s} \tau^{-\gamma-\frac{3}{4}+\frac{3}{2q}}d\tau\\
&\leq  C_{21}|t-s|^{\frac{1}{4}-\gamma+\frac{3}{2q}},
\end{split}
\ee
where we used the following estimates: for generic functions $g$ and $h$,
\be\label{try-1}
\begin{split}
&\|(-\Delta)^\gamma (e^{(t-t_0)\Delta} -e^{(s-t_0)\Delta})(g-\bar g)\|_{L^q}\\
&=\|(-\Delta)^\gamma\int_{s-t_0}^{t-t_0} \frac{d}{d\tau} e^{\tau \Delta}(g-\bar g)d\tau\|_{L^q}\\
&\leq\int_{s-t_0}^{t-t_0}\| (-\Delta)^{1+\gamma}   e^{\tau \Delta}(g-\bar g)\|_{L^q}d\tau\\
&\leq\int_{s-t_0}^{t-t_0} \tau^{-1-\gamma-\frac{n}{2}(\frac{1}{p}-\frac{1}{q})}  e^{-\mu \tau} \|g-\bar g\|_{L^p}d\tau\\
&\leq |t-s|\|g-\bar g\|_{L^p}, \ \ \forall t\geq s\geq 1+t_0,
\end{split}
\nonumber
\ee
and
\be\label{try-2}
\begin{split}
&\|(-\Delta)^\gamma (e^{(t-t_0)\Delta} -e^{(s-t_0)\Delta})\nabla\cdot h\|_{L^q}\\
&=\|(-\Delta)^\gamma\int_{s-t_0}^{t-t_0} \frac{d}{d\tau} e^{\tau \Delta}\nabla\cdot hd\tau\|_{L^q}\\
&\leq\int_{s-t_0}^{t-t_0}\| (-\Delta)^{1+\gamma}   e^{\tau \Delta}\nabla\cdot h\|_{L^q}d\tau\\
&\leq\int_{s-t_0}^{t-t_0} \tau^{-\frac{3}{2}-\gamma-\frac{n}{2}(\frac{1}{p}-\frac{1}{q})}  e^{-\mu \tau} \|h\|_{L^p}d\tau\\
&\leq |t-s|\|h\|_{L^p}, \ \ \forall t\geq s\geq 1+t_0.
\end{split}
\nonumber
\ee
In the above, $0<\mu<\lambda_1$ with $\lambda_1$ being the first positive eigenvalue to $-\Delta$ under homogeneous Neumann boundary condition. Noticing the continuous embedding $D((-\Delta)^\gamma)\hookrightarrow C^\sigma(\bar{\Omega})$ for any $\sigma\in (0, 2\gamma-\frac{3}{q})$, cf. \cite{HW05, Xiang15}, then the H\"{o}lder estimate \eqref{holder-es+} follows from \eqref{holder1-es} and \eqref{holder2-es}.
\end{remark}
	
\begin{proof}[Proof of Theorem \ref{theorem4}]
 Now, we can extract a subsequence $\{\epsilon\}_{i\geq1}\subset (0, \epsilon_0)$ such that, by taking along the limit $\epsilon=\epsilon_i\searrow 0$, we obtain a weak solution $(u,v,w)$ for the original system \eqref{main equation} in the sense of Definition \ref{weak-definition}, see also Lemma \ref{weak-solution-proposition}. Furthermore, the preceding  higher order estimates in \eqref{u-epsilon-regularity} of Lemma \ref{u-epsilon-infty-lemma} and the Arzel\`{a}-Ascoli theorem imply,  for $\delta\in (0,\sigma)$ and for  $t_0=3+t_\eta$ with $t_\eta$ provided by Lemma \ref{rho-l1.5-lemma}, as $\epsilon=\epsilon_i\searrow 0$, that
 \be\label{wv-ca-c2-cov}
 (u_\epsilon, v_\epsilon)\rightarrow (u, v) \text{ in }  C_{\text{loc}}^{\delta, \frac{\delta}{2}}(\bar\Omega\times [t_0, \infty))\times C_{\text{loc}}^{2+\delta, 1+\frac{\delta}{2}}(\bar\Omega\times [t_0, \infty)).
 \ee
 Therefore, by \eqref{uvw-com-con-es}, for any $t>t_0$, we infer
 \begin{align*}
 \begin{split}
 \|u(t)-\bar u_0\|_{L^\infty}&\leq \|u_\epsilon(t)-u(t)\|_{L^\infty}+\|u_\epsilon(t)-\bar u_0\|_{L^\infty}\\
 &\leq \|u_\epsilon(t)-u(t)\|_{L^\infty}+C_1e^{-\mu (t-t_0)},
 \end{split}
 \end{align*}
 which, upon sending $\epsilon=\epsilon_i\searrow 0$, along with  \eqref{wv-ca-c2-cov} shows
 \be\label{weak-u-con}
 \|u(t)-\bar u_0\|_{L^\infty}\leq C_1e^{-\mu (t-t_0)},  \ \   \forall t\geq t_0.
 \ee
 Similarly, noting
$$
 \|v(t)-\bar u_0\|_{W^{1,\infty}}\leq \|v_\epsilon(t)-v(t)\|_{W^{1,\infty}}+\|v_\epsilon(t)-\bar u_0\|_{W^{1,\infty}},
 $$
 we again use \eqref{uvw-com-con-es} and \eqref{wv-ca-c2-cov} to deduce that
 \be\label{weak-v-con}
 \|v(t)-\bar u_0\|_{L^\infty}\leq C_2e^{-\mu (t-t_0)},  \ \   \forall t\geq t_0,
 \ee
Next, by the compact Sobolev embedding of $W^{1,6}(\Omega)$ into $L^\infty(\Omega)$, for all $t\geq t_0$, we conclude from \eqref{uvw-com-con-es} that $w_\epsilon(t)\rightarrow w(t)$ in  $L^\infty(\Omega)$ as $\epsilon=\tilde \epsilon_i\searrow 0$ and then
\be\label{weak-w-con}
 \|w(t)\|_{L^\infty}\leq C_3e^{-\mu (t-t_0)},  \ \   \forall t\geq t_0.
 \ee
 The reflexivity of $L^6(\Omega)$ and \eqref{uvw-com-con-es} enable us to infer further that
 $\nabla  w_\epsilon(t) \rightharpoonup \nabla w(t)$ in  $L^6(\Omega)$ as $\epsilon=\hat \epsilon_i\searrow 0$. Then, for $i=1,2,\cdots, n$,  we can choose $\phi_i\in L^{6/5}$ with $\|\phi_i\|_{L^{6/5}}=1$ to satisfy
 $$
 \|w_{x_i}(t)\|_{L^6}=\int_\Omega w_{x_i}(t)\phi_i.
 $$
 Notice that
 $$
\int_\Omega w_{\epsilon, x_i}(t)\phi_i\leq \|w_{\epsilon,x_i}(t)\|_{L^6}\|\phi_i\|_{L^{6/5}}= \|w_{\epsilon,x_i}(t)\|_{L^6}.
$$
Taking the limit inferior as $\epsilon=\hat \epsilon_i\searrow 0$, the weak convergence allows us  to infer
$$
\|\nabla w(t)\|_{L^6}\leq \liminf_{\epsilon=\hat \epsilon_i\searrow 0} \|\nabla w_\epsilon (t)\|_{L^6},
$$
which along with \eqref{uvw-com-con-es} again  entails that
$$
\|\nabla w(t)\|_{L^6}\leq C_4e^{-\mu (t-t_0)}, \ \  \forall t\geq t_0.
$$
This in conjunction with  \eqref{weak-u-con}, \eqref{weak-v-con} and \eqref{weak-w-con} yields the desired exponential convergence  in  \eqref{3D-convergence2}.

 Finally, those  convergence stated  in \eqref{3D-convergence2} and \eqref{wv-ca-c2-cov} together imply the eventual regularities in \eqref{reg-weak2}.
\end{proof}

\section*{Acknowledgments}
This work was partially supported by NSF of China (Nos. 12071476  and 12371214).

\section*{Data Availability}
Data sharing is not applicable to this article as no datasets were generated or analyzed during
the current study.
\section*{Declarations}
The authors declare that they have no conﬂict of interest.

\end{document}